\documentclass{article}
\usepackage[utf8]{inputenc}
\usepackage{amsmath}
\usepackage{amsfonts}
\usepackage{amssymb}
\usepackage{theorem,pifont}
\usepackage{graphicx}
\usepackage{color}

\title{The ISS framework for time-delay systems: a survey}
\author{Antoine Chaillet, Iasson Karafyllis, Pierdomenico Pepe, and Yuan Wang}
\date{\today}

\newcommand\mR{\mathbb{R}}
\newcommand\mRp{\mathbb{R}_{\geq 0}}
\newcommand\mN{\mathbb{N}}
\newcommand\cK{\mathcal{K}}
\newcommand\cN{\mathcal{N}}
\newcommand\cKL{\mathcal{KL}}
\newcommand\cX{\mathcal{X}}
\newcommand\cU{\mathcal{U}}
\newcommand\esssup{\textrm{ess\,sup}}

\newtheorem{theorem}{Theorem}
\newtheorem{corollary}{Corollary}
\newtheorem{conj}{Conjecture}
\newtheorem{definition}{Definition}
\newtheorem{example}{Example}
\newtheorem{remark}{Remark}
\newtheorem{lemma}{Lemma}
\newtheorem{proposition}{Proposition}
\newtheorem{stass}{Standing Assumption}

\definecolor{darkpastelgreen}{rgb}{0.01, 0.75, 0.24}
\definecolor{darkorange}{rgb}{1.0, 0.55, 0.0}

\begin{document}
\sloppy
\maketitle

\textbf{Abstract.} At the occasion of Eduardo D. Sontag's $70^{\textrm{th}}$ birthday, we provide here an overview of the tools available to study input-to-state stability (ISS) and related notions for time-delay systems. After a hopefully pedagogical presentation of the main differences with respect to the finite-dimensional theory, we review basic stability concepts for input-free time-delay systems, as well as instruments to guarantee them in practice, including the Lyapunov-Krasosvkii, Lyapunov-Razumikhin, and Halanay approaches. We then consider the influence of inputs through the notions of ISS, integral ISS, and input-to-output stability and provide both Lyapunov-like and solutions-based characterizations of these properties. We also show how these notions can be helpful for the stability analysis of interconnected systems, whether in cascade or in feedback form. We finally provide a list of questions which remain open until now.

\tableofcontents

\section{Introduction}

In the late 1980's, Eduardo D. Sontag introduced a concept that profoundly changed the way to approach stability and robustness of dynamical systems: the input-to-state stability (ISS) property \cite{SON1}. This property essentially imposes that the norm of any solution is bounded by a decaying term of the initial state's norm plus a term involving the amplitude of the applied input. Since this latter term is continuous and zero at zero, ISS ensures solutions' boundedness in response to any bounded input, small steady-state error for inputs of sufficiently small magnitude, and a vanishing state in response to any vanishing input.

One of the key reasons of the success of the ISS property lies in its Lyapunov characterization. It was shown in \cite{SONWANSCL} that ISS is equivalent to the existence of a Lyapunov function candidate whose derivative along the system's solution is upper bounded by a negative unbounded dissipation rate involving the state norm plus a continuous term involving the input norm. The simplicity of this characterization  and its resemblance to Lyapunov tools used for input-free systems not only provide an easy way to check ISS in practice, but was also at the basis of several further developments in terms of analysis of interconnected systems  and control design, as reviewed in  \cite{cetraro,Dashkovskiy:2011iv}. %\AC{We should cite A. Mironchenko's book once it is published.}

A decade after the birth of ISS, Eduardo D. Sontag introduced a weaker robustness property known as integral input-to-state stability (iISS, \cite{IISS}). Rather than assessing the input's influence through its amplitude, it takes into account the energy it feeds to the system. The Lyapunov characterization of iISS turns out to be similar to that of ISS, at the notable exception that the dissipation rate is no longer requested to be unbounded, but just positive definite \cite{ANG-bigIISS}. In particular, the corresponding Lyapunov function is no longer guaranteed to decay when the state is large, even when the applied input is of small magnitude. Nevertheless, iISS does ensure some robustness features, by guaranteeing that the state vanishes in response to any input with bounded energy (as measured through a specific nonlinear gain).

Both ISS and iISS imply that the origin of the system is globally asymptotically stable when the input is identically zero (0-GAS). In some applications, this turns out to be a too demanding requirement, either because the considered application is concerned with the behavior of particular state variables with no real interest for the other variables or because these extra variables simply do not have the requested stable behavior. A typical illustration is adaptive control, in which only the state variables are requested to be properly controlled whereas the parameter estimation does not need to be precise. To overcome this limitation, the concept of input-to-output stability (IOS) was introduced by Eduardo D. Sontag and the last author of the present survey \cite{SONWANIOS}. This property embeds all state variables of interest into an output and requests that only this output satisfies an ISS-like estimate. Accordingly, the characterization of this property does not require that the considered Lyapunov function be positive definite or radially unbounded in all state variables, but only in terms of the selected output \cite{SONWANIOS-LYA}.

In the ISS constellation, we may also mention the input/output-to-state stability property (IOSS), which relates the state norm to the magnitude of the system's input and output \cite{KRISONWAN}. In this setup, the output is no longer meant to be an error signal one wishes to attenuate, but rather a signal available for measurements from which one would like to reconstruct the full state. The IOSS property can thus be seen as a robust version of zero-state detectability \cite{SONSCLOSS}. The Lyapunov characterization of IOSS imposes a dissipation rate involving the state norm plus two positive terms: one in the input norm and the other in the output norm \cite{KRISONWAN}.

The whole ISS framework was originally built for finite-dimensional systems, meaning for systems ruled by an ordinary differential equation. Nevertheless, as recently reviewed in \cite{mironchenko2020input}, it has progressively been extended to the infinite dimensional case. The present survey focuses on a particular class of infinite-dimensional systems: time-delay systems. One motivation to study such systems is obvious and lies in the pervasive nature of delays in physical processes (transport phenomena, finite propagation velocity of information, non-instantaneous reaction, mechanical slack, sampled-data control, control over digital networks, etc.). The other motivation is of a more theoretical nature, as the mathematical peculiarities of time-delay systems allow to derive results that would be unreachable in a more general infinite-dimensional setting. These two motivations explain the number of monographs on the subject, including \cite{KRA59, Halanay:1966wl, Burton85, HALU93, KIM1999, MAHMOUD, ERNEUX, KAJIbook11, Gu:2012vm, KHARITONOV2013, KOMY13, Fridman:2014vo, MINI08, CORBRIAT2015, LIUFRIXIA2021}. 

In order to make this survey accessible to non-experts in time-delay systems, we start in Section \ref{pp:sec:1} by providing some mathematical background and intuition about the considered class of systems. We highlight the key difference with respect to ordinary differential equations, namely that the state is no longer a point of a finite-dimensional vector space, but rather a history segment (a function). We also recall basic results ensuring existence, uniqueness and regularity of solutions, and introduce two notions of forward completeness which happen to coincide in finite dimension. We finally discuss how to assess the evolution of a functional (meaning a function of the history segment) along the system's solutions.

In Section \ref{sec:no-input}, we review some stability notions of autonomous time-delay systems, thus first disregarding the influence of inputs. In particular, we cover both asymptotic and exponential stability and present tools to establish them in practice. These tools heavily rely on Lyapunov-Krasovkii functionals (LKF), which essentially play the same role for time-delay systems as Lyapunov functions in finite dimension. Beyond the fact that LKFs involve history segments, the main difference with Lyapunov functions stands in the way they may dissipate along solutions: either as function of the the current solution's norm (point-wise dissipation), in terms of the LKF itself (LKF-wise dissipation), or in terms of the norm of the history segment (history-wise dissipation). They may also be sandwiched between functions of the norm of the history segment (coercive LKF), but their practical use often makes it convenient to allow for more general lower-bounds. All these possibilities are thoroughly discussed. We also provide some results about output stability analysis, in which a stable behavior is expected only for a part of the state variables (or, more generally, for a specific output).

Section \ref{iasson:sec-iss} is devoted to the ISS property. We show that Eduardo D. Sontag's original concept smoothly translates to time-delay systems. We provide classes of systems (namely, linear or globally Lipschitz ones) for which ISS can easily be derived based on internal stability properties. We also provide some LKF characterizations of ISS, as well as sufficient conditions for ISS by using Razumikhin's or Halanay's approaches, and give sufficient conditions under which a point-wise dissipation is enough to guarantee ISS. We finally present some extensions of the solutions-based characterizations of ISS available in finite dimension \cite{SONWANTAC}.

Section \ref{sec-iISS} has the same objectives than Section \ref{iasson:sec-iss}, but for the iISS property. We show in particular that iISS can be established no matter how it dissipates along solutions (point-wise, LKF-wise or history-wise) and whether or not it is coercive. Solutions-based characterizations of iISS are also given.

In Section \ref{yw: sec-ios}, we review some results pertaining to input-to-output stability properties. In particular, we provide an LKF characterization of the IOS property, by  assuming either that the state norm is upper-bounded by a function of the output norm modulo a constant, or that the bounded input-bounded state property holds. We also say a word about the integral extension of IOS.

In Section \ref{sec:interconnection}, we make use of the ISS formalism to study stability of interconnected systems. Although the literature on feedback-interconnected systems is quite vast, we focus here on the interconnection of two systems only and provide both LKF-based and solutions-based small gain conditions to preserve ISS. For cascade interconnections, we show that ISS is preserved with no additional requirement. We also provide growth rate constraints to ensure that the cascade of two iISS systems is itself iISS.

We hope that this survey demonstrates that the ISS framework is now a mature subject for time-delay systems, with plenty of tools already available to study analysis and robustness. Nevertheless, some important questions remain open in this field of research: in Section \ref{sec:open}, we list some of them, explain their relevance, and sometimes suggest research directions to solve them.

Although the ISS framework was at the basis of several control designs for time-delay systems, we have decided not to present them here. In particular, for control methodologies that account for delays in the input/output channel, the reader is invited to consult the monographs \cite{ZHONG2006, RICOCAMACHO2007, ZHANGXIE2007, KRSTICBOOK2009, BEKIARISKRSTIC2013, KARKRS2017,    ZHUQIMACHEN18,LIUFRIXIA2021}, as well as the tutorial chapter on stabilization of delay systems \cite{KARMALMAZPEP2016}. 
%as well as the papers \cite{mazenc2006backstepping,mazenc2008further}
%\PP {there are many papers for control with input/output delays in both the linear and the nonlinear case.... we should pay some attention here if we cite just a few ones and not others.} \AC{You are right, it is just that I feel unfair not to cite Mazenc's work\ldots} \PP{I Would add also recent book \cite{LIUFRIXIA2021} and the tutorial on stabilization of delay systems \cite{KARMALMAZPEP2016}}

\section{Background}\label{pp:sec:1}

\subsection{Notation}\label{sec-Notation}

In this survey, $\mathbb R$ stands for the set of real numbers,  whereas $\mathbb N$ stands for the set of non-negative integers. Given $a\in\mR$, $\mR_{\geq a}:=\{x\in\mR\,:\, x\geq a\}$, and similarly for $\mR_{\leq a}$ and $\mathbb N_{\geq a}$. $\overline {\mathbb R}$ stands for the extended real line $[-\infty,+\infty]$. The symbol $\circ$ denotes the composition of functions. The symbol $^\top$ denotes the transpose of a matrix. $I$ denotes the identity matrix of appropriate dimension. Given $n\in\mathbb N_{\geq 1}$ and a continuously differentiable $V:\mR^n\to\mRp$, $\nabla V:\mR^n\to\mR^n$ denotes its gradient. The symbol $\vert\cdot \vert$ stands for the Euclidean norm of a real vector, or the induced Euclidean norm of a matrix. Given a non-empty compact set $\Omega\subset\mR^n$, ${\rm dist}(x, \Omega)$ denotes the distance between $x\in\mR^n$ and $\Omega$, namely: ${\rm dist}(x,\Omega):=\inf_{z\in\Omega}|x-z|$. Given a non-empty (possibly unbounded) interval $\mathcal I\subset\mR$ and a Lebesgue  measurable signal $u:\mathcal I\to\mR^m$, $\sup_{t\in \mathcal I} |u(t)|\in[0,+\infty]$ denotes the essential supremum of $u$, that is
\begin{align*}
	\sup_{t\in \mathcal I}\vert u(t)
	\vert:=\inf\left\{ \overline u\geq 0\,:\, \lambda\left(\{t\in \mathcal I:\vert
	u(t)\vert >\overline u      \} \right)=0 \right\},
\end{align*}
where $\lambda$ denotes the Lebesgue measure. We also let $\Vert
	u\Vert:=\sup_{t\in \mathcal I}\vert u(t) \vert$. The signal $u$ is said to be locally essentially bounded if it is essentially bounded on any bounded subset of $\mathcal I$, and essentially bounded if $\|u\|<+\infty$.
	
Given $\Delta\ge 0$, $\mathcal{X}$ denotes the space of the continuous functions mapping the interval $[-\Delta,0]$ into $\mathbb R$. Given $n\in\mN_{\geq 1}$, $T\in (0,+\infty]$, $t\in[0,T)$, and a continuous function $x:[-\Delta,T)\to \mathbb R^n$, $x_{t}\in\mathcal X^n$ denotes the history segment at time $t$ and is defined as
	$x_{t}(\tau):=x(t+\tau)$ for all $\tau\in[-\Delta,0]$. The set $\mathcal W^n\subset \mathcal X^n$ denotes the space of absolutely continuous functions mapping the
	interval $[-\Delta,0]$ into $\mathbb R^{n}$ with essentially bounded derivative. Given any $k\in\mN_{\geq 1}$, $C^k([-\Delta,0];\mR^n)\subset\mathcal W^n$ denotes the set of all functions $\phi\in\mathcal X^n$ with continuous derivatives up to order $k$. A function $z:\mathbb R_{\ge 0}\to \mathbb R^n$ is said to be piece-wise continuous if, given any $t_2>t_1\geq 0$, it is continuous on $[t_1,t_2]$, except possibly at a finite number of points. The symbol $\mathcal U$ denotes the set of the Lebesgue measurable and locally essentially bounded functions $u:\mRp\rightarrow \mathbb R$. Given $m\in\mathbb N_{\geq 1}$ and $u\in\mathcal U^m$, $u_{\mathcal I}$ denotes the restriction of $u$ to $\mathcal I$, namely $u_{\mathcal I}:\mathcal I\to\mR^m$ is defined as $u_{\mathcal I}(t):=u(t)$ for all $t\in \mathcal I$.

A function $\alpha:\mRp\to\mRp$ is said to be of class $\mathcal N$ if it is continuous, non-decreasing, and satisfies $\alpha(0)=0$. It is said to be of class $\mathcal P$ if it is continuous and satisfies $\alpha(0)=0$ and $\alpha(s)>0$ for all $s>0$. It is said to be of class $\cK$ if $\alpha\in\mathcal P$ and it is increasing. It it said to be of class $\cK_\infty$ if $\alpha\in\cK$ and $\lim_{s\to+\infty}\alpha(s)=+\infty$. A function $\ell:\mRp\to\mRp$ is said to be of class $\mathcal L$ if it is continuous, non-increasing, and satisfies $\lim_{s\to+\infty} \ell(s)=0$. A function $\beta:\mathbb R_{\ge 0}\times \mathbb R_{\ge 0}\to \mathbb R_{\ge 0}$
	is said to be of class $\mathcal {KL}$ if, for each fixed $t\ge 0$, $\beta(\cdot,t)\in\cK$ and, for each fixed $s\ge 0$, $\beta (s,\cdot)\in\mathcal L$.

\newpage
\subsection{List of acronyms}

For ease of reference, we list below all the acronyms used throughout this survey.

\hspace{-2cm}
\begin{tabular}{llr}
%	\item RFDE: retarded functional differential equation
0-GAS & globally asymptotically stable in the absence of an input & Section \ref{sec:ISS:def} \\
0-GES & globally exponentially stable in the absence of an input & Section \ref{sec:ISS:def}\\
AG & asymptotic gain & Equation \eqref{eq-antoine-AG}\\
AS & asymptotic stability or asymptotically stable & Definition \ref{pp:defstabetc}\\
BECS & bounded energy-converging state & Equation \ref{eq_antoine_BECS}\\
BIBS & bounded input-bounded state & Equation \eqref{eq-antoine-BIBS}\\
CICS & converging input-converging state & Equation \eqref{eq-antoine-CICS}\\
ES & exponential stability or exponentially stable & Definition \ref{pp:defstabetc}\\
FC & forward completeness or forward complete & Definition \ref{def_FC}\\
GAS & global asymptotic stability or globally asymptotically stable & Definition \ref{pp:defstabetc}\\
GAOS & global asymptotic output stability or globally asymptotically output stable & Definition \ref{yw:def-os}\\
GES & global exponential stability or globally exponentially stable & Definition \ref{pp:defstabetc}\\
GOS& global output stability or globally output stable & Equation \eqref{eq-antoineOGS}\\
iIOS & integral input-to-output stability or integral input-to-output stable & Definition \ref{yw:def-iios}\\
iISS & integral input-to-state stability or integral input-to-state stable & Definition \ref{def_antoine_iISS}\\
IOS & input-to-output stability or input-to-output sable & Definition \ref{yw:def-ios}\\
ISS & input-to-state stability or input-to-state stable & Definition \ref{def_Iasson_ISS}\\
LISS & local input-to-state stability or locally input-to-state stable & Definition \ref{def_Iasson_ISS}\\
LKF & Lyapunov-Krasovskii functional candidate & Definition \ref{def-LKF}\\
%OAtt & output attractivity & Equation \eqref{eq-antoine-OAtt}\\
OL-iIOS & output-Lagrange iIOS & Definition \ref{yw:def-iios}\\
OL-IOS & output-Lagrange IOS & Definition \ref{yw:def-ios}\\
RFC & robust forward completeness or robustly forward complete & Definition \ref{def_RFC}\\
SI-iIOS & state-independent iIOS & Definition \ref{yw:def-iios}\\
SI-IOS & state-independent IOS & Definition \ref{yw:def-ios}\\
UAG & uniform asymptotic gain & Equation \eqref{eq-antoine-UAG}\\
UBEBS & uniform bounded energy-bounded state & Definition \ref{def_antoine_UBEBS}\\
UGOS & uniform global output stability or uniformly globally output stable & Equation \eqref{eq-antoine-BIBO}\\
UGS & uniform global stability or uniformly globally stable& Definition \ref{yw:UGS}\\
%UBIBS & uniform bounded input-bounded state & Definition %\ref{yw:bibs}\\
ULIM & uniform limit & Equation \eqref{eq_antoine_19}\\
ULS & uniform local stability or uniformly locally stable &  Equation \eqref{eq-antoine-ULS}
\end{tabular}

\subsection{Considered class of systems}

%
% \AC{I feel that Section 2 lacks references.}
%

The class of systems studied in this survey  is that of nonlinear time-delay systems, namely functional differential equations of the type
\begin{align}\label{pp:RFDE}
\dot x(t)=f(x_t,u(t)),
\end{align}
 where $t\geq 0$ is the time variable,  $f$ is a continuous map from $\mathcal {X}^n\times \mathbb R^m$ to $\mathbb R^n$, $x(t)\in \mathbb R^n$ is the internal variable, $\dot x(t)$ is the right-hand derivative\footnote{which coincides with the derivative of $t\mapsto x(t)$ at each time $t$ where $x$ is differentiable.} of $x(t)$ with respect to $t$, $u\in\mathcal U^m$ is the input, and $x_t:[-\Delta,0]\to \mathbb R^n$ is the solution's history, as defined in Section \ref{sec-Notation}, where
$\Delta\geq 0$ denotes the maximum time delay involved. 

\begin{figure}[h!]
    \centering
    \includegraphics[width=0.8\textwidth]{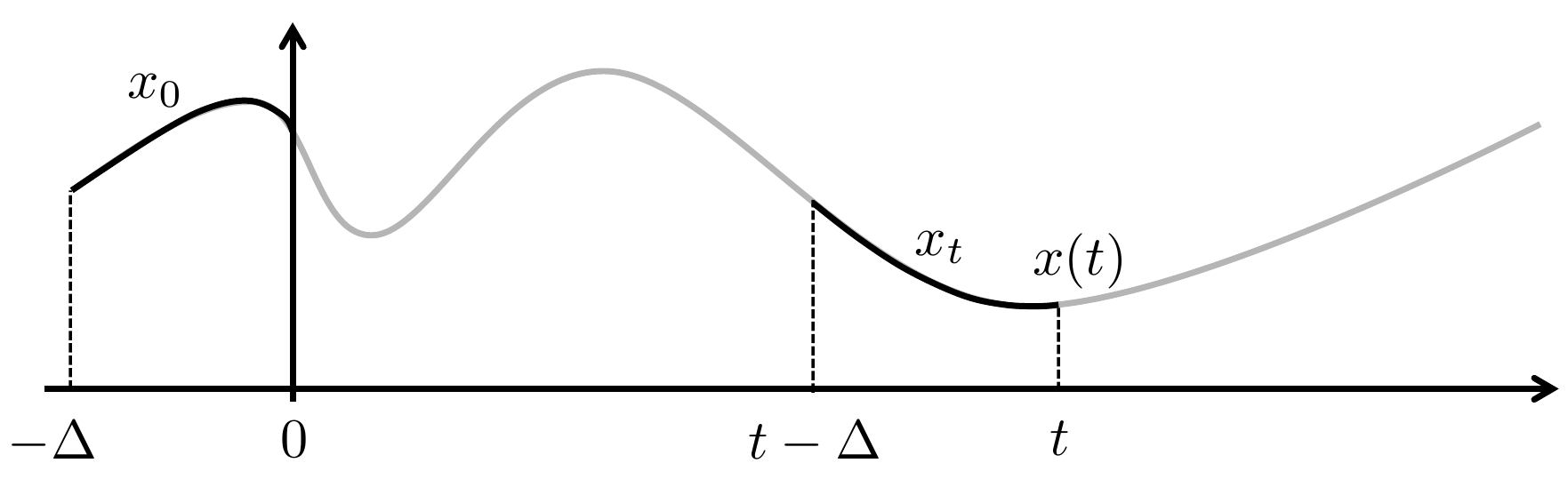}
    \caption{Time evolution of the solution of (\ref{pp:RFDE}).}
    \label{fig-1}
\end{figure}

The state space of the system (\ref{pp:RFDE}) is $\mathcal X^n$: the status of the system at some time $t\geq 0$ is not captured by current value of the solution $x(t)\in\mR^n$ only, but rather from the history segment $x_t\in\mathcal X^n$ as depicted by Figure \ref{fig-1}. This constitutes the main difference with respect to finite-dimensional systems: the state is no longer a point in $\mR^n$, but rather a function in $\mathcal X^n$, thus explaining the infinite dimension of \eqref{pp:RFDE}. 

The class of systems covered by \eqref{pp:RFDE} is rather wide. It encompasses single or multiple discrete delays (possibly non-commensurate) as well as distributed time delays. It also covers systems whose dynamics depend on the maximum or minimal value of the internal variable over some given time window.

We stress that \eqref{pp:RFDE} is time-invariant, in the sense that the vector field $f$ does not evolve with time. Time-varying systems are not considered in this survey, but extensive work on uniform and nonuniform stability properties for time-varying systems can be found in \cite{KAJIbook11,KAPEJI08,Karafyllis:2008hc} and references therein.

For the system \eqref{pp:RFDE}, the notion of solution is defined as follows.

\begin{definition}[Solution]
Given an initial state $x_0\in \mathcal X^n$ and an input signal $u\in \mathcal U^m$, a \emph{solution} of \eqref{pp:RFDE} denotes any function $x:[0,t_{\max})\to \mathbb R^n$, with  $t_{\max}\in(0,+\infty]$, which is locally absolutely continuous and satisfies \eqref{pp:RFDE} almost everywhere in $[0,t_{\max})$ or, equivalently, satisfies the integral equation
\begin{equation*}
    x(t)=x_0(0)+\int_{0}^tf(x_{\tau},u(\tau))d\tau,\quad \forall t\in [0,t_{\max}). 
\end{equation*}
\end{definition}

\subsection{Existence, uniqueness, regularity of solutions}

Some regularity conditions need to be imposed on the vector field $f$ in order to ensure existence and uniqueness of solutions. Due to its infinite-dimensional nature, several notions of Lipschitz continuity can be envisioned.

\begin{definition}[Lipschitz continuity] \label{pp:deflipschitzrfde}
	A map $f:\mathcal {X}^n\times \mathbb R^m\to \mathbb R^n$ is said to be:	
	\begin{itemize}
	\item \emph{locally Lipschitz} if, for any $\phi\in\mathcal X^n$ and any $v\in \mathbb R^m$, there exist $\delta,L>0$ such that, for all $\phi_1,\phi_2\in\mathcal X^n$ and all $v_1,v_2\in \mR^m$ satisfying $\|\phi_i-\phi\|\leq \delta$ and $|v_i-v|\leq \delta$, $i\in\{1,2\}$, 
		\begin {align}\label{eq_antoine_113}
		\vert f(\phi_1, v_1)-f(\phi_2,v_2)\vert\le L(\Vert \phi_1-\phi_2\Vert+\vert v_1-v_2\vert)
		\end {align}
	\item \emph{Lipschitz on compact sets} if, for any compact set $\mathcal C\subset \mathcal X^n \times \mathbb R^m$, there exists $L>0$ such that \eqref{eq_antoine_113} holds for all $(\phi_1,v_1), (\phi_2,v_2)\in \mathcal C$
		%\begin {align*}
		%\vert f(\phi_1, v_1)-f(\phi_2,v_2)\vert\le L(\Vert \phi_1-\phi_2\Vert+\vert v_1-v_2\vert)
		%\end {align*}

		\item \emph{Lipschitz on bounded sets} if, for any $r> 0$, there
		exists $L> 0$ such that \eqref{eq_antoine_113} holds for all $\phi_1,\phi_2\in \mathcal X^n$ with $\|\phi_1\|\leq r$ and $\|\phi_2\|\leq r$ and all $v_1,v_2\in \mR^m$ with $|v_1|\leq r$ and $|v_2|\leq r$
		%\begin {align*}
		%\vert f(\phi_1, v_1)-f(\phi_2,v_2)\vert\le L(\Vert
		%\phi_1-\phi_2\Vert+\vert v_1-v_2\vert)
		%\end {align*}
		
		\item \emph{globally Lipschitz} if there exists $L> 0$ such that \eqref{eq_antoine_113} holds for all $\phi_1,\phi_2\in \mathcal {X}^n$ and all $v_1,v_2\in\mathbb R^m$.
		%\begin {align*}
		%\vert f(\phi_1, v_1)-f(\phi_2,v_2)\vert\le L(\Vert \phi_1-\phi_2\Vert+\vert v_1-v_2\vert).
		%\end {align*}
	\end{itemize}
\end{definition}

Analogous definitions of Lipschitz properties hold for a functional $V:\mathcal X^n\to \mathbb R_{\ge 0}$, by simply removing the requirements concerning the second entry.

In finite dimension, local Lipschitz, Lipschitz on compact sets, and Lipschitz on bounded sets are all equivalent properties.  In the infinite dimensional case, Lipschitz on bounded sets clearly implies Lipschitz on compact sets. The following recent result, established in \cite[Theorem 3.1]{XULIUFEN20} in a wider infinite-dimensional context, shows that local Lipschitz also turns out to be equivalent to Lipschitz on compact sets for the considered class of functions. Nevertheless, as shown in \cite[Theorem 3.2]{XULIUFEN20}, local Lipschitz does not imply the Lipschitz on bounded sets property.
 
\begin{theorem}[Link between Lipschitz conditions]\label{thm-Lipschitz}
 A map $f:\mathcal {X}^n\times \mathbb R^m\to \mathbb R^n$ is \emph{Lipschitz on 
 compact sets} if and only if it is \emph{locally Lipschitz}.
 \end{theorem}

Throughout this paper, the following properties are assumed on the dynamics.

\begin{stass}\label{stand_ass_1}
The vector field $f$ involved in \eqref{pp:RFDE} is Lipschitz on bounded sets and satisfies $f(0,0)=0$,
\end{stass}

This assumption ensures that $x(\cdot)\equiv 0$ is the solution
corresponding to zero initial state and zero input (sometimes referred to as the trivial solution). It also ensures existence, uniqueness, continuous dependency and maximum continuation? of solutions, as established in \cite[Theorems 2.1, 2.2, 2.3, 3.1, 3.2, Section 2.6]{HALU93} and \cite[Theorems 2.1, 2.2]{KOMY13}. 

\begin{theorem}[Existence, uniqueness, continuity, maximum continuation of solutions]\label{pp:thm:exist}
The following results hold:
\begin{itemize}
\item for any initial state $x_0\in \mathcal {X}^n$ and
any input $u\in\mathcal U^m$, (\ref{pp:RFDE}) admits a unique locally absolutely
continuous solution $x(\cdot,x_0,u)$ on a maximal time interval $[0,t_{\max})$ with
$t_{\max}\in(0,+\infty]$
\item if $t_{\max}<+\infty$, then the solution is unbounded
in $[0,t_{\max})$
\item for any $T\in (0,t_{\max})$ and any $\varepsilon>0$ there exist $\delta> 0$ such that,
for any $z_0 \in \mathcal X^n$ with $\|x_0-z_0\|\leq \delta$ and any $v\in\mathcal U^m$ with
$\esssup_{t\in [0,T]}\vert v(t)-u(t)\vert\le \delta$,
the solution $x(\cdot,z_0,v)$ of (\ref{pp:RFDE}) exists on $[0,T]$ and satisfies
\begin{align*}
\vert x(t,z_0,v)-x(t,x_0,u)\vert \leq \varepsilon, \qquad \forall t\in [0,T].
\end{align*}
%\item if the vector field $f$ in (\ref{pp:RFDE}) is \emph{globally Lipschitz}, then $t_{\max}=+\infty$. 
\end{itemize}
\end{theorem}

In most of the survey, we will write $x(t,x_0,u)\in\mR^n$ the maximal solution of \eqref{pp:RFDE} at time $t\geq 0$ corresponding to initial state $x_0\in\mathcal X^n$ and input
$u\in\mathcal U^m$. Similarly, $x_t(x_0,u)\in\mathcal X^n$ will denote the corresponding history segment, expressed in $\mathcal X^n$. The following statement, given in \cite[Lemma 2.1]{HALU93} and \cite[Lemma 4]{Pepe2007}, ensures some regularity of $x_t$ with respect to $t$. 

\begin{theorem}[Regularity with respect to time] \label{pp:lemmahalepeperfde} For any initial state $x_0\in \mathcal {X}^n$ and any input $u\in \mathcal U^m$, let $t_{\max}\in(0,+\infty]$ denote the maximal time of existence of the corresponding solution of (\ref{pp:RFDE}). Then %the solution $x(t,x_0,u)$, $t\in [0,b)$, of the RFDE (\ref{pp:RFDE}) is such that 
the function $t\mapsto x_t(x_0,u)$ is continuous on $[0,t_{\max})$. 
Moreover, if $x_0\in \mathcal W^n$ (and, a fortiori, if $x_0\in C^1([-\Delta,0];\mR^n)$), then $t\mapsto x_t(x_0,u)$ is locally absolutely continuous on $[0,t_{\max})$.
\end{theorem}

The above result states in particular that local absolute continuity of $t\mapsto x_t$ can be ensured provided that we restrict the considered class of initial states to absolutely continuous signals. As will become clearer later, this turns out particularly useful when studying the evolution of a functional $V$ along the system's solutions. The following result shows that restricting the class of considered initial states to $\mathcal W^n$ actually comes with no loss of generality for all the stability and robustness properties covered by this survey.

\begin{lemma}[Restricting the class of initial states]\label{lem-dense}
Let $\rho:\mRp\times \mathcal X^n\times \mathcal U^m\to\mR$ be such that, for each $t\geq 0$ and each $u\in\mathcal U^m$, the map $x_0\mapsto \rho(t,x_0,u)$ is continuous. Assume that there exists a continuous function $W:\mR^n\to\mRp$ such that, for all $x_0\in C^1([-\Delta,0];\mR^n)$ and all $u\in\mathcal U^m$,
\begin{align}\label{eq_antoine_109}
    W(x(t,x_0,u))\leq \rho(t,x_0,u)
\end{align}
for all $t$ in the maximal interval of existence of $x(\cdot,x_0,u)$. Then \eqref{eq_antoine_109} also holds for all $x_0\in\mathcal X^n$ and all $u\in\mathcal U^m$ over the maximal interval of existence of $x(\cdot,x_0,u)$.
\end{lemma}

The proof of this result consists in invoking the density of $C^1([-\Delta,0];\mR^n)$ in $\mathcal X^n$ and the continuity of solutions with respect to the initial state (Theorem \ref{pp:thm:exist}): see \cite[Lemma 2.6]{Karafyllis:2008hc} and \cite[Proposition 3]{Pepe2007}.

\subsection{Forward completeness}

While Theorem \ref{pp:thm:exist} ensures existence of solutions over some time interval $[0,t_{\max})$, it does not guarantee that $t_{\max}=+\infty$, as imposed by the following property.

\begin{definition}[FC]\label{def_FC}
The system \eqref{pp:RFDE} is said to be \emph{forward complete (FC)} if, for any initial state $x_0\in \mathcal X^n$ and any input $u\in \mathcal U^m$, the corresponding solution exists on $\mathbb R_{\ge 0}$.
\end{definition}

A seemingly more conservative property that has serious consequences for time-delay systems is the following RFC property \cite[Definition 2.1]{KAJIbook11}, sometimes also known as the \emph{Bounded Reachability Sets (BRS)} property \cite[Definition 4]{mironchenko2017characterizations}.

\begin{definition}[RFC]\label{def_RFC} The system \eqref{pp:RFDE} is said to be \emph{robustly forward complete (RFC)} if it is FC and, for every $T,r>0$,
\begin{equation*} %\label{pp:Iasson__14_}
\sup \left\{\, \left\| x_{t} (\phi ,u)\right\| \, :\, \phi \in \mathcal X^n,\ \left\| \phi \right\| \le r, \ t\in [0,T],\ u\in \mathcal U^m,\ \|u\|\leq r \right\}<+\infty.
\end{equation*}
\end{definition}

In other words, the RFC property imposes that, starting from any bounded set of initial states and considering inputs in any given bounded set, solutions can only evolve in a bounded region over any finite time interval. In finite dimension, FC and RFC turn out to be equivalent properties \cite[Corollary III.4]{SONWANTAC}. It is not known at this stage whether the same holds for time-delay systems: see Section \ref{sec:open:UGAS} for further discussions on this important question.

An easy (though conservative) way to ensure RFC is to request that $f$ is globally Lipschitz, but more interesting conditions for RFC will be presented in Section \ref{sec-RFC-cond}.

\begin{proposition}[Global Lipschitz $\Rightarrow$ RFC]\label{prop-RFC-Lipschitz} If the map $f$ is globally Lipschitz in $\mathcal X^n\times \mathbb R^m$, then \eqref{pp:RFDE} is RFC.
\end{proposition}

\subsection{Functional derivatives along solutions}

In view of Theorem \ref{pp:lemmahalepeperfde}, it can easily be checked that, given any functional $V:{\cal X}^n\to  \mRp$, Lipschitz on bounded sets, any initial
state $x_0\in \mathcal W^n$ and any input $u\in\mathcal U^m$, the map $t\mapsto V(x_t)$ is locally absolutely continuous over the maximal interval of existence of the solution $x(\cdot)$ (see \cite[Theorem 5, Remark 6]{Pepe2007} or \cite[Lemma 2.5]{Karafyllis:2008hc}). 
This regularity is needed, for instance, in order to guarantee that $t\mapsto V(x_t)$ is non-increasing when its derivative is almost everywhere non-positive. %Recall the Cantor function which is continuous (but not absolutely continuous), non-negative and monotonic in $[0,1]$, it's upper-right Dini derivative is almost everywhere equal to $0$, and the value at $1$ is greater than the value at $0$.  

The following functional derivative will be extensively used in this survey. It allows to compute the derivative of a functional along a system's solution without needing to consider the solution of the system, not even formally \cite{driver1962existence}.

\begin{definition}[Driver's derivative]\label{pp:defdipiuv} Let $V:\mathcal {X}^n\to \mathbb R_{\ge 0}$ be a continuous functional. Its \emph{Driver's derivative} $D^+V:\mathcal X^n\times \mR^n\to\overline\mR$ is defined, for all $\phi \in \mathcal X^n$ and all $w\in \mathbb R^n$, as
\begin{equation*}%\label{pp:definitionvpuntorfde}
D^+ V(\phi,w):=
\limsup_{h\to 0^+} \frac {V(\phi^{}_{h,w})-V(\phi)}{h},
\end{equation*}
where, for each $h\in [0,\Delta)$, $\phi^{}_{h,w} \in \mathcal {X}^n$ is given by
\begin{equation*}
\phi^{}_{h,w} (s) :=  \left \{ \begin{array}{cl}  \phi(s+h), & \textrm{ if }s\in [-\Delta,-h),    \\
\phi(0)+w(h+s), & \textrm{ if } s\in [-h,0].
\end{array} \right .
\end{equation*}
\end{definition}

In practice, we often want to study Driver's derivative in the direction imposed by the vector field of \eqref{pp:RFDE}, in which case we simply let $w=f(\phi,v)$, with $\phi\in \mathcal X^n$ and $v\in \mathbb R^m$, in the above definition.

\begin{remark}\label{pp:remarkdriverbyallfrozen}
Notice that Driver's derivative can be computed as the upper-right Dini derivative along a dynamical system with constant right-hand side. More precisely,
\begin{equation*}
D^+ V(\phi,w)=
\limsup_{h\to 0^+} \frac{V(z_h)-V(\phi)}{h} ,\end{equation*}
where, for $h\in [0,\Delta)$, $z_{h} \in \mathcal {X}^n$ is the solution of the particular functional differential equation
\begin{align*}
\dot z(t)&=w\\
z_0&=\phi
\end{align*}
So, for instance, given the solution $x_t$ of \eqref{pp:RFDE}, $t\in [0,t_{\max})$,  corresponding to an initial state $x_0\in \mathcal X^n$ and an input $u\in \mathcal U^m$, we have, for $t\in [0,t_{\max})$, whenever $u(t)$ is defined (i.e., almost everywhere in $[0,t_{\max}))$, 
\begin{equation}
D^+ V(x_t,f(x_t,u(t)))=
\limsup_{h\to 0^+} \frac {1}{h} \left ( V(z_h)-V(x_t) \right
),\end{equation}
where, for $h\in [0,\Delta)$, $z_{h} \in \mathcal {X}^n$ is the solution of the particular system (with ``frozen'' right-hand side) 
\begin{align*}
\dot z(\tau)&=f(x_t,u(t)), \qquad  \tau\ge 0\\
z_0&=x_t.
\end{align*}
\end{remark}

It turns out that Driver's derivative provides useful information of the upper-right Dini derivative on the function $t\mapsto V(x_t)$, as explained next (see \cite{driver1962existence}, \cite[Theorem 2]{pepe2007liapunov}).

\begin{theorem}[Link between Dini and Driver] \label{pp:lemmadriverrfde}
Let $V:\mathcal {X}^n\to \mathbb R_{\ge 0}$ be a locally Lipschitz functional.
Given $x_0\in\mathcal X^n$ and $u\in\mathcal U^m$, let $x:[0,t_{\max})\to\mR^n$, $t_{\max}\in (0,+\infty]$, denote the corresponding maximal solution of (\ref{pp:RFDE}). Let
$
\nu:[0,t_{\max})\to \mRp$ be the function defined by $\nu(t):=V(x_t)$ and $\mathcal D^{+}\nu:[0,t_{\max})\to\overline\mR$ denote its upper right-hand Dini derivative, that is
\begin{equation}\label{pp:derivatawrfde}
\mathcal D^+\nu(t):=\limsup_{h\to 0^+} \frac {\nu(t+h)-\nu(t)}{h},\quad \forall t\in[0,t_{\max}).
\end{equation}
Then, for almost all $t\in[0,t_{\max})$, it holds that
\begin{equation}\label{pp:relazionesolitarfde}
\mathcal D^+\nu(t)=D^{+}V\big(x_t,f(x_t,u(t))\big).%=D_1^{+}V(x_t,u(t)).
\end{equation}
Moreover, \eqref{pp:relazionesolitarfde} holds for all $t\in[0,t_{\max})$ provided that the input $u\in\mathcal U^m$ is piece-wise continuous and right-continuous. 
\end{theorem}

In other words, by estimating a suitable upper bound for the Driver's derivative of $V$ along the vector field $f(\phi,v)$, with any $\phi \in \mathcal X^n$, and any $v\in \mathbb R^m$,  one obtains the same upper bound on the Dini derivative of $V(x_t)$ almost everywhere on $[0,t_{\max})$, by simply considering  $\phi=x_t$ and $v=u(t)$.

\begin{remark}[Absolutely continuous initial states] Recall that any locally absolutely continuous function is differentiable almost everywhere and that, wherever its derivative exists, it is equal to its upper-right Dini derivative. Hence, in view of Theorem \ref{pp:lemmahalepeperfde}, the upper-right Dini derivative in the above statement can be replaced by the standard derivative when restricting the class of initial states $x_0$ to $\mathcal W^n$ (and, even more so, to $C^1([-\Delta,0];\mR^n)$).
\end{remark}
%Then, for instance, the non-increasing property of the function $t\to V(x_t)$, $t\in [0,b)$, is guaranteed by it's Dini derivative to be almost everywhere non-positive, provided that the function $t\mapsto V(x_t)$ is locally absolutely continuous on $[0,b)$. 

The following example describes how to compute Driver's derivative of a particular functional and illustrates the fact that, along solutions, it coincides almost everywhere with the upper-right Dini derivative.

\begin{example}[Dini Vs. Driver]\label{exa1}
A standard functional used in the stability analysis of time-delay systems is
\begin{equation}\label{eq_antoine_23}
    V(\phi):=\phi(0)^{\top}P\phi(0)+\int_{-\Delta}^0e^{c\tau}\phi(\tau)^{\top}Q\phi(\tau)d\tau, \quad \forall \phi\in\mathcal X^n,
\end{equation}
where $P,Q\in \mathbb R^{n\times n}$ denote symmetric positive definite matrices and $c\geq 0$. This functional is Lipschitz on bounded sets, hence locally Lipschitz. Indeed, given any $r>0$, consider any $\phi_1,\phi_2\in\mathcal X^n$ satisfying $\|\phi_1\|\leq r$ and $\|\phi_2\|\leq r$. 
Then, by the mean value theorem, it holds that
\begin{align*}
    \left|V(\phi_1)-V(\phi_2)\right| =&\, \Big| \phi_1(0)^{\top}P\phi_1(0)-\phi_2(0)^{\top}P\phi_2(0) \\
   &\, +\int_{-\Delta}^0 e^{c\tau}\phi_1(\tau)^{\top}Q\phi_1(\tau)d\tau-\int_{-\Delta}^0 e^{c\tau}\phi_2(\tau)^{\top}Q\phi_2(\tau)d\tau\Big|\nonumber 
   \\
   \le&\, \max_{\lambda\in [0,1]}\left | 2\big(\lambda \phi_1(0)+(1-\lambda)\phi_2(0)\big)^{\top}P (\phi_1(0)-\phi_2(0))\right |\nonumber \\ &\,+\int_{-\Delta}^0
    \left | \phi_1(\tau)^{\top}Q\phi_1(\tau)-\phi_2(\tau)^{\top}Q\phi_2(\tau)\right |d\tau
        \nonumber \\ 
\le&\, \max_{\lambda\in [0,1]}\left | 2\big(\lambda \phi_1(0)+(1-\lambda)\phi_2(0)\big)^{\top}P (\phi_1(0)-\phi_2(0))\right |\nonumber \\ &\,+\int_{-\Delta}^0
    \max_{\lambda\in [0,1]}\left | 2\big(\lambda \phi_1(\tau)+(1-\lambda)\phi_2(\tau)\big)^{\top}Q (\phi_1(\tau)-\phi_2(\tau))\right
        |d\tau\nonumber \\ 
        \leq&\, 2r|P| |\phi_1(0)-\phi_2(0)|+2r\Delta|Q|\|\phi_1-\phi_2\|\\
    \leq&\, 2r\max\{|P|, \Delta|Q|\}\|\phi_1-\phi_2\|,
\end{align*}
thus showing that $V$ is Lipschitz on bounded sets. We show now that, for any $\phi\in \mathcal X^n$ and any $w\in \mathbb R^n$,
\begin{equation}\label{pp:exespressionedriver}
    D^+V(\phi,w)=2\phi(0)^{\top}Pw+\phi(0)^{\top}Q\phi(0)-\phi(-\Delta)^{\top}Q\phi(-\Delta)-c\int_{-\Delta}^0e^{c\tau}\phi(\tau)^{\top} Q\phi(\tau)d\tau.
\end{equation}

To that aim, we rely on Remark \ref{pp:remarkdriverbyallfrozen}. Using the Leibniz integral rule for the derivation under the integral sign applied to the function $h\mapsto \int_{-\Delta+h}^he^{c(\tau-h)}z(\tau)^{\top}Qz(\tau)d\tau$ and recalling that $z_0=\phi$, this remark ensures that
\begin{align*}
D^+V(\phi,w)=&\,\limsup_{h\to 0^+}\frac{1}{h}\left(V(z_h)-V(\phi)\right )\nonumber \\ =&\,\limsup_{h\to 0^+}\frac{1}{h}\Bigg(z_h(0)^{\top}Pz_h(0)+\int_{-\Delta}^0e^{c\tau}z_h(\tau)^{\top}Qz_h(\tau)d\tau\\
&\,-\phi(0)^{\top}P\phi(0)-\int_{-\Delta}^0e^{c\tau}\phi(\tau)^{\top}Q \phi(\tau)d\tau\Bigg)\nonumber \\ 
=&\,
\limsup_{h\to 0^+}\frac{1}{h}\Bigg(z(h)^{\top}Pz(h)+\int_{-\Delta}^0e^{c\tau}z(h+\tau)^{\top}Qz(h+\tau)d\tau\\
&\,-z(0)^{\top}Pz(0)-\int_{-\Delta}^0e^{c\tau}\phi(\tau)^{\top}Q \phi(\tau)d\tau\Bigg )\nonumber \\ 
=&\,
\limsup_{h\to 0^+}\frac{1}{h}\Bigg(z(h)^{\top}Pz(h)-z(0)^\top Pz(0)\\
&\,+\int_{-\Delta+h}^he^{c(s-h)}z(s)^{\top}Qz(s)ds-\int_{-\Delta}^0e^{c\tau}\phi(\tau)^{\top}Q \phi(\tau)d\tau\Bigg)\nonumber \\ 
=&\,
2z(0)^{\top}P\dot z(0)+z(0)^{\top}Qz(0)-e^{-c\Delta}z(-\Delta)Qz(-\Delta)\\
&\,-c\int_{-\Delta}^0e^{c\tau}z(\tau)^{\top}Qz(\tau)d\tau
\nonumber \\ 
=&\,2\phi(0)^{\top}Pw +\phi(0)^{\top}Q\phi(0)-e^{-c\Delta}\phi(-\Delta)^{\top}Q\phi(-\Delta)\\
&\,-c\int_{-\Delta}^0e^{c\tau}\phi(\tau)^{\top}Q\phi(\tau)d\tau,
\end{align*}
thus confirming \eqref{pp:exespressionedriver}. In particular, for $w=f(\phi,v)$, with $\phi\in \mathcal X^n$ and $v\in \mathbb R^m$, we obtain
\begin{align}
    D^+V(\phi,f(\phi,v))=&\,2\phi(0)^{\top}Pf(\phi,v)+\phi(0)^{\top}Q\phi(0)-\phi(-\Delta)^{\top}Q\phi(-\Delta)\nonumber\\
    &\,-c\int_{-\Delta}^0e^{c\tau}\phi(\tau)^{\top}Q\phi(\tau)d\tau.\label{pp:exespressionedriver2}
\end{align}
Given any initial state $x_0\in \mathcal X^n$ and any input $u\in \mathcal U^m$, let $x(\cdot)$ denote the corresponding solution of \eqref{pp:RFDE}, defined on $[0,t_{\max})$ with $t_{\max}\in(0,+\infty]$. Then it holds from \eqref{pp:exespressionedriver2} that, for all $t\in[0,t_{\max})$ where $u$ is defined (meaning at least almost everywhere in $[0,t_{\max})$),
\begin{align}
    D^+V(x_t,f(x_t,u(t)))=&\,2x_t(0)^{\top}Pf(x_t,u(t))+x_t(0)^{\top}Qx_t(0) \nonumber\\
    &\,-x_t(-\Delta)^{\top}Qx_t(-\Delta)
        -c\int_{-\Delta}^0e^{c\tau}x_t(\tau)^{\top}Qx_t(\tau)d\tau
        \nonumber \\  
        =&\,2x(t)^{\top}Pf(x_t,u(t))+x(t)^{\top}Qx(t)\label{pp:derivdrisol}\\
        &\,-x(t-\Delta)^{\top}Qx(t-\Delta)
                -c\int_{t-\Delta}^t e^{c(s-t)}x(s)^{\top}Qx(s)ds.\nonumber
\end{align}
\noindent 
On the other hand, let 
$
\nu:[0,t_{\max})\to \mRp$ be the function defined by $\nu(t):=V(x_t)$ for all $t\in [0,t_{\max})$ and let $\mathcal D^+\nu:[0,t_{\max})\to \overline {\mathbb R}$ denote its upper-right Dini derivative as defined in \eqref{pp:derivatawrfde}. Observing that, from \eqref{eq_antoine_23}, 
\begin{align*}
    \nu(t)&= x(t)^{\top}Px(t)+\int_{-\Delta}^0e^{c\tau}x(t+\tau)^{\top}Qx(t+\tau)d\tau\\
    & =  x(t)^{\top}Px(t)+\int_{t-\Delta}^te^{c(s-t)}x(s)^{\top}Qx(s)ds,
\end{align*}
and recalling that $x(\cdot)$ satisfies \eqref{pp:RFDE} almost everywhere on $[0,t_{\max})$, we obtain, by the Leibniz integral rule for the derivation under the integral sign, that, for almost all $t\in[0,t_{\max})$,
\begin{align}
\mathcal D^+\nu(t)=&\,2x(t)^{\top}P\mathcal D^+x(t)+x(t)^{\top}Qx(t)-x(t-\Delta)^{\top}Qx(t-\Delta)\nonumber \\
&\,-c\int_{t-\Delta}^te^{c(\tau-t)}x(\tau)^{\top}Qx(\tau)d\tau\nonumber \\
=&\,2x(t)^{\top}P\dot x(t)+x(t)^{\top}Qx(t)-x(t-\Delta)^{\top}Qx(t-\Delta)\nonumber \\
&\,-c\int_{t-\Delta}^te^{c(\tau-t)}x(\tau)^{\top}Qx(\tau)d\tau\nonumber \\
=&\,2x(t)^{\top}Pf(x_t,u(t)) +x(t)^{\top}Qx(t)-x(t-\Delta)^{\top}Qx(t-\Delta))\nonumber \\
&\,-c\int_{t-\Delta}^te^{c(\tau-t)}x(\tau)^{\top}Qx(\tau)d\tau.\label{pp:comp_ex_derivatadini}
\end{align}
Comparing \eqref{pp:derivdrisol} and  \eqref{pp:comp_ex_derivatadini}, we confirm the expectations from Theorem \ref{pp:lemmadriverrfde}, namely that, for almost all $t\in [0,t_{\max})$,
\begin{equation}\label{pp:ex_dini_eq_driver}
   \mathcal D^+\nu(t)= D^+V(x_t,f(x_t,u(t))).
\end{equation}
If the equation \eqref{pp:RFDE} held everywhere on $[0,t_{\max})$, for instance in the case of piece-wise and right-continuous input signal, then \eqref{pp:ex_dini_eq_driver} would also hold everywhere on $[0,t_{\max})$. 
\end{example}

\subsection{Lyapunov-like conditions for robust forward completeness}\label{sec-RFC-cond}

As we will see throughout this survey, Driver's derivative provides a useful way to ensure several stability and robustness properties. It is also at the basis of the following two sufficient conditions for robust forward completeness. The first one relies on a functional which diverges at most exponentially along the system's solutions, which is reminiscent of the Lyapunov characterization of forward completeness proposed in \cite{ANGSON-UFC} for finite-dimensional systems.

\begin{theorem}[LKF-wise condition for RFC]\label{theo-RFC-1}
Assume that there exist a functional $V:\mathcal X^n\to\mRp$ which is Lipschitz on bounded sets, constants $a,c,\bar c\geq 0$, and functions $\underline\alpha,\overline\alpha\in\cK_\infty$, and $\gamma \in\mathcal N$ such that, for all $\phi\in\mathcal X^n$ and all $v\in\mathbb R^m$,
\begin{align}
\underline \alpha(|\phi(0)|) &\leq V(\phi)\leq \overline \alpha(\|\phi\|)+\bar c\label{eq_antoine_105} \\
D^+V(\phi,f(\phi,v))&\leq aV(\phi)+\gamma(|v|)+c. \label{eq_antoine_106}
\end{align}
Then the system \eqref{pp:RFDE} is RFC.
\end{theorem}

Condition \eqref{eq_antoine_105} imposes that $V$ can only vanish when $\phi(0)=0$ and that it tends to infinity whenever $|\phi(0)|$ tends to infinity. Note that $V$ is not requested to be zero for $\phi=0$. In view of Theorem \ref{pp:lemmadriverrfde}, condition \eqref{eq_antoine_106} imposes that, along the system's solutions, $V$ grows at most exponentially fast.

The second sufficient condition for RFC allows to bound Driver's derivative of the considered functional by a term involving the whole state history norm, but it comes at the price of requiring a particular lower bound on the functional. Its proof is provided in \cite{CHKAPEWA22}.

\begin{theorem}[History-wise condition for RFC]\label{theo-RFC-2}
Assume that the exist a functional $V:\mathcal X^n\to\mRp$ which is Lipschitz on bounded sets, $a,c,\bar c\geq 0$, $\alpha,\overline\alpha\in\cK_\infty$ and $\gamma \in\mathcal N$ such that, for all $\phi\in\mathcal X^n$ and all $v\in\mR^m$,
\begin{align}
\alpha(|\phi(0)|) &\leq V(\phi)\leq \overline \alpha(\|\phi\|)+\bar c\label{eq_antoine_107}\\
D^+V(\phi,f(\phi,v))&\leq a \alpha(\|\phi\|)+\gamma(|v|)+c.\label{eq_antoine_108}
\end{align}
Then the system \eqref{pp:RFDE} is RFC.
\end{theorem}

From the regularity assumptions made on both $V$ and $f$, the upper bound \eqref{eq_antoine_108} constitutes a mild assumption. The key requirement in the above result therefore lies in the fact $\alpha$ can also be picked as a lower bound on $V$, as imposed by \eqref{eq_antoine_107}.

\section{Input-free properties}\label{sec:no-input}

In this section, we consider an input-free version of \eqref{pp:RFDE}:
\begin{eqnarray}\label{pp:RFDE_autonomous}
\dot x(t)=f_0(x_t).
\end{eqnarray}
Such a class of systems can be seen as a particular case of \eqref{pp:RFDE} by letting $f_0(\phi):=f(\phi,0)$ for all $\phi \in\mathcal X^n$. In line with Standing Assumption \ref{stand_ass_1}, $f_0:\mathcal X^n\to \mathbb R^n$ is assumed to be Lipschitz on bounded sets and to satisfy $f_0(0)=0$.

\subsection{Definitions and equivalent formulations}\label{sec-no-input-def}

For such autonomous systems, the definitions of FC and RFC follow readily from Definitions \ref{def_FC} and \ref{def_RFC} by simply considering $u\equiv 0$. We next recall classical stability notions.

\begin {definition}[Stability properties]\label{pp:defstabetc} The origin of  (\ref{pp:RFDE_autonomous})
is said to be:
\begin{itemize}
\item \emph{asymptotically stable (AS)} if there exist $r>0$ and $\beta\in\cKL$ such that, for all $x_0\in {\cal X}^n$ with $\|x_0\|\leq r$, the corresponding solution of (\ref{pp:RFDE_autonomous}) satisfies
\begin{equation}\label{eq_antoine_3}
\vert x(t,x_0)\vert \le \beta(\Vert x_0\Vert,t),\quad \forall
t\ge 0
\end{equation}
\item \emph{exponentially stable (ES)} if there
exist $r,\lambda>0$ and $k\ge 1$ such that, for all $x_0\in
{\cal X}^n$ with $\|x_0\|\leq r$, the corresponding solution of (\ref{pp:RFDE_autonomous}) satisfies
\begin{equation}\label{eq_antoine_4}
\vert x(t,x_0)\vert \le k\Vert x_0\Vert e^{-\lambda t}, \quad
\forall t\ge 0
\end{equation}
\item \emph{globally asympotically stable (GAS)} if there exist $\beta\in\cKL$ such that \eqref{eq_antoine_3} holds for all $x_0\in {\cal X}^n$
\item \emph{globally exponentially stable (GES)} if there exist $\lambda>0$ and
$k\ge 1$ such that \eqref{eq_antoine_4} holds for all $x_0\in {\cal
X}^n$.
\end{itemize}
As far as global properties are concerned, the equilibrium is necessarily unique, so we will say with a slight abuse of terminology that the system is GAS (or GES) if its origin is GAS (or GES).
\end{definition}

In view of Lemma \ref{lem-dense}, the state estimates \eqref{eq_antoine_3}-\eqref{eq_antoine_4} need only to be checked for $x_0\in C^1([-\Delta,0],\mR^n)$ to conclude the corresponding stability property. This proves particularly handy in practice as the solutions $t\mapsto x_t$ from continuously differentiable initial states turn out to be locally absolutely continuous (see Theorem \ref{pp:lemmahalepeperfde}), and consequently, for functionals $V:\mathcal X^n\to \mathbb R_{\ge 0}$ which
are Lipschitz on bounded sets, the function $t\mapsto V(x_t)$ results locally absolutely continuous as well. 

It is worth stressing that the above definitions of AS and GAS imply some uniformity with respect to $x_0$. More precisely, the decay rate to the origin and the magnitude of transient overshoot are requested to be uniform over bounded sets of initial states. This explains why these definitions are often referred to as \emph{uniform} (global) asymptotic stability. For the sake of homogeneity with respect to the finite-dimensional terminology, we have decided not to stress this uniformity explicitly. As we discuss in Section \ref{sec:open:UGAS}, while such a uniformity comes for free in finite dimension, it is not yet known whether the same holds for time-delay systems.

A partial answer to this question lies in the following result, presented in \cite[Section 30]{KRA59}, which shows that, as far as local properties are concerned, the above uniform version of AS turns out to be equivalent to stability combined with local attractivity of the origin. This result also sheds some light about the $\cKL$ formulation of AS and some ``$\varepsilon-\delta$'' formulations employed in classical textbooks such as \cite{HALU93}, \cite{KOMY13}. 

\begin{theorem}[AS reformulation]\label{pp:equivASsolita}
The origin of (\ref{pp:RFDE_autonomous}) is AS if and only if the following two properties hold:
\begin{itemize}
\item \emph{(stability)} for every $\varepsilon>0$ there exists $\delta>0$
such that, for all $x_0\in {\cal X}^n$ with $\|x_0\|\leq \delta$, the corresponding solution  of (\ref{pp:RFDE_autonomous}) satisfies $\vert x(t,x_0)\vert
\leq\varepsilon$ for all $t\ge 0$
\item \emph{(attractivity)} there exists $r>0$ such that, for any $x_0\in {\cal X}^n$ with $\|x_0\|\leq r$, the
corresponding solution of (\ref{pp:RFDE_autonomous}) 
satisfies $ \lim_{t\to +\infty}x(t,x_0)=0$.
\end{itemize}
\end{theorem}

%\begin{theorem}[AS reformulation]\label{pp:equivASsolita}
%The following statements are equivalent:
%\begin{itemize}
%\item[a)] the origin of (\ref{pp:RFDE_autonomous}) is AS;
%\item[b)] for every $\epsilon>0$ there exists $\delta>0$
%such that, for all $x_0\in {\cal X}^n$ with $\|x_0\|\leq \delta$, the corresponding solution  of (\ref{pp:RFDE_autonomous}) satisfies $\vert x(t,x_0)\vert
%\leq\epsilon$ for all $t\ge 0$, and any of the two following properties hold:
%\begin{itemize}
%\item[$b_1$)] there exists $r>0$ such that, for every $\epsilon>0$, there exists $T\geq 0$ such that, for all $x_0\in {\cal X}^n$ with $\|x_0\|\leq r$, the
%corresponding solution of (\ref{pp:RFDE_autonomous}) 
%satisfies
%$\vert x(t,x_0)\vert \le\epsilon$ for all $t\ge T$.
%\item[$b_2$)] there exists $r>0$ such that, for any $x_0\in {\cal X}^n$ with $\|x_0\|\leq r$, the
%corresponding solution of (\ref{pp:RFDE_autonomous}) 
%satisfies
%$\lim_{t\to +\infty}x(t,x_0)=0$.
%\end{itemize}
%\end{itemize}
%\end{theorem}

In the above statement, the second item requests that any solution starting in the ball of radius $r$ eventually converge to the origin. Note that no uniformity requirement is made on the convergence rate of solutions. On the other hand, Definition \ref{pp:defstabetc} readily ensures that $|x(t,x_0)|\leq \beta(r,t)$ for all $\|x_0\|\leq r$ and all $t\geq 0$, thus imposing a common minimal convergence rate to all solutions starting in the ball of radius $r$.

When global stability properties are considered, it is not known yet whether GAS is equivalent to stability and global attractivity: see Conjecture \ref{conj:GAS} in Section \ref{sec:open}. Nevertheless, the following result, recently established in  \cite[Corollary 1]{karafyllis2022global}, shows that this indeed holds if the system is RFC.

\begin{theorem}[GAS reformulation]\label{pp:equivGASsolita} Assume that (\ref{pp:RFDE_autonomous}) is RFC. Then it is GAS if and only if the following two properties hold:
\begin{itemize}
    \item \emph{(stability)} for every $\varepsilon>0$ there exist $\delta>0$ such that, for any $x_0\in \mathcal X^n$ with $\Vert x_0\Vert\le \delta$, the corresponding solution of \eqref{pp:RFDE_autonomous} satisfies  $\vert x(t,x_0)\vert\le \varepsilon$ for all $t\ge 0$
\item \emph{(global attractivity)} for any $x_0\in {\cal X}^n$, the
corresponding solution of \eqref{pp:RFDE_autonomous} satisfies
$\lim_{t\to +\infty}x(t,x_0)=0$.
\end{itemize}
\end{theorem}
%
% \AC{I find the above statement complicated and hard to follow. I would 
% suggest to replace it by ``Assuming RFC, GAS is equivalent to stability + 
% global attractivity'' (in mathematical words, of course). Agreed, and then we % make a remark saying that RFC holds under Lagrange stability.}
%
\begin{remark}[RFC under Lagrange stability]\label{rmk-Lagrange} Notice that, if the Lagrange stability holds, then the RFC property holds. We recall that Lagrange stability means that, for every $\delta>0$ there exists $\varepsilon>0$ such that, for any 
$x_0\in \mathcal X^n$ with $\Vert x_0\Vert\le \delta$, the corresponding solution  of \eqref{pp:RFDE_autonomous} satisfies the inequality $\vert 
x(t,x_0)\vert\le \varepsilon$ for all $t\ge 0$ (notice the reverted order of quantifiers as compared to the first item of Theorem \ref{pp:equivGASsolita}). In particular, if $\delta$ can be chosen arbitrarily large for sufficiently large $\varepsilon$ in the first item of Theorem \ref{pp:equivGASsolita}, then the Lagrange stability, and consequently the RFC property, hold. 
\end{remark}

%\begin{theorem}\label{pp:lemmaequivcc1}
%The origin of (\ref{pp:RFDE_autonomous}) is GAS if and only if there exist $\beta\in\cKL$ such that, for any $x_0\in W^{1,\infty}$, the corresponding solution of (\ref{pp:RFDE_autonomous}) satisfies $|x(t,x_0)|\leq\beta(\|x_0\|,t)$ for all $t\ge 0$. 
%\end{theorem}

%\begin{remark}
%The proof of the above result being based on density arguments, similar statements also hold true for the AS, ES and GES properties. See \cite[Proposition 3]{Pepe2007} and its proof.
%\end{remark}

It is worth stressing that the stability notions of Definition \ref{pp:defstabetc} turn out to be equivalent when $f_0$ is linear. To see this, consider the linear time-invariant system
\begin{equation} \label{pp:LTI__14_}
\dot{x}(t)=A x_{t}
\end{equation}
where $A:\mathcal X^n\to\mR^n$ is a bounded linear operator. Then, we have the following result from \cite[Corollary 6.1, p. 215]{HALU93}.

\begin{theorem}[Linear case]\label{pp:LTI}
The linear system \eqref{pp:LTI__14_} is GES if and only if its origin is AS.
\end{theorem}

The class of systems \eqref{pp:LTI__14_} benefit from a wide range of analysis and control tools especially designed for linear dynamics: see for instance the textbooks \cite{BOUKASLIU2002, Fridman:2014vo, Gu:2012vm, GUANICARB, KHARITONOV2013, Kuang93,MAHMOUD,Niculescu01,MINI08}. 
Some linearization procedures allow to conclude ES of the nonlinear system \eqref{pp:RFDE_autonomous} through the study of an associated linear system: see \cite[Theorem 33.2]{KRA59}, \cite{KOMY13}, \cite[Theorem 5.2]{HAIPEP21}. One of these approaches is through the Fr\'echet derivative.

\begin{definition}[Fr\'echet derivative] The function $f_0$ in \eqref{pp:RFDE_autonomous} is said to be Fr\'echet differentiable at the origin if there exists a bounded linear operator $A:\mathcal X^n\to \mathbb R^n$ (the Fr\'echet derivative at the origin) such that
\begin{align*}
\lim_{\Vert\phi\Vert\to 0^+}\frac{\vert f_0(\phi)-A\phi\vert }{\Vert \phi\Vert}=0.
%\lim_{\phi\in \mathcal X^n, \ \Vert \phi\Vert\ne 0, \ \Vert\phi\Vert\to 0}\frac{\vert f_0(\phi)-A\phi\vert }{\Vert \phi\Vert}=0.
\end{align*}
\end{definition}

\begin{theorem}[ES through linearization]\label{theo-Frechet}
Let the function $f_0$ in \eqref{pp:RFDE_autonomous} be Fr\'echet differentiable at the origin and let $A:\mathcal X^n\to \mathbb R^n$ denote its Fr\'echet derivative at the origin. Then the origin of the nonlinear system  \eqref{pp:RFDE_autonomous} is ES if and only if the origin of the linear system $\dot x(t)=Ax_t$ is AS.
\end{theorem}

The following simple example illustrates how to apply this result in practice.

\begin{example}[ES through linearization]
Consider the scalar nonlinear system
\begin{align}\label{eq_antoine_110}
    \dot x(t)=-x(t)+x(t-\Delta)^q,
\end{align}
with $\Delta>0$ and $q>1$. For this system, $f_0(\phi)=-\phi(0)+\phi(-\Delta)^q$ for all $\phi\in \mathcal X$. The Fr\'echet derivative of $f_0$ at the origin is given by $A\phi=-\phi(0)$ for all $\phi\in\mathcal X$. The corresponding linear system is described by the equation $\dot x(t)=-x_t(0)=-x(t)$ and is finite dimensional. Its origin can easily be checked to be AS. We conclude with Theorem \ref{theo-Frechet} that the origin of the nonlinear system \eqref{eq_antoine_110} is ES.  
\end{example}

%
% \AC{Add a  result (first order approximation). Also say that we can benefit from the 
% books by Niculescu, Fridman, etc.}
%

\subsection{Lyapunov-like conditions}\label{sec-no-input-LKF}

A very useful way to study stability properties relies on the Lyapunov-Krasovskii approach, introduced in \cite{KRA59}.

\begin{definition}[LKF, coerciveness]\label{def-LKF}
A functional $V:\mathcal X^n\to\mRp$ is said to be a \emph{Lyapunov-Krasovskii functional candidate (LKF)} if it is Lipschitz on bounded sets and there exist $\underline\alpha,\overline\alpha\in\cK_\infty$ such that, for all $\phi\in\mathcal X^n$,
\begin{align}\label{pp:lowerupperboundpoint}
\underline\alpha(|\phi(0)|)\leq V(\phi)\leq \overline\alpha(\|\phi\|).
\end{align}
It is said to be a \emph{coercive LKF} if there exist $\underline\alpha,\overline\alpha\in\cK_\infty$ such that, for all $\phi\in\mathcal X^n$,
\begin{equation} \label{pp:Coer__3_}
\underline\alpha(\|\phi\|)\leq V(\phi)\leq \overline\alpha(\|\phi\|).
\end{equation}
We call these LKFs \emph{local} if \eqref{pp:lowerupperboundpoint}-\eqref{pp:Coer__3_} hold only for $\|\phi\|\leq r$ for some $r>0$.
\end{definition}

Condition \eqref{pp:lowerupperboundpoint} imposes that $V(0)=0$, that $V$ is positive whenever $\phi(0)\neq 0$, and that it grows to infinity as $|\phi(0)|\to+\infty$. Condition \eqref{pp:Coer__3_} is a stronger requirement as it imposes additionally that $V$ does not vanish unless $\phi$ is identically zero. 

Coercive LKFs provide more information on the system, but they are usually more difficult, and often less intuitive, to handle in practice. 
The following technical lemma, taken from \cite{karafyllis2010necessary}, provides a constructive way to get a coercive LKF based on the knowledge of a non-coercive one. As we will see through Example \ref{exa-dissipations}, this lemma in turn provides interesting properties on its Driver's derivative, which can be useful in stability analysis.  

\begin{lemma}[Construction of coercive LKFs]\label{means} Let $W:\mathbb R^{n}\to \mRp$ be a continuously differentiable function. Then, given any $c >0$, the functional $V:\mathcal X^n\to\mRp$ defined for all $\phi \in \mathcal X^{n} $ as
\begin{equation*} %\label{lemma__1_} 
V(\phi ):=\mathop{\max }\limits_{\tau \in [-\Delta ,0]} W(\phi (\tau ))e^{2c\tau} 
\end{equation*} 
is Lipschitz on bounded sets and satisfies the following implications for all $w\in {\mathbb R}^{n} $: 
\begin{align*} 
W(\phi (0))<V(\phi )\quad &\Rightarrow\quad D^{+} V(\phi ,w)\le -2c V(\phi ) %\label{lemma__2_} 
\\
W(\phi (0))=V(\phi )\quad & \Rightarrow\quad D^{+} V(\phi ,w)\le \max \left\{-2c V(\phi )\, ,\, \nabla W(\phi (0))w\right\}. %\label{lemma__3_} 
\end{align*} 
\end{lemma}

The following result shows that asymptotic stability can be characterized in LKF terms. See \cite[Theorems 30.2, 31.3]{KRA59} for the AS property and \cite[Theorem 2.3]{Pepe:2013it}-\cite[Proposition 2]{CHGOPE21} for the GAS property. See also \cite{HALU93,KOMY13}.

\begin{theorem}[LKF characterization of GAS/AS]\label{pp:teoremafondamentalekargas}
The following statements are equivalent:
\begin{itemize}
\item \eqref{pp:RFDE_autonomous} is GAS
\item there exist an LKF $V:\mathcal X^n\to\mRp$ and $\sigma\in\mathcal{KL}$ such that, for all $\phi\in\mathcal X^n$, 
\begin{align}\label{eq_antoine_5bis}
D^+V(\phi,f_0(\phi))\le -\sigma(|\phi(0)|,\|\phi\|)
\end{align}
\item there exist an LKF $V:\mathcal X^n\to\mRp$ and $\alpha\in\mathcal P$ such that, for all $\phi\in\mathcal X^n$, 
\begin{align}\label{eq_antoine_5}
D^+V(\phi,f_0(\phi))\le -\alpha(|\phi(0)|)
\end{align}
\item there exist an LKF $V:\mathcal X^n\to\mRp$ such that, for all $\phi\in\mathcal X^n$, 
\begin{align}\label{eq_antoine_6}
D^+V(\phi,f_0(\phi))\le -V(\phi)
\end{align}
\item there exist a coercive LKF $V:\mathcal X^n\to\mRp$ and $\alpha\in\mathcal K_\infty$ such that, for all $\phi\in\mathcal X^n$, 
\begin{align}\label{eq_antoine_7}
D^+V(\phi,f_0(\phi))\le -\alpha(\|\phi\|).
\end{align}
\end{itemize}
Moreover, the same equivalences hold for AS if \eqref{eq_antoine_5}, \eqref{eq_antoine_6} and \eqref{eq_antoine_7} hold with a local LKF and only for $\|\phi\|\leq r$ for some $r>0$.
\end{theorem}

Condition \eqref{eq_antoine_5} is probably the handiest way to establish GAS (or AS) as the considered LKF is not requested to be coercive and its dissipation rate is only required to involve the current value of the solution's norm (point-wise dissipation). It can actually be even relaxed to \eqref{eq_antoine_5bis}, in which the dissipation rate is allowed to be smaller when $\|\phi\|$ gets larger ($\cKL$ dissipation), which turns out particularly useful when dealing with LKFs of the form $V(\phi)=\ln(1+W(\phi))$ where $W$ denotes another LKF. It is interesting to notice that any of these two conditions ensures the existence of an LKF which dissipates exponentially along the systems solution (condition \eqref{eq_antoine_6}), and even the existence of a coercive LKF that dissipates in terms of the whole state history (history-wise dissipation, condition \eqref{eq_antoine_7}), which may prove useful for further robustness analysis.

%
% \AC{Replace the 3rd item by $D^+V(\phi,f_0(\phi))\le -V(\phi)$?}
%
%\begin{theorem} The system described by (\ref{pp:RFDE}) is $0$-GAS if and only if
%there exist a locally Lipschitz functional $V:\mathcal {X}^n\to
%\mathbb R_{\ge 0}$ and functions $\underline \alpha$, $\overline \alpha$ of class
%${\cal K}_{\infty}$,  $\alpha$ of class ${\cal K}$, such that, $\forall \phi \in \mathcal {X}^n$,
%the following inequalities hold:
%\smallskip
%\item {i)} $\underline \alpha(\vert \phi(0)\vert)\le V(\phi)\le \overline \alpha(\Vert \phi\Vert_{\infty})$;
%\smallskip
%\item {ii)} $D^+V(\phi,0)\le -a_3(\vert \phi(0)\vert)$
%\smallskip
%
%\end{theorem}

%\begin{theorem}[LKF characterization of ES]\label{pp:infty} The following statements are equivalent:
%\begin{itemize}
%\item [1)] the system described by (\ref{pp:RFDE_autonomous}) is ES;
%\item [2)] there exist positive reals $H, \underline a, \overline a, a, p$ and a LKF $V:\mathcal {X}^n\to \mathbb R_{\ge 0}$ such that, for all $\phi \in \mathcal {X}^n$ satisfying $\|\phi\|\leq H$, the following inequalities hold
%\begin{align*}
%\underline a\Vert \phi\Vert^p\le V(\phi)\le %\overline a\Vert \phi\Vert^p \\ \\
%D^+V(\phi,f_0(\phi))\le -a\Vert \phi \Vert^p.
%\end{align*}
%\end{itemize}
%\end{theorem}

Exponential stability can also be characterized using LKFs, as stated next. See \cite[Lemma 33.1]{KRA59} for the ES property and \cite[Theorems 2.3, 2.4, 2.5]{Pepe:2013it} or \cite{haidar2021lyapunov} for the GES property.

\begin{theorem}[LKF characterizations of GES/ES]\label{pp:theoremgesinfty} The following statements are equivalent:
\begin{itemize}
\item (\ref{pp:RFDE_autonomous}) is GES
\item there exist an LKF $V:\mathcal {X}^n\to \mathbb R_{\ge 0}$ and $\underline a,\overline a, p>0$ such that, for all $ \phi \in \mathcal {X}^n$, 
\begin{align}
\underline a|\phi(0)|^p\le V(\phi)\le \overline a\Vert \phi\Vert^p \label{eq_antoine_8}\\
D^+V(\phi,f_0(\phi))\le -aV(\phi) \label{eq_antoine_9}
\end{align}
\item there exist an LKF $V:\mathcal {X}^n\to \mathbb R_{\ge 0}$ and $\underline a,\overline a, a, p>0$ such that, for all $\phi \in \mathcal {X}^n$, 
\begin{align}
\underline a\Vert \phi\Vert^p\le V(\phi)\le \overline a\Vert \phi\Vert^p \label{eq_antoine_10o}\\
D^+V(\phi,f_0(\phi))\le -a\Vert \phi \Vert^p. \label{eq_antoine_11}
\end{align}
\end{itemize}
Moreover, the same equivalence holds for ES if the above bounds hold only for $\|\phi\|\leq r$ for some $r>0$.
\end{theorem}

Here again, conditions \eqref{eq_antoine_8}-\eqref{eq_antoine_9} are usually easier to invoke in practice as they do not require a coercive LKF. Nevertheless, the dissipation is requested to involve the whole LKF itself. As far as local properties are concerned, a point-wise dissipation is also sufficient to ensure ES, as stated in the following result, which can be derived from \cite{CHORPE19} or \cite{CHKAPEWA22}.

\begin{theorem}[ES under point-wise dissipation]\label{pp:ESteoremapointwise}
The origin of (\ref{pp:RFDE_autonomous}) is ES if and only if there exist an LKF $V:\mathcal {X}^n\to \mathbb R_{\ge 0}$ and $\underline a,\overline a, a,r>0$ such that, for all $\phi \in \mathcal {X}^n$ with $\|\phi\|\leq r$,
\begin{align*}
\underline a\vert \phi(0)\vert^2\le V(\phi)\le \overline a\Vert \phi\Vert^2 \\
D^+V(\phi,f_0(\phi))\le -a\vert \phi(0) \vert^2.
\end{align*}
\end{theorem}

As discussed in Section \ref{sec:open:GES}, it is not known yet whether GES can be established through a point-wise dissipation in general. Nevertheless, additional growth conditions have been proposed in \cite{CHKAPEWA22} to make this possible.

\begin{theorem}[GES under point-wise dissipation]\label{pp:gesteoremapointwise}
Assume that there exist a functional $V:\mathcal X^n\to \mRp$ which is Lipschitz on bounded sets, $\underline a,\overline a,a>0$, and $\varepsilon\ge 0$ such that, for all $\phi \in \mathcal X^n$ and $v\in \mR^n$,
\begin{align}
\underline a \left|\phi (0)\right|^{2} &\le V(\phi )\le \overline a \left\| \phi \right\| ^{2}  \label{eq_antoine_12} \\ 
D^{+} V(\phi ,f_0(\phi))&\le -a \left|\phi (0)\right|^{2} +\varepsilon\left\| \phi \right\| ^{2}.  \label{eq_antoine_13}
\end{align} 
Assume further that there exists a symmetric positive definite matrix $P\in \mR^{n\times n} $ and a constant $c>0$ such that any of the two following conditions is satisfied:
\begin{align} 
\phi(0)^\top P\, f_0(\phi)&\leq c\left\| \phi \right\| ^{2},\quad \forall \phi\in\mathcal X^n \label{eq-antoine-110}\\
\phi(0)^\top P\, f_0(\phi)&\ge -c\left\| \phi \right\| ^{2},\quad \forall \phi\in\mathcal X^n.\label{eq-antoine-111}
\end{align} 
Then there exists $\bar \varepsilon >0$ such that, if $\varepsilon\in[0,\bar \varepsilon)$, the system \eqref{pp:RFDE} is GES.
\end{theorem}

This result provides a point-wise LKF condition for GES when focused on the case when $\varepsilon=0$. A similar result (for $\varepsilon=0$) was originally given in \cite{CHORPE19} under the assumption that $f_0$ is globally Lipschitz. It is worth noting that any of the conditions \eqref{eq-antoine-110}-\eqref{eq-antoine-111} is trivially satisfied with $P=I$ if $f_0$ is globally Lipschitz. Nevertheless, they are far from being confined to such a requirement (see the examples in \cite{CHKAPEWA22}). Interestingly, Theorem \ref{pp:gesteoremapointwise} allows for the presence of a positive term $\varepsilon\|\phi\|^2$, which may prove useful in establishing GES under modeling uncertainty. This feature would be trivial if the considered LKF were coercive and with a history-wise dissipation, but becomes more interesting here (as \eqref{eq_antoine_13} does not even guarantee that $D^{+} V(\phi ,f_0(\phi))\leq 0$ if $\varepsilon>0$). Note finally that \cite{CHKAPEWA22} provides explicit estimates of $\bar \varepsilon$, which are not reported here.

\begin{example}[Neuronal population]\label{exa_neuro_1}
Consider the following scalar system 
\begin{align}\label{eq_antoine_14}
    \dot x(t)=-x(t)+g\big(x(t-\Delta)\big),
\end{align}
where $g:\mR\to\mR$ denotes a globally Lipschitz function with Lipschitz constant $\ell>0$ and satisfying $g(0)=0$. Such class of systems is sometimes employed to model the dynamics of a neuronal population, in which case the delay $\Delta\geq 0$ reflects the non-instantaneous propagation of spikes along the axons. We claim that, regardless of the value of $\Delta$, the origin of this system is GES provided that
\begin{align}\label{eq_antoine_15}
    \ell <1.
\end{align}
To see this, consider the LKF defined as
\begin{align*}
V(\phi) := \frac{1}{2}\left(\phi(0)^2 +  \int_{-\Delta}^0 \phi(\tau)^2 d\tau\right), \quad \forall \phi\in\mathcal X.
\end{align*}
Such a functional clearly satisfies condition \eqref{eq_antoine_12} since
\begin{align*}
    \frac{1}{2}|\phi(0)|^2\leq V(\phi)\leq \frac{1+\Delta}{2}\|\phi\|^2.
\end{align*}
Defining $f_0$ as the right-hand side of \eqref{eq_antoine_14} and reasoning as in Example \ref{exa1}, it holds that
\begin{align*}
    D^+V(\phi,f_0(\phi)) &= -\phi(0)^2+\phi(0)g\big(\phi(-\Delta)\big)+\frac{1}{2}\left(\phi(0)^2-\phi(-\Delta)^2\right).
\end{align*}
Using the assumptions made on $g$, it holds that $|g(z)|\leq \ell|z|$ for all $z\in\mR$. Recalling that $ab\leq (\ell a^2+b^2/\ell)/2$ for all $a,b\in\mR$, it follows that
\begin{align*}
    D^+V(\phi,f_0(\phi)) &\leq -\phi(0)^2+\ell|\phi(0)||\phi(-\Delta)|+\frac{1}{2}\left(\phi(0)^2-\phi(-\Delta)^2\right)\\
    &\leq -\frac{1}{2}\phi(0)^2+\frac{\ell}{2}\left(\ell\phi(0)^2+\frac{\phi(-\Delta)^2}{\ell}\right)-\frac{1}{2}\phi(-\Delta)^2\\
    &\leq -\frac{1}{2}\left(1-\ell^2\right)\phi(0)^2.
\end{align*}
In view of \eqref{eq_antoine_15}, this makes condition \eqref{eq_antoine_13} fulfilled. Finally, since $g$ is globally Lipschitz, any of the two conditions \eqref{eq-antoine-110}-\eqref{eq-antoine-111} is fulfilled with $P=I$ and GES follows from Theorem \ref{pp:gesteoremapointwise}. Note that, for this system, an LKF-wise dissipation rate can easily be obtained by considering the alternative LKF
\begin{align*}
W(\phi) := \frac{1}{2}\left(\phi(0)^2 +  \int_{-\Delta}^0 e^{\rho\tau}\phi(\tau)^2 d\tau\right),
\end{align*}
for some $\rho>0$: this is formally shown in \cite[Lemma 1]{ORCHDESI21}. Nevertheless, the point-wise dissipation required by \eqref{eq_antoine_13} significantly simplifies the analysis.
\end{example}

The following theorem provides an alternative methodology for checking the GAS property by means of the Razumikhin methodology. See \cite{Razumikhin56}, \cite[Theorem 4.2]{HALU93}, and \cite[Proposition 4.2]{Karafyllis:2008hc}. %The conditions provided in next theorem imply the RFC property. 
\begin{theorem}[GAS through Lyapunov-Razumikhin]\label{pp:razum-GAS}
Assume that there exist a continuously differentiable function $V:\mathbb R^n\to \mathbb R_{\ge 0}$ and functions $\underline \alpha, \overline \alpha \in \mathcal K_{\infty}$ such that, for all $z\in\mathbb R^n$,
\begin{align*}
    \underline \alpha(\vert z\vert)\le V(z)\le \overline \alpha (\vert z\vert ).
\end{align*}
Assume further that there exist $\alpha\in\mathcal P$ and $\rho\in\cK_{\infty}$, such that, for all $\phi\in\mathcal X^n$,
\begin{align*}
V(\phi(0))\ge \rho\left(\max_{\tau\in [-\Delta,0]}V(\phi(\tau))\right)\quad \Rightarrow\quad \nabla V(\phi(0))f_0(\phi)\le -\alpha(\vert \phi(0)\vert).
\end{align*}
Then the system \eqref{pp:RFDE_autonomous} is GAS provided that $\rho(s)<s$ for all $s>0$.
\end{theorem}

GES can also be established using the Lypunov-Razumikhin approach, as stated in the following result, whose proof can be obtained by invoking \cite[Lemma 2.14]{KAJIbook11}. We stress that we could not find any published version of this result, although similar ones were provided in \cite{wang2005exponential,dashkovskiy2012stability} in the context of time-delay systems with resets.

\begin{theorem}[GES through Lyapunov-Razumikhin]\label{pp:razum-GES}
Assume that there exist a continuously differentiable function $V:\mathbb R^n\to \mathbb R_{\ge 0}$ and $\underline a,\overline a, p>0$ such that, for all $z\in\mathbb R^n$,
\begin{align*}
    \underline a\vert z\vert^p\le V(z)\le \overline a\vert z\vert^p.
\end{align*}
Assume further that there exist $a>0$ and $\rho\in(0,1)$ such that, for all $\phi\in\mathcal X^n$,
\begin{align*}
V(\phi(0))\ge \rho\max_{\tau\in [-\Delta,0]}V(\phi(\tau))\quad \Rightarrow\quad \nabla V(\phi(0))f_0(\phi)\le -aV(\phi(0)).
\end{align*}
Then the system \eqref{pp:RFDE_autonomous} is GES.
\end{theorem}

%\begin{theorem}[GES through Lyapunov-Razumikhin]\label{pp:razum-GES}
%Assume that there exist a continuously differentiable function $V:\mathbb R^n\to \mathbb R_{\ge 0}$ and $\underline a,\overline a, p>0$ such that, for all $z\in\mathbb R^n$,
%\begin{align*}
%    \underline a\vert z\vert^p\le V(z)\le \overline a\vert z\vert^p.
%\end{align*}
%Assume further that there exists $a>0$ such that, for all $\phi\in\mathcal X^n$,
%\begin{align*}
%V(\phi(0))\ge \max_{\tau\in [-\Delta,0]}V(\phi(\tau))e^{a\tau}\quad \Rightarrow\quad \nabla V(\phi(0))f_0(\phi)\le -aV(\phi(0)).
%\end{align*}
%Then the system \eqref{pp:RFDE_autonomous} is GES.
%\end{theorem}

A related way to establish GES is through Halanay's inequality \cite[p. 378]{Halanay:1966wl}.
\begin{theorem}[GES through Halanay]\label{pp:Halanay-GES}
Assume that there exist a continuously differentiable function $V:\mathbb R^n\to \mathbb R_{\ge 0}$ and $\underline a,\overline a, p>0$ such that, for all $z\in\mathbb R^n$,
\begin{align*}
    \underline a\vert z\vert^p\le V(z)\le \overline a\vert z\vert^p.
\end{align*}
Assume further that there exists $a,b>0$ such that, for all $\phi\in\mathcal X^n$,
\begin{align*}
\nabla V(\phi(0))f_0(\phi)\leq -aV(\phi(0))+b\max_{\tau\in [-\Delta,0]}V(\phi(\tau)).
\end{align*}
Then the system \eqref{pp:RFDE_autonomous} is GES provided that $a>b$.
\end{theorem}

This result has more recently been extended to GAS by relying on a nonlinear version of Halanay's inequality \cite[Corollary 1]{Pepe_Halanay22}.
\begin{theorem}[GAS through Halanay]\label{pp:Halanay_GAS}
Assume that there exist a continuously differentiable function $V:\mathbb R^n\to \mathbb R_{\ge 0}$ and $\underline \alpha,\overline \alpha\in\cK_\infty$ such that, for all $z\in\mathbb R^n$,
\begin{align*}
    \underline \alpha(\vert z\vert)\le V(z)\le \overline \alpha(\vert z\vert).
\end{align*}
Assume further that there exists $\alpha\in \mathcal P$ and  $\gamma\in \mathcal K$ such that, for all $\phi\in\mathcal X^n$,
\begin{align*}
\nabla V(\phi(0))f_0(\phi)\leq -\alpha(V(\phi(0)))+\gamma\left(\max_{\tau\in [-\Delta,0]}V(\phi(\tau))\right).
\end{align*}
Then the system \eqref{pp:RFDE_autonomous} is GAS provided that $\alpha-\gamma\in\mathcal P$ 
%Let there exist a continuously differentiable function $V:\mathbb R^n\to \mathbb R_{\ge 0}$, functions $\underline \alpha$, $\overline \alpha$ of class $\mathcal K_{\infty}$, functions $\alpha$ and $\gamma$ of class $\mathcal P$ and $\mathcal K$, respectively, with $\alpha-\gamma$ of class $\mathcal P$, such that, for any $x\in \mathbb R^n$, $\phi\in \mathcal X^n$, the following inequalities hold:
%\begin{itemize}
%\item[$i)$]\ $\underline \alpha(\vert x\vert)\le V(x)\le \overline \alpha(\vert x\vert)$;
%\item [$ii)$]\ $\nabla V(\phi(0))f_0(\phi)\le -\alpha(V(\phi(0)))+\gamma\left(\max_{\tau\in [-\Delta,0]}V(\phi(\tau))\right)$,
%\end{itemize}
%Then the origin of \eqref{pp:RFDE_autonomous} is GAS.
\end{theorem}

\subsection{Output stability properties}\label{yw:sec-autonomous}

\newcommand{\R}{\mathbb R}

In Sections \ref{sec-no-input-def} and \ref{sec-no-input-LKF}, we have presented tools to study stability of the origin in terms of all state variables. In several applications, this constitutes a too demanding requirement, for instance when only particular state variables converge to the zero. This is particularly true when dealing with adaptive control, in which case the state of the closed-loop system is made of both the state of the plant and the parameter estimation errors: the control objective is then to drive the plant state to the origin, but it is not necessarily guaranteed that the estimation error will vanish. Other applications arise from tracking and output regulation problems where only the tracking error is requested to converge to zero. Similarly, in observer design, the aim is to cancel the observation error while the actual state may not have any prescribed behavior. 
%Another useful setting in which such a limitation occurs is for stability with respect to a set where the distance to the set is treated as the output. 
%The possibility to address stability of only part of the state variables can also prove useful when studying time-varying systems $\dot x(t) = f(t, x_t)$ by treating the time variable $t$ as an extra state variable with $\dot t=1$, which, of course, will never converge to the origin. 

This kind of properties can be analyzed using output stability notions, which impose stability properties on an output $y$ only, rather than the full state $x$. In this section, we therefore consider systems as in \eqref{pp:RFDE_autonomous} with an output map:
\begin{subequations}\label{yw:sys-auto}
\begin{align}
\dot x(t) &= f_0(x_t), \label{yw:sys-auto1}\\
y(t) &= h(x_t), \label{yw:sys-auto2}
\end{align}
\end{subequations}
where $y(t)\in\R^p$ represents the output we are interested in. 
%\YW{Theorem 19 used the assumption that $f(0)=0$, so this remark should be dropped. In view of the possible applications listed above, we do not request in this section that $f_0(0)=0$. Nevertheless, we still assume that $f_0$ is Lipschitz on bounded sets.}\AC{Do you suggest to assume that Standing Assumption 1 holds in this section? This is problematic, as it does not allow to cover time-varying systems by letting $\dot t=1$\ldots Probably not a big deal: we can assume Standing Assumption 1 in the whole paper.} \YW{I think it is better to still use the Standing Assumption.  This will cause an issue for application to time varying systems, but that's just a very specific case. To keep the homogeneity, we perhaps omit the time varying application? }
Beyond Standing Assumption \ref{stand_ass_1}, we also assume throughout this section that the system is forward complete, which can be established for instance with Theorem \ref{theo-RFC-1} or \ref{theo-RFC-2}. The output map $h:\mathcal X^n\rightarrow\R^p$ in (\ref{yw:sys-auto2}) is assumed to be Lipschitz on bounded sets with $h(0)=0$. One may consider a more general output map $h:\mathcal X^n\to\mathcal Y$, where $\mathcal Y$ is a normed linear space with a norm $\|\cdot\|_{\mathcal Y}$.  However, in the context of output stability, one can in general reformulate the problem by considering the output map $\tilde h:\mathcal X^n\rightarrow\R^p$ given by $\tilde h(\phi) := \|h(\phi)\|_{\mathcal Y}$. Given any $x_0\in\mathcal X^n$, we let $y(t,x_0) := h(x(t, x_0))$.

%Note that it follows from the Lipschitz condition on $h$ that there exists $\pi\in\cK_\infty$ such that
%\[
%|h(\phi)|\le \pi(\|\phi\|), \quad\forall\,\phi\in\mathcal X^n.
%\]

In some applications, the considered output depends only on the current value of the solutions, namely $h(x_t) = h_0(x(t))$ for some continuous map $h_0:\mR^n\to\mR^p$. This constitutes a particular case of the following notion, introduced in \cite[Definition 3.8]{KAPEJI08}.

\begin{definition}[Equivalence to finite-dimensional output map]
\label{yw:fd-output}
Suppose there exists a continuous map $h_0: \R^n\rightarrow\R^p$ with $h_0(0)=0$ and $\underline\alpha,\overline\alpha\in{\mathcal K}_\infty$ such that, for all $\phi\in\mathcal X^n$,
\[
\underline\alpha(|h_0(\phi(0))|)\le |h(\phi)|\le \overline\alpha\left(\sup_{\tau\in[-\Delta, 0]}|h_0(\phi(\tau))|\right).
\]
Then we say that $h:\mathcal X^n\rightarrow\R^p$ is \emph{equivalent to the finite-dimensional mapping} $h_0$.
\end{definition}

For finite dimensional systems, output-stability notions were considered in \cite{SONWANIOS} and \cite{TEEPRAconverse} in the ISS framework. Such properties were also referred to as \emph{stability with respect to two measures} \cite{LALI93} or
\emph{partial asymptotic stability}  \cite{Vorotnikov93}. For delay systems, the following property was introduced in \cite{KAPEJI08}.

\begin{definition}[GAOS]\label{yw:def-os}
A forward complete system (\ref{yw:sys-auto}) is said to be
{\it globally asymptotically output stable (GAOS)} if there exists $\beta\in \mathcal{ KL}$ such that, for all $x_0\in\mathcal X^n$,
\begin{equation}\label{yw:gaos}
|y(t,x_0)|\le \beta(\|x_0\|, t), \quad\forall\,t\ge 0.
\end{equation}
\end{definition}

Here again, Lemma \ref{lem-dense} ensures that GAOS holds if the state estimate \eqref{yw:gaos} is satisfied for all $x_0\in C^1([-\Delta,0],\mR^n)$. Clearly, GAOS boils down to GAS in the particular case when $h(\phi)=\phi(0)$. The GAOS property readily implies that the output transient overshoot is arbitrarily small for initial states sufficiently close to the origin. It also guarantees that the output norm converges to zero from any initial state. More precisely, GAOS ensures the following properties:
\begin{itemize}
%\item \emph{Output global stability (OGS)}:
\item \emph{global output stability (GOS)}:
there exists some $\zeta\in\cK$ such that, for all $x_0\in\cX^n$,
\begin{align}\label{eq-antoineOGS}
|y(t, x_0)|\le\zeta(\|x_0\|),\quad\forall\,t\ge 0
\end{align}
\item \emph{global output attractivity}: for each $x_0\in\cX^n$,
\begin{align}\label{eq-antoine-OAtt}
\lim_{t\to+\infty}y(t, x_0) = 0.
\end{align}
\end{itemize}
%However, the converse implication does not hold even for finite dimensional delay-free system: see for instance the counter-example in \cite{ORCHSI-lcss20}. Accordingly, the following non-uniform notion of GAOS was proposed in \cite{KACH20}.

%\begin{definition}[Non-uniform GAOS]\label{yw:def-non-unif-GAOS}
%A forward complete system (\ref{yw:sys-auto}) is said to be \emph{non-uniformly GAOS} if it satisfies the GOS and OAtt properties.
%\end{definition}

Unlike the $\cKL$ estimate \eqref{yw:gaos}, global output attractivity does not guarantee that the rate at which the output converges is uniform over bounded sets of initial states. In view of Remark \ref{rmk-Lagrange}, the GAOS property coincides with the combination of GOS and global output attractivity in the case when $h(\phi)=\phi(0)$ (since GOS implies in particular Lagrange stability in that case). A natural question is then whether Theorem \ref{pp:equivGASsolita} extends to GAOS, namely if non-uniform GAOS implies GAOS under the RFC assumption. A negative answer to this question was provided through a counterexample in \cite{ORCHSI-lcss20}.

%\YW{ A minor conjecture is that for a RFC system, the OL-GOS and the O-Att may lead to the OL-version of the GAOS property:
%\[
%|y(t, \phi)|\le \beta\left(|h(\phi)|, \frac{t}{1+\kappa(\|\phi\|)}\right) \ (\beta\in\cKL, \kappa\in\cK),
%\]
%where OL-GOS mean:
%\[
%|y(t, \phi)|\le\zeta(|h(\phi)|)\quad\forall\,t\ge 0.
%\]
%This should be an output version of (c) $\Rightarrow$ (a) in Theorem \ref{pp:equivASsolita}.)
%}\AC{Worth considering in the open questions section.}

%Some more types of output stability properties can be found in the context of robust output stability or input-to-output stability and will be detailed in Section \ref{yw: sec-ios}.

%For an autonomous system as in \eqref{yw:sys-auto}, the RFC condition becomes
%\[
%\sup\left\{
%\|x_t(x_0)\|: \ x_0\in\mathcal X^n, \ \|x_0\|\le r, \ t\in[0, T]\right\} < +\infty,\quad \forall T, r>0.
%\]
The following LKF characterization of GAOS was stated in \cite[Theorem 4.2]{KAPEJI08}.

\begin{theorem}[LKF characterization of GAOS]\label{yw:theorem-gaos}
Assume that the system \eqref{yw:sys-auto} is RFC. Then the following statements are equivalent:
\begin{itemize}
\item[$i)$] \eqref{yw:sys-auto} is GAOS
\item[$ii)$] there exist a functional $V:\mathcal X^n\rightarrow\R_{\ge 0}$ which is Lipschitz on bounded sets, $\underline{\alpha}, \overline{\alpha}\in {\mathcal K}_\infty$, and $\alpha\in\mathcal P$ such that, for all $\phi\in\mathcal X^n$, 
\begin{align}
\underline{\alpha}(|h(\phi)|)\le V(\phi)\le \overline{\alpha}(\|\phi\|) \label{yw:os-eq9}\\
D^+V(\phi, f_0(\phi))\le -\alpha(V(\phi)) \label{yw:os-lyap2n}
\end{align}
\item[$iii)$] there exist a functional $V:\mathcal X^n\rightarrow\R_{\ge 0}$ which is Lipschitz on bounded sets and $\underline{\alpha}, \overline{\alpha}\in {\mathcal K}_\infty$ such that, for all $\phi\in\mathcal X^n$, \eqref{yw:os-eq9} holds and
\begin{align}
D^+V(\phi, f_0(\phi))\le -V(\phi). \label{yw:os-eq10}
\end{align}
\end{itemize}
Moreover,
%\begin{itemize}
%\item[$i)$] if $h:\mathcal X^n\rightarrow\R^m$ is equivalent to a finite-dimensional continuous mapping
%$h_0: \R^n\rightarrow\R^p$, then inequalities \eqref{yw:os-eq9} can be replaced by
%\begin{equation}\label{yw:os-eq7}
%\underline{\alpha}(|h_0(\phi(0))|)\le V(\phi)\le \overline{\alpha}(\|\phi\|),
%\quad
%\forall\,\phi\in\mathcal X^n.
%\end{equation}
%\item[$ii)$] 
if there exist $\eta\in{\mathcal K}_\infty$ and $R\ge 0$ such that
\begin{equation}\label{yw:os-eq8}
\eta(|\phi(0)|)\le V(\phi) + R, \quad\forall\,\phi\in \mathcal X^n,
\end{equation}
then the RFC condition is not needed.
%\end{itemize}
\end{theorem}

It is worth stressing that \eqref{yw:os-eq9} requires a lower bound on $V$ only in terms of the output norm; hence, $V$ is not requested to be an LKF in the sense of Definition \ref{def-LKF}.

This theorem states that GAOS is ensured provided that this function $V$ admits an LKF-wise dissipation (condition \eqref{yw:os-lyap2n}). Condition \eqref{yw:os-eq10} then states that there exists a functional $V$ that dissipates exponentially along the system's solutions.

Condition \eqref{yw:os-eq8} imposes that the considered functional $V$ is lower bounded by a function of the full state (rather than the output only) modulo a constant $R\geq 0$. It turns out particularly useful when considering practical stability properties, meaning stability of a ball of radius $r>0$ centered at the origin (rather than the origin itself). In this case, the output can be picked as the distance with respect to this ball, namely $y(t)=\max\{|x(t)|-r,0\}$, meaning that $h(\phi)=\max\{|\phi(0)|-r,0\}$ for all $\phi\in\mathcal X^n$.

\begin{remark}\label{yw:remark-gaos1}
In Theorem \ref{yw:theorem-gaos}, the RFC assumption can be relaxed to FC for the implications of 
$ii)$ $\Leftrightarrow$ $iii)$ and $iii)$ $\Rightarrow$ $i)$. 
\end{remark}

In view of Theorem \ref{pp:teoremafondamentalekargas}, it seems natural to ask whether GAOS would also hold under a $\cKL$ dissipation of the form 
\begin{align}
D^+V(\phi, f_0(\phi))\le -\sigma(V(\phi),\|\phi\|) \label{yw:os-lyap2nbis}
\end{align}
for some $\sigma\in\cKL$, meaning that the dissipation rate is allowed to be smaller when $\|\phi\|$ gets larger. The following example shows that is not the case in general, even in finite dimension.

\begin{example}[RFC + $\cKL$ dissipation $\nRightarrow$ GAOS]
Consider the finite-dimensional system
\begin{subequations}\label{eq-antoine-113}
\begin{align}
    \dot x_1(t) &= x_1(t)\\
    \dot x_2(t) &= -\frac{x_2(t)}{1+x_1(t)^2}.
\end{align}
\end{subequations}
RFC of this system can easily be shown with Theorem \ref{theo-RFC-1}. Letting $V(\phi):=\phi_2(0)^2$ for all $\phi=(\phi_1,\phi_2)^\top\in\mathcal X^2$, it readily holds that 
\begin{align*}
    |\phi_2(0)|^2\leq V(\phi)\leq \|\phi\|^2
\end{align*}
thus fulfilling \eqref{yw:os-eq9} with the output map $h(\phi):=\phi_2(0)$. Moreover, defining $f_0$ as the vector field of \eqref{eq-antoine-113}, it holds that
\begin{align*}
    D^+V(\phi,f_0(\phi))=-\frac{2\phi_2(0)^2}{1+\phi_1(0)^2}\leq -\frac{2V(\phi)}{1+\|\phi\|^2},
\end{align*}
thus establishing \eqref{yw:os-lyap2nbis} with $\sigma(r,s):=2r/(1+s^2)$ for all $r,s\geq 0$. Nevertheless, the system is not GAOS. To see this, consider any initial state satisfying $x_0(0)=(1,1)^\top$. Then it holds that $x_1(t)=e^t\geq t$ for all $t\geq 0$. Consequently,
\begin{align*}
    x_2(t) = x_2(0)\exp\left(-\int_0^t \frac{d\tau}{1+x_1(\tau)^2}\right)\geq \exp\left(-\int_0^t \frac{d\tau}{1+\tau^2}\right)\geq e^{-\pi/2}.
\end{align*}
Thus, the considered output does not converge to zero, meaning that the system is not GAOS.
\end{example}

Nevertheless, we will see in Theorem \ref{yw:yw-ios-thm} that such a $\cKL$ dissipation does ensure GAOS provided that the RFC assumption is replaced by the stronger requirement that $|x(t,x_0)|\leq \eta(\|x_0\|)$ for all $t\geq 0$ for some $\eta\in\cK_\infty$, meaning that the whole state is bounded and that the origin is stable (which is clearly violated in the above example).

All the conditions of Theorem \ref{yw:theorem-gaos} impose to obtain negative terms in $V$ in its Driver's derivative. In the analysis of output stability properties, this sometimes constitutes a difficult task. It was shown through a counter-example in \cite{ORCHSI-lcss20} that the existence of a functional $V$ satisfying
\begin{align*}
    D^+V(\phi,f_0(\phi))\leq -\alpha(|h(\phi)|)
\end{align*}
with some $\alpha\in\cK_\infty$ is not enough to conclude GAOS (even in a finite-dimensional context). In other words, a dissipation rate involving merely the output norm is not enough for GAOS: while such a condition does ensure GOS and global output attractivity, the uniformity requirement on bounded sets of initial states is not met. The following result, proposed in \cite{KACH20}, provides an additional condition under which GAOS can be derived from a dissipation rate involving the output only.

\begin{theorem}[GAOS under non-LKF-wise dissipation]\label{yw:thm2-gaos}
Let \eqref{yw:sys-auto} be forward complete and assume there exist two functionals $V, W:\mathcal X^n\rightarrow\R_{\ge 0}$ that are Lipschitz on bounded sets, %and satisfying $\sup_{\|\phi\|\leq s}V(\phi)+ W(\phi) <+\infty$ for all $s\ge 0$. 
$\underline\alpha\in\mathcal K_\infty$, and $\alpha\in\mathcal P$ such that, for all $\phi\in\mathcal X^n$,
\begin{align}
\underline \alpha(|h(\phi)|)&\le W(\phi)
\label{yw:gaos-eq05}\\
D^+V(\phi, f_0(\phi))&\le-\alpha(W(\phi)).
\end{align}
Then, under the condition that
\begin{align}
D^+W(\phi, f_0(\phi))\le 0, \quad\forall\,\phi\in\mathcal X^n,\label{yw:gaos-eq6}
\end{align}
the system \eqref{yw:sys-auto} is GAOS. 
\end{theorem}

Observe that in the special case when $W=V$, this result is consistent with the implication of $iii)$ $\Rightarrow$ $i)$ in Theorem \ref{yw:theorem-gaos}. Nevertheless, the above result does not require that $W=V$: it just needs to be lower-bounded in terms of the output norm, as imposed by \eqref{yw:gaos-eq05} (in particular, $W(\phi)$ could be the output norm itself). Theorem \ref{yw:thm2-gaos} thus states that GAOS does hold if $W$ does not increase along the system's solutions, as imposed by \eqref{yw:gaos-eq6}. This additional flexibility with respect to Theorem \ref{yw:theorem-gaos} finds applications in adaptive control \cite{KACH20}.

\section{Input-to-state stability}\label{iasson:sec-iss}

In this section we consider systems of the form \eqref{pp:RFDE} to evaluate the impact of the input $u$ on the system's behavior. Standing Assumption 1 is assumed to hold, although many of the results presented here can be proved under less demanding regularity assumptions for the map $f$. For example, some results below hold for mappings $f$ which satisfy a right-hand side Lipschitz condition on bounded sets of $\mathcal X^n$: see \cite[Assumption (S1)]{Karafyllis:2008hc}.

\subsection{Definition and equivalent formulations}\label{sec:ISS:def}

The first definition of ISS for delay systems of the form \eqref{pp:RFDE} was given in \cite{Teel98}. The standard definitions of LISS and ISS are given below.

\begin{definition}[ISS, LISS]\label{def_Iasson_ISS}
System \eqref{pp:RFDE} is said to be \emph{input-to-state stable (ISS)} if there exist $\beta \in \mathcal {KL}$ and $\mu \in \mathcal N$ such that, for all $x_0\in \mathcal X^n$ and all $u\in \mathcal U^m$, 
\begin{equation} \label{Iasson__3_}
\left|x(t,x_0,u)\right|\le \beta \left(\left\| x_0 \right\| ,t\right)+
\mu \left(\left\| u_{[0,t]} \right\| \right),\quad \forall t\geq 0.
\end{equation}
 The function $\mu$ is then called an \emph{ISS gain}. We say that system \eqref{pp:RFDE} is \emph{locally input-to-state stable (LISS)} if \eqref{Iasson__3_} holds only for $\|x_0\|\leq r$ and $\|u\|\leq r$ for some $r>0$.
\end{definition}

As in the finite-dimensional case \cite{SON1,SONTAG2000}, the estimate \eqref{Iasson__3_} expresses the asymptotically vanishing effect of the initial conditions (via the term $\beta \left(\left\| x_0 \right\| ,t\right)$) and the persistent effect of the external input $u$ (via the term $\mu \left(\left\| u_{[0,t]} \right\| \right)$). 

It is worth stressing that the ISS gain $\mu$ is here considered to be of class $\mathcal N$, namely to be continuous, non-decreasing and zero at zero. The ISS literature often imposes that $\mu\in\cK$ or even $\cK_\infty$. The key property is that the ISS gain vanishes continuously at zero, which ensures that inputs of small magnitude produce small steady-state errors. The advantage of considering $\mu\in\mathcal N$ is that it allows for ISS gains that are identically zero, which proves useful in particular to assess robust versions of GAS and to unify results on feedback and cascade interconnections (see Section \ref{sec:interconnection}). Note however that any function of class $\mathcal N$ can be upper bounded by a $\cK$ or $\cK_\infty$ function, thus making the state estimate \eqref{Iasson__3_} also satisfied with these classes of ISS gains. Moreover, as in the finite-dimensional case, ISS automatically guarantees the following properties: 
\begin{itemize}
\item \emph{0-GAS}: the input-free system $\dot x(t)=f(x_t,0)$ is GAS
\item \emph{bounded input-bounded state (BIBS)} property:  
given any $x_0\in\cX^n$ and any $u\in\cU^m$,
\begin{align}\label{eq-antoine-BIBS}
    \|u\|<+\infty \quad \Rightarrow\quad \sup_{t \ge 0 } \left|x(t,x_0,u)\right| <+\infty
\end{align}
\item \emph{asymptotic gain (AG)} property: there exists $\mu\in\mathcal N$ such that, for all $x_0\in\cX^n$ and all $u\in\cU^m$:
\begin{align}\label{eq-antoine-AG}
\limsup_{t\to+\infty}|x(t, x_0, u)|\le \mu(\|u\|)
\end{align}
\item \emph{converging input-converging state (CICS)} property: given any $x_0\in\mathcal X^n$ and any $u\in\mathcal U^m$, 
\begin{align}\label{eq-antoine-CICS}
    \lim_{t\to +\infty } |u(t)|=0\quad \Rightarrow\quad \lim_{t\to +\infty } \left|x(t,x_0,u)\right|=0.
\end{align}
\end{itemize}

\begin{figure}[h!]
    \centering
    \includegraphics[width=0.8\textwidth]{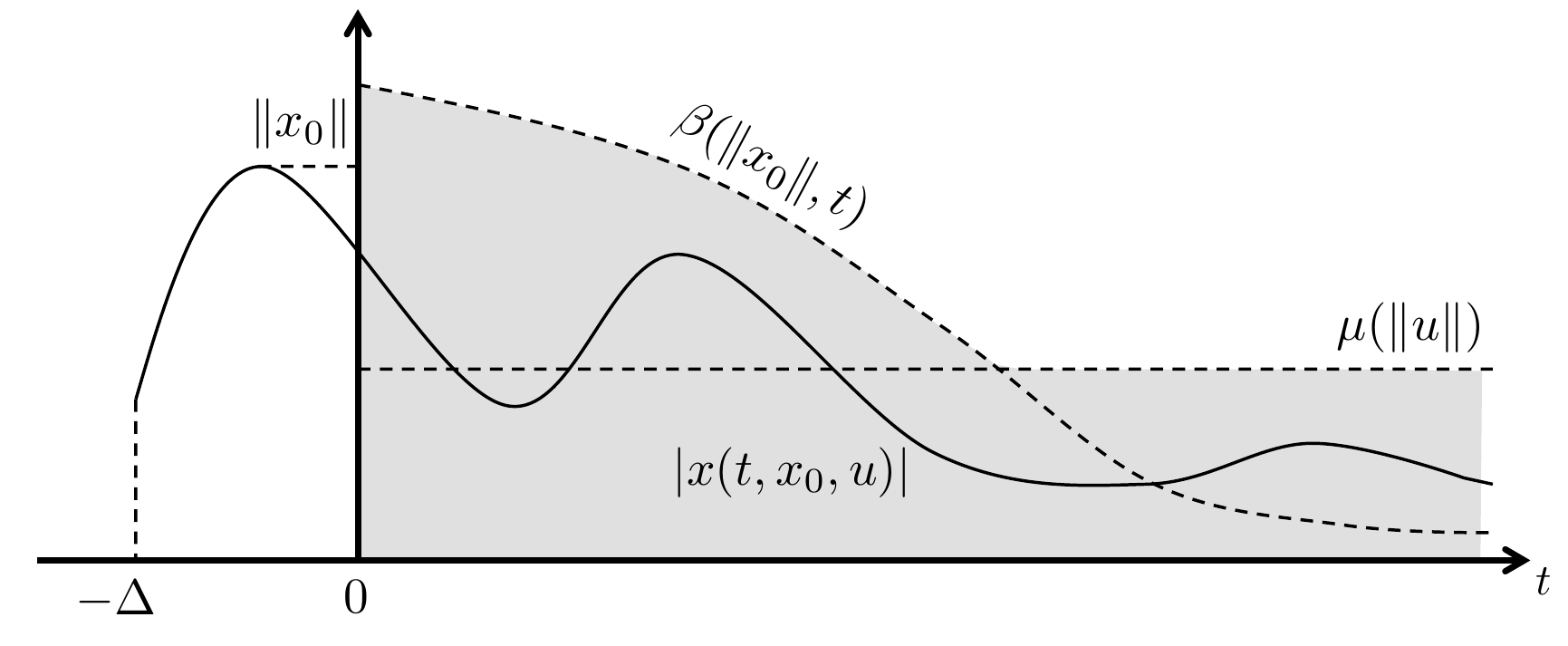}
    \caption{Schematic representation of the time evolution of an ISS system's solution.}
    \label{fig-2}
\end{figure}

 It should be noticed that the notion of ISS can be formulated in various equivalent ways. An example is the so-called ``max-formulation'', which replaces \eqref{Iasson__3_} by
\begin{align*}
 |x(t,x_0,u)|\leq \max \left\{\beta \left(\left\| \phi \right\| ,t\right),\mu \left(\left\| u_{[0,t]} \right\| \right)\right\},\quad \forall t\geq 0.
\end{align*}
This alternative formulation clearly illustrates the combination of transient behavior (mostly dictated by the initial state) and the steady-state one (only dictated by the input), as depicted by Figure \ref{fig-2}.

The following alternative formulation proves particularly useful for time-delay systems.

\begin{proposition}[History norm formulation] System \eqref{pp:RFDE} is ISS if and only if there exist $\beta \in \mathcal {KL}$ and $\mu \in \mathcal N$ such that, for all $x_0\in\mathcal X^n$ and all $u\in\mathcal U^m$,
\begin{equation*} %\label{Iasson__4_}
\left\| x_{t}(x_0,u) \right\| \le \beta \left(\left\| x_0 \right\| ,t\right)+ \mu \left(\left\| u_{[0,t]} \right\| \right),\quad \forall t\geq 0.
\end{equation*}
\end{proposition}

The proof of this result is straightforward, by taking the $\sup$ of \eqref{Iasson__3_} over the interval $[t-\Delta,t]$. The following formulation of the ISS property is less immediate and allows to focus on a more restricted class of initial states and inputs.

\begin{proposition}[Restricting initial states and inputs]\label{prop:dense} System \eqref{pp:RFDE} is ISS if and only if there exist $\beta \in \mathcal {KL}$ and $\mu \in \mathcal N$ such that estimate \eqref{Iasson__3_} holds for all initial states $x_0 \in C^{1} \left([-\Delta ,0],\mR^{n} \right)$ and all inputs $u\in \mathcal U^m$ which are piecewise-constant and right-continuous on $\mathbb R_{\ge 0}$.
\end{proposition}

The fact that we can restrict the considered initial states to continuously differentiable segments was already anticipated by Lemma \ref{prop:dense}. This result also shows that we can focus on piece-wise constant and right-continuous inputs with no loss of generality. The above formulation is particularly useful when studying the time evolution of functionals along the solutions of \eqref{pp:RFDE}. More specifically, given a functional $V: \mathcal X^n\to \mathbb R $ which is locally Lipschitz (or Lipschitz on bounded sets), for every initial condition $x_0 \in C^{1} \left([-\Delta ,0],\mathbb R^{n} \right)$ and every input $u\in \mathcal U^m$, the mapping $t\mapsto V(x_{t}(\phi,u) )$ is locally absolutely continuous for all $t\ge 0$ for which the solution of \eqref{pp:RFDE} is defined, thus allowing an easier integration of its Dini derivative.

The equivalence shown in Proposition \ref{prop:dense} was firstly formulated in \cite{HAIPEP21} in the context of switching retarded systems. See in particular Theorem 3.1 in that reference. 
%
% , where ISS of switching retarded systems with respect to % locally essentially bounded Lebesgue measurable input and % switching signals is proved to be equivalent to ISS with 
% respect to piecewise constant, right-continuous inputs. 
%
The equivalence for the case of piecewise continuous right-continuous inputs was established in \cite{Karafyllis:2008hc}.

The next formulation uses a different way of expressing the effect of the external input. 

\begin{proposition}[Fading memory estimate]\label{prop:fading} System \eqref{pp:RFDE} is ISS if and only if there exist $\beta\in\mathcal{KL}$ and $\mu\in\mathcal N$ such that, for all $x_0\in\mathcal X^n$ and all $u\in \mathcal U^m$,
\begin{equation}\label{iasson-iss-eq1}
|x(t, x_0, u)|\le\max\left\{\beta(\|x_0\|, t)\,, \, \sup_{\tau\in[0,t]}
\beta\bigl(\mu(|u(\tau)|), \ t-\tau\bigr)\right\},\quad \forall t\geq 0.
\end{equation} 
\end{proposition}

Estimate \eqref{iasson-iss-eq1} illustrates the fading memory effect: the solution's behavior is more influenced by recent values of the input than past ones. Fading memory ISS estimates were first used in \cite{PRAWAN} for finite-dimensional systems and are particularly useful for the derivation of small-gain conditions to analyze stability of interconnected systems (Section \ref{sec:interconnection}). In an infinite-dimensional context, fading memory estimates were recently utilized in the book \cite{karafyllis2019input} for the case where $\beta (s,t)=ks e^{-\lambda t}$ and $\mu(s)=\mu_0 s$, for some $k,\lambda,\mu_0>0$ (a property referred to as \emph{exp-ISS}). See also the discussion in \cite{KAJIbook11}. The proof of Proposition \ref{prop:fading} relies on the LKF characterizations of ISS (see Theorem \ref{theo_antoine_ISS_char} below) and can be found in \cite{Karafyllis:2008hc}.

The LISS property happens to be less interesting than expected, as it simply boils down to internal asymptotic stability, as proved in \cite[Theorem 6]{PEPAPAGA17}. 

\begin{theorem}[LISS $\Leftrightarrow$ 0-AS]\label{theo:LISS} System \eqref{pp:RFDE} is LISS if and only if the origin of the zero-input system $\dot x(t)=f(x_t,0)$ is AS.
\end{theorem}

This result shows the intrinsic robustness of asymptotic stability. It extends to time-delay systems the original result \cite[Lemma I.1]{SONWANTAC} established for finite-dimensional systems. An extension to more general classes of infinite-dimensional systems was also provided in \cite{mironchenko2016local}. Theorem \ref{theo:LISS} is a purely local result, both in the state and the input, and does not extend to global properties even in finite dimension \cite{SOKR03,TEHE04}. Indeed, 0-GAS (or even 0-GES) does not necessarily imply ISS, as illustrated by the following example.  

\begin{example}[0-GES $\nRightarrow$ ISS]
The system $\dot{x}(t)=-x(t)+x(t-1)u(t)$ admits unbounded solutions. For instance, for $x_0(\tau)=e^\tau$ for all $\tau\in[-1,0]$ and for the constant input $u\equiv 2e$, the solution reads $x(t)=e^t$ for all $t\geq 0$. Consequently, it does not satisfy the BIBS property and thus cannot be ISS. Nevertheless, the corresponding input-free system $\dot x(t)=-x(t)$ is GES. %Notice also that, for $x_0\equiv 0$, the solution is identically zero, regardless of the applied input.
\end{example}

There is anyway a class of systems for which ISS can be deduced from the stability properties of the corresponding input-free system. This is the class of globally Lipschitz systems. More precisely, we have the following
from \cite[Theorem 3.2]{Yeganefar:2008it}.

\begin{theorem}[0-GES $\Rightarrow$ ISS for globally Lipschitz systems]\label{theo_GES_Lipschitz} Suppose that there exist constants $L>0$, $p\in [0,1)$ such that, for all $\phi,\psi\in\mathcal X^n$ and all $v\in\mR^m$,
\begin{align} 
\left|f(\phi ,0)-f(\psi ,0)\right|&\le L\left\| \phi -\psi \right\| \label{GL__1_}\\
\left|f(\phi ,v)-f(\phi ,0)\right|&\le L\max \left\{\left\| \phi \right\| ^{p} ,1\right\}\left|v\right|. \label{GL__2_}
\end{align}
Then, provided that the origin of the input-free system $\dot{x}(t)=f(x_{t} ,0)$ is GES, the system \eqref{pp:RFDE} is ISS.
\end{theorem}

\begin{example}[Neuronal population, continued]
Let us go back to the neural population model considered in Example \ref{exa_neuro_1}, namely:
\begin{align}\label{eq_antoine_14bis}
    \dot x(t)=-x(t)+g\big(x(t-\Delta)+u(t)\big),
\end{align}
where $g:\mR\to\mR$ denotes a globally Lipschitz function with Lipschitz constant $\ell>0$ and satisfying $g(0)=0$. Note that, as compared to \eqref{eq_antoine_14}, we have added an input $u\in\mathcal U$, which could for instance model neuronal activity from other brain structures. We have seen in Example \ref{exa_neuro_1} that the origin of this system is GES if $\ell<1$, regardless of the value of the delay $\Delta\geq 0$. The vector field for this system reads $f(\phi,v)=-\phi(0)+g\big(\phi(-\Delta)+v\big)$ for all $\phi\in \mathcal X$ and all $v\in\mR$. Using the Lipschitz property of $g$, it holds for all $\phi,\psi\in\mathcal X$ that
\begin{align*}
    \left|f(\phi ,0)-f(\psi ,0)\right|&= \left|-\phi(0)+g\big(\phi(-\Delta)\big)+\psi(0)-g\big(\psi(-\Delta)\big)\right|\\
    &\leq |\phi(0)-\psi(0)|+\ell|\phi(-\Delta)-\psi(-\Delta)|\\
    &\leq (\ell+1)\|\phi-\psi\|,
\end{align*}
thus making \eqref{GL__1_} fulfilled with $L=\ell+1$. In addition, given any $\phi\in\mathcal X$ and any $v\in \mR$,
\begin{align*}
    \left|f(\phi,v)-f(\phi,0)\right| &= \left|-\phi(0)+g\big(\phi(-\Delta)+v\big)+\phi(0)-g\big(\phi(-\Delta)\big)\right|\\
    &\leq \ell \left|\phi(-\Delta)+v-\phi(-\Delta)\right|\\
    &\leq \ell |v|,
\end{align*}
thus establishing \eqref{GL__2_} with $p=0$ and the constant $L$ chosen above. We conclude with Theorem \ref{theo_GES_Lipschitz} that the system  \eqref{eq_antoine_14bis} is ISS if $\ell<1$, no matter the value of the delay $\Delta$.
\end{example}

Theorem \ref{theo_GES_Lipschitz} covers in particular linear time-invariant systems of the form
\begin{equation} \label{LTI__3_}
\dot{x}(t)=Ax_{t} +Bu(t)
\end{equation}
where $A$ is a bounded linear operator from $\mathcal X^n$ to $\mathbb R^{n}$ and $B\in\mR^{n\times m}$. Combining the above result with Theorem \ref{pp:LTI}, we obtain the following statement (see \cite{Pepe:2006ju}).

\begin{corollary}[AS $\Rightarrow$ ISS for linear systems] If the origin of the linear system $\dot x(t)=Ax_t$ is AS, then \eqref{LTI__3_} is ISS.
\end{corollary}

\subsection{Lyapunov-like conditions}

The notion of ISS Lyapunov function, that plays a crucial role for the characterization of the ISS property in the finite-dimensional case \cite{SONWANSCL}, can be extended to time-delay systems. However, in the infinite-dimensional case the situation is more complicated as an LKF may be requested to have different types of dissipation properties.

\begin{definition}[ISS LKF]\label{def_iasson_ISS_LKF}
For system \eqref{pp:RFDE}, an LKF $V:\mathcal X^n\to\mRp$ is said to be
\begin{itemize}
\item an \emph{ISS LKF with history-wise dissipation} if there exist $\alpha \in\cK_\infty$ and $\gamma \in \mathcal N$ such that
\begin{align} \label{Iasson__9_}
D^{+} V\left(\phi ,f(\phi,v)\right)\le -\alpha\left(\left\| \phi \right\| \right)+\gamma \left(\left|v\right|\right)
\end{align}
\item an \emph{ISS LKF with LKF-wise dissipation} if there exist $\alpha \in\cK_\infty$ and $\gamma \in \mathcal N$ such that
\begin{align} \label{Iasson__8_}
D^{+} V\left(\phi ,f(\phi,v)\right)\le -\alpha\left(V(\phi )\right)+\gamma \left(\left|v\right|\right)
\end{align}
\item an \emph{ISS LKF with point-wise dissipation} if there exist $\alpha \in\cK_\infty$ and $\gamma \in \mathcal N$ such that
\begin{align} \label{Iasson__7_}
D^{+} V\left(\phi ,f(\phi,v)\right)\le -\alpha\left(\left|\phi (0)\right|\right)+\gamma \left(\left|v\right|\right)
\end{align}
\item an \emph{ISS LKF in implication form} if there exist $\alpha \in \mathcal P$ and $\chi\in\mathcal N$ such that
 \begin{align} \label{Iasson__10_}
V\left(\phi \right)\ge \chi \left(\left|v\right|\right)\quad \Rightarrow\quad D^+V(\phi,f(\phi,v))\le -\alpha\left(\left|\phi (0)\right|\right)
\end{align}
where \eqref{Iasson__9_}-\eqref{Iasson__10_} are all meant to hold for all $\phi\in\mathcal X^n$ and all $v\in\mathbb R^m$.
\end{itemize}
\end{definition}

It is clear that the following implications hold:
\begin{align*}
\eqref{Iasson__9_} \quad \Rightarrow \quad  &\eqref{Iasson__8_} \quad \Rightarrow \quad \eqref{Iasson__7_} \\ \eqref{Iasson__8_} \quad &\Rightarrow \quad \eqref{Iasson__10_}.
\end{align*}

In a finite-dimensional context, the distinction between all these notions of ISS Lyapunov functions is irrelevant. This is due to the fact that, for delay-free systems, any Lyapunov function candidate is both upper and lower bounded by functions of the state norm. This fact ensures that a dissipation rate involving the whole Lyapunov function is qualitatively equivalent to a dissipation rate involving the state norm. In other words, all the above concepts are equivalent for finite-dimensional systems. However, this is not the case for delay systems. The following example shows that the difference in the dissipation rate appears even in the simplest example of a linear scalar system with a single discrete delay. 

\begin{example}[Different types of dissipation] \label{exa-dissipations}
\noindent Consider the scalar linear delay system
\begin{equation} \label{Ex__1_} 
\dot{x}(t)=-cx(t)+x(t-\Delta )+u(t) 
\end{equation} 
where %$x(t)\in {\mathbb R}$, $u$ and
$c=\left(3+e^{2\Delta}\right)/2$. Consider the LKFs defined for all $\phi \in \mathcal X$ as
\begin{align} 
V_{1} (\phi )&:=\frac{1}{2}\max_{\tau \in [-\Delta ,0]} \phi(\tau )^2e^{2\tau} \label{Ex__2_} \\
V_{2} (\phi )&:=\frac{1}{2} \phi(0)^2+\frac{e^{2\Delta}}{2} \int _{-\Delta }^{0}e^{2\tau}\phi(\tau )^2d\tau    \label{Ex__3_} \\
V_{3} (\phi )&:=\frac{1}{2} \phi(0)^2+\frac{1}{2} \int _{-\Delta }^{0}\phi(\tau )^2d\tau. \label{Ex__4_}  
\end{align} 
By virtue of Lemma \ref{means}, $V_1$ is Lipschitz on bounded sets. So are the functionals $V_{2}$ and $V_{3}$, as seen in Example \ref{exa1}. The functionals also satisfy the following inequalities for all $\phi \in \mathcal X$:
\begin{align} 
\frac{e^{-2\Delta}}{2} \left\| \phi \right\| ^{2} &\le V_{1} (\phi )\le \frac{1}{2} \left\| \phi \right\| ^{2}  \label{Ex__5_} \\
\frac{1}{2} \left|\phi (0)\right|^{2} &\le V_{2} (\phi )\le \frac{1+\Delta e^{2\Delta}}{2} \left\| \phi \right\| ^{2}  \label{Ex__6_} \\ 
\frac{1}{2} \left|\phi (0)\right|^{2} &\le V_{3} (\phi )\le \frac{1+\Delta }{2} \left\| \phi \right\| ^{2}. \label{Ex__7_}  
\end{align} 
In other words, all of them are LKFs (and $V_1$ is actually a coercive LKF). Using Lemma \ref{means}, it holds that
\begin{align*}
    &D^{+} V_{1} (\phi ,f(\phi ,v))\le \\ &\left\{\begin{array}{cl}-2V_{1} (\phi ),\quad &\textrm{if } V_1(\phi)> \phi(0)^2/2\\
    \max \left\{-2V_{1} (\phi ),-c\phi(0)^2+\phi (0)\phi (-\Delta )+\phi (0)v\right\},\quad &\textrm{if } V_{1} (\phi )=\phi (0)^2/2,
    \end{array}\right.
    \end{align*}
 where $f$ denotes the functional defining the right-hand side of \eqref{Ex__1_}. It follows in the latter case that
\begin{align*} %\label{Ex__8_} 
 D^{+} V_{1} (\phi ,f(\phi ,v))&\le \max \left\{-2V_{1} (\phi ),-c\phi(0)^2+\phi (0)\phi (-\Delta )+\phi (0)v\right\} \\ 
 &\le \max \left\{-2V_{1} (\phi ),-\left(c-1\right)\phi(0)^2+\frac{1}{2} \phi (-\Delta )^2+\frac{1}{2} v^{2} \right\} \\ 
 &\le \max \left\{-2V_{1} (\phi ),-2\left(c-1\right)V_{1} (\phi )+e^{2\Delta} V_{1} (\phi )+\frac{1}{2} v^{2} \right\}\\
 &\leq-V_{1} (\phi )+\frac{1}{2} v^{2}.
\end{align*} 
We conclude that the following estimate holds for all $\phi \in \mathcal X$ and all $v\in {\mathbb R}$:
\begin{equation} \label{Ex__9_} 
D^{+} V_{1} (\phi ,f(\phi ,v))\le -V_{1} (\phi )+\frac{1}{2} v^{2},
\end{equation} 
meaning that $V_1$ is a coercive ISS LKF with LKF-wise dissipation. 
It then follows from \eqref{Ex__5_} and \eqref{Ex__9_} that 
\begin{align}
    D^{+} V_{1} (\phi ,f(\phi ,v))\le -\frac{e^{-2\Delta}}{2} \left\| \phi \right\| ^{2} +\frac{1}{2} v^{2},
\end{align}
meaning that $V_{1} $ is also a coercive ISS LKF with history-wise dissipation. 

On the other hand, observing that $\phi (0)\phi (-\Delta )\le \frac{1}{2} \phi(0)^2+\frac{1}{2} \phi(-\Delta )^2$ and $\phi (0)v\le \frac{1}{2} \phi(0)^2+\frac{1}{2} v^{2} $, $V_2$ satisfies for all $\phi \in \mathcal X$ and all $v\in {\mathbb R}$:
\begin{align*} %\label{Ex__10_} 
D^{+} V_{2} (\phi ,f(\phi ,v))&\le -\frac{1}{2} \phi(0)^2-e^{2\Delta}\int _{-\Delta }^{0}e^{2\tau}\phi(\tau )^2d\tau  +\frac{1}{2} v^{2} \\
&\le -V_{2} (\phi )+\frac{1}{2} v^{2}. 
\end{align*} 
In other words, $V_{2} $ is an ISS LKF with LKF-wise dissipation. However notice that $V_{2} $ is not an ISS LKF with history-wise dissipation, as the negative terms appearing in its Driver's derivative cannot be upper-bounded by a term of the form $-\alpha(\|\phi\|)$ with $\alpha\in\cK_\infty$.

Finally, observing that $\phi (0)\phi (-\Delta )\le \frac{1}{2} \phi(0)^2+\frac{1}{2} \phi (-\Delta )^2$ and $\phi (0)v\le \frac{1}{2} \phi (0)^2+\frac{1}{2} v^{2} $, it holds for all $\phi \in \mathcal X$ and all $v\in {\mathbb R}$ that
\begin{align*} %\label{Ex__11_} 
D^{+} V_{3} (\phi ,f(\phi ,v))\le -\frac{e^{2\Delta}}{2} \phi(0)^2+\frac{1}{2} v^{2},
\end{align*} 
meaning that $V_{3} $ is an ISS LKF with point-wise dissipation. However, notice that $V_{3} $ has neither a history-wise nor a  LKF-wise dissipation, as the only negative term appearing in its Driver's derivative involves $|\phi(0)|^2$, which cannot be lower-bounded by any $\cK_\infty$ function of $V_3(\phi)$ or $\|\phi\|$.
\end{example}

The following characterization of ISS is a combination of \cite[Theorem 3.3]{Karafyllis:2008hc} and \cite[Theorem 2]{Kankanamalage:2017ug}.

\begin{theorem}[LKF characterizations of ISS]\label{theo_antoine_ISS_char} The following properties are equivalent:
\begin{itemize}
\item[$i)$] \eqref{pp:RFDE} is ISS
\item[$ii)$] \eqref{pp:RFDE} admits a coercive ISS LKF with history-wise dissipation
\item[$iii)$] \eqref{pp:RFDE} admits an ISS LKF with LKF-wise dissipation
\item[$iv)$] \eqref{pp:RFDE} admits an ISS LKF in implication form.
\end{itemize}
\end{theorem}

Note that, the LKF in Item $ii)$ being coercive, its history-wise dissipation is equivalent to an LKF-wise dissipation. It should also be noted that in \cite{Kankanamalage:2017ug} the function $\alpha \in \mathcal P$ appearing in the implication form \eqref{Iasson__10_} (Item $iv)$) is required to be of class $\mathcal \cK_\infty$ but the result can be generalized to cover the case where $\alpha\in\mathcal P$.

The definition of LKF adopted in this survey requests a $\cK_\infty$ lower bound in terms of $|\phi(0)|$. As stated in \cite[Theorem 4.4]{jacob2020noncoercive}, ISS can actually be demonstrated under an even weaker requirement for systems satisfying the RFC property, namely that $V$ vanishes only when $\phi$ is identically zero.

\begin{theorem}[ISS through non-coercive LKF]Assume that \eqref{pp:RFDE} is RFC and that there exist a functional $V: \mathcal X^n\to \mathbb R_{\ge 0}$ which is Lipschitz on bounded sets of $\mathcal X^n$, $\alpha, \overline\alpha \in \mathcal K_{\infty }$ and $\gamma \in \mathcal N$ such that, for all $\phi\in\mathcal X^n$ and all $v\in\mR^m$,
\begin{align} 
0<V(\phi )&\le \overline\alpha \left(\left\| \phi \right\| \right),\quad \forall \phi\neq 0
\label{Iasson__16_}\\
 D^{+} V\left(\phi ,f(\phi,v) \right) &\le -\alpha\left(\left\| \phi \right\| \right)+\gamma \left( \left| v \right| \right). \label{Iasson__17_}
\end{align}
Then the system \eqref{pp:RFDE} is ISS.
\end{theorem}

The above result was stated and proved in \cite{jacob2020noncoercive} in an actually wider infinite-dimensional context and considering a continuous functional $V: \mathcal X^n\to \mathbb R_{\ge 0}$ (which required the replacement of the Driver's derivative in \eqref{Iasson__17_} by a Dini derivative along the solutions of \eqref{pp:RFDE}). 

Among the different types of dissipation proposed in Definition \ref{def_iasson_ISS_LKF}, the point-wise one is the least demanding. As discussed in Section \ref{sec:open:pointwiseISS}, it is not known yet whether or not the existence of an ISS LKF with point-wise dissipation is enough to guarantee ISS. The only results in that direction require an additional assumption on the growth of the LKF or on the vector field itself. The following result is a slight generalization of the corresponding result in \cite{CHPEMACH17}.

\begin{theorem}[ISS through point-wise dissipation rate]\label{theo_antoine_ISS_pointwise2}
Consider an LKF $V:\mathcal X^n\to \mRp$ for which there exist $\underline\alpha,\overline\alpha_1,\overline\alpha_2,\overline\alpha_3\in\cK_\infty$ such that, for all $\phi\in\mathcal X^n$,
\begin{align}\label{eq_antoine_1}
\underline\alpha(|\phi(0)|)\leq V(\phi)\leq \overline\alpha_1(|\phi(0)|)+\overline\alpha_2\left(\int_{-\Delta}^0 \overline\alpha_3(|\phi(\tau)|)d\tau\right).
\end{align}
Assume that $V$ is an ISS LKF with point-wise dissipation for \eqref{pp:RFDE}, meaning satisfying \eqref{Iasson__7_} for some $\alpha \in \cK_\infty$ and $\gamma \in \mathcal N$. Then, under the condition that
\begin{align}\label{eq_antoine_2}
\liminf_{s\to +\infty} \frac{\alpha(s)}{\overline\alpha_3(s)}>0,
\end{align}
the system \eqref{pp:RFDE} is ISS.
\end{theorem}

Per se, \eqref{eq_antoine_1} constitutes a mild requirement as several typical LKFs used in practice do satisfy such a bound (see for instance Example \ref{exa1}). The main requirement in the above statement is \eqref{eq_antoine_2}, which imposes that the term $\overline\alpha_3$ under the integral sign in the bound \eqref{eq_antoine_1} is dominated at infinity by the dissipation rate $\alpha$. It should be noted that growth rate conditions on the mapping $f$ have also been proposed in \cite{CHPEMACH17} under which a point-wise ISS dissipation is enough to guarantee ISS. When focusing on quadratic LKFs, the above result boils down to the following.

\begin{corollary}[ISS through quadratic ISS LKF with point-wise dissipation]\label{corol_antoine_corol1}
Assume that there exist symmetric positive definite matrices $P_1,P_2,Q\in\mR^{n\times n}$ and $\gamma \in \mathcal N$ such that the functional defined as
\begin{align*}
V(\phi):=\phi(0)^{\top}P_1\,\phi(0)+\int_{-\Delta}^0 \phi(\tau)^\top P_2\,\phi(\tau)d\tau,\quad \forall \phi\in\mathcal X^n,
\end{align*}
satisfies, for all $\phi\in\mathcal X^n$ and all $v\in\mR^m$,
\begin{align*}
D^+V(\phi,f(\phi,v))\leq -\phi(0)^\top Q\,\phi(0)+\gamma(|v|).
\end{align*}
Then the system \eqref{pp:RFDE} is ISS.
\end{corollary}

Another method to establish ISS is through the so-called Lyapunov-Razumikhin approach. This alternative to the Lyapunov-Krasovskii approach owns the advantage to rely on functions rather than functionals, which sometimes simplifies the analysis. The following result originally appeared in \cite{Teel98}, but can also be found in \cite[Proposition 4.1]{Karafyllis:2008hc}.

\begin{theorem}[ISS through Razumikhin approach]\label{thm-iss-raz}
Assume that there exist $\rho\in \mathcal K_{\infty }$, $\alpha\in\mathcal P$, $\gamma \in \mathcal N$, and a positive definite and radially unbounded function $V\in C^{1} (\mathbb R^{n} ;\mathbb R_{\ge 0} )$ satisfying, for all $\phi\in\mathcal X^n$ and all $v\in\mR^m$,
\begin{align}\label{Iasson__13_}
V(\phi (0))\ge \max\left\{\rho\left(\max_{\tau\in [-\Delta ,0]} V(\phi (\tau))\right)\,,\,\gamma \left(\left|v\right|\right)\right\}\quad \Rightarrow \quad\nabla V(\phi (0))f(\phi ,v)\le -\alpha(|\phi (0)|).
\end{align}
If $\rho(s)<s$ for all $s>0$, then the system \eqref{pp:RFDE} is ISS.
\end{theorem}

Theorem \ref{thm-iss-raz} can be interpreted as a small-gain theorem for systems in the following form:

\begin{equation}\label{yw:sys-sg}
\dot x(t) = F(x(t), x_t, u(t)),
\end{equation}
where $F:\R^n\times\cX^n\times\R^m\rightarrow\R^n$ is Lipschitz on bounded sets of $\R^n\times\cX^n\times\R^m$ with $F(0, 0, 0)=0$. For such a system, $x_t$ can be treated as an external input by considering the auxiliary system
\begin{equation}\label{yw:sys-sgw}
\dot x(t) = F(x(t), w_t, u(t)),
\end{equation}
where $w_t\in\mathcal X^n$. Note then that \eqref{yw:sys-sg} is the closed-loop system of \eqref{yw:sys-sgw} under the feedback law $w_t=x_t$. Assume that there exist $\rho\in \mathcal K_{\infty }$, $\alpha\in\mathcal P$, $\gamma \in \mathcal N$, and a positive definite and radially unbounded function $V\in C^{1} (\mathbb R^{n} ;\mathbb R_{\ge 0} )$ satisfying
\begin{align}\label{yw:sg-dis1}
V(x)\ge \max\left\{\rho\left(\max_{\tau\in [-\Delta ,0]} V(\phi(\tau)) \right) , \gamma \left(\left|v\right|\right)\right\} \quad \Rightarrow \quad
\nabla V(x)f(x, \phi, v)\le -\alpha(|x|)
\end{align}
for all $(x, \phi, v)\in\R^n\times\cX^n\times\R^m$. Then \eqref{Iasson__13_} is automatically satisfied. Indeed, \eqref{yw:sg-dis1} is requested to hold for all $(x, \phi, v)\in\R^n\times\cX^n\times\R^m$; while \eqref{Iasson__13_} is only required when $x=\phi(0)$. A motivation for considering \eqref{yw:sys-sgw} is that it converts a delay system as in \eqref{yw:sys-sg} into a delay-free system with $w$ treated as a disturbance. In applications, it can be more difficult to design feedback laws based on the history of the trajectories of 
\eqref{yw:sys-sg} than for \eqref{yw:sys-sgw} where $w$ is treated as disturbance or uncertainty. Such an approach was used, for instance, in \cite{jankovic2001control} and \cite{Teel98}.

It should be stressed that in many of the results that were presented in this section, the proofs provide explicit formulas for the estimation of the ISS gain $\mu$. This is important for two reasons. First, we can directly estimate the sensitivity of the system with respect to external disturbances. Second, we can apply the derived ISS estimate by using a small-gain result in order to study the stability properties of a more complicated system (see Section \ref{sec:ISS:smallgain}).

\subsection{Solutions-based conditions}

Many solutions-based conditions for ISS were proposed in the finite-dimensional case \cite{SONWANTAC}. Such conditions enabled researchers to obtain less demanding conditions for the derivation of ISS estimates and to deeply understand the nature and link between different stability notions. Some of the solutions-based conditions that were used in the finite-dimensional case were extended to cover time-delay systems. In particular, the following result provides ISS characterizations under the RFC assumption. Its proof can be found in \cite{Mironchenko:2017wb}.

\begin{theorem}[Solutions-based characterizations of ISS]\label{theo-sol-char-ISS}
Assume that \eqref{pp:RFDE} is RFC. Then the following statements are equivalent:
\begin{enumerate}
\item[$i)$] \eqref{pp:RFDE} is ISS;
\item[$ii)$] \eqref{pp:RFDE} satisfies the \emph{uniform asymptotic gain property (UAG)}, that is,
there exists some $\mu \in \mathcal N$ such that for every $\varepsilon , r>0$ there exists $T >0$ such that, for all $u\in \mathcal U^m$, % with $\|u\|<+\infty$,
\begin{align}\label{eq-antoine-UAG}
\sup \left\{\, \left| x(t,x_0 ,u)\right| \, :\, x_0 \in \mathcal X^n,\
 \left\| \phi \right\| \le r, \ t\ge T \right\}\le \varepsilon +\mu \left(\|u\|\right)
\end{align}
\item[$iii)$] \eqref{pp:RFDE} satisfies the following two properties:
\begin{itemize}
\item
its origin is \emph{uniformly locally stable (ULS)}, meaning that there exist $\sigma \in \cK_\infty$, $\mu \in \mathcal N$ and $r>0$ such that, for all $x_0 \in \mathcal X^n$ and all $u\in \mathcal U^m$ with $\left\| x_0 \right\| \le r$ and $\|u\|\le r$,
\begin{align}\label{eq-antoine-ULS}
|x(t;x_0,u)| \le \sigma \left(\left\| x_0 \right\| \right)+\mu \left(\|u\|\right),\quad \forall t\geq 0
\end{align}
\item
it owns the \emph{uniform limit (ULIM)} property, meaning that there exists $\mu\in\mathcal N$ such that, for every $\varepsilon, r>0$, there exists $T >0$ such that, for all $x_0 \in \mathcal X^n$ and all $u\in \mathcal U^m$ with $\left\| \phi \right\| \le r$,
\begin{align}\label{eq_antoine_19}
\min_{t\in[0,T]} \left\| x_{t} (x_0 ,u)\right\|\le \varepsilon +\mu \left(\|u\|\right).
\end{align}
\end{itemize}
\end{enumerate}
\end{theorem}

This result is reminiscent of its finite-dimensional counterpart presented in \cite{SONWANTAC}, where it is shown that, for delay-free systems, ISS is equivalent to UAG, meaning that any solution eventually converges to a neighborhood of the origin whose size in ``proportional'' to $\|u\|$ and that the rate of this convergence is uniform on bounded sets of initial states and inputs, as requested by Item $ii)$. In that reference, it is also shown that, for finite-dimensional systems, ISS is equivalent to stability of the origin of the input-free system plus the limit property, which essentially requires that any solution eventually visits a neighborhood of the origin whose size is small when the input is small. Here, the requirements in Item $iii)$ are more demanding, as ULS imposes a stability notion for the system with inputs and as the ULIM property requests that the time needed to visit the considered neighborhood of the origin is uniform on bounded sets of initial states and inputs.

It is also worth stressing that the ULIM property imposes that the whole history segment eventually approaches the $\mu(\|u\|)$-neighborhood of the origin. To date, it is not known whether ISS would also hold under the combination of ULS and the point-wise version of ULIM, in which \eqref{eq_antoine_19} would be replaced by
\begin{align*}
\min_{t\in[0,T]} \left| x(t,x_0 ,u)\right|\le \varepsilon +\mu \left(\|u\|\right).
\end{align*}
See Section \ref{sec:open:pointwiseISS} for more discussion on that matter.

It is well-known that the ISS property in the finite-dimensional case is equivalent to an integral-to-integral estimate \cite[Theorem 1]{IISS}. The same is true for the time-delay case but for a more demanding integral-to-integral estimate. From Theorem \ref{theo_antoine_ISS_char}, we know that if \eqref{pp:RFDE} is ISS then there exists a coercive LKF $V: \mathcal X^n\to \mathbb R_{\ge 0}$ with LKF-wise dissipation. In particular, for all $x_0 \in C^1([-\Delta,0],\mR^n)$ and all $u\in \mathcal U^m$, the map $t\mapsto V(x_t(x_0,u))$ is locally absolutely continuous and satisfies, for almost all $t \geq 0$,
\begin{equation*} %\label{DV__2_}
\frac{d}{dt}V(x_t(x_0,u))\le -\alpha (\|x_t(x_0,u)\|)+\gamma \left(\left|u(t)\right|\right),
\end{equation*}
for some $\alpha\in\cK_\infty$ and some $\gamma\in\mathcal N$. By integration and considering the bounds provided in \eqref{pp:Coer__3_}, we obtain the following integral-to-integral estimate:
\begin{equation} \label{Iasson__25_}
\int_{0}^{t}\alpha\left(\left\| x_{\tau} (x_0 ,u)\right\| \right)d\tau \le \overline\alpha \left(\left\| x_0 \right\| \right)+\int _{0}^{t}\gamma \left(\left|u(\tau)\right|\right)d\tau   ,\quad \forall t\geq 0.
\end{equation}
Approximating an arbitrary initial condition $x_0 \in \mathcal X^n$ by an initial condition in $C^1([-\Delta,0],\mR^n)$, we can establish that estimate \eqref{Iasson__25_} actually holds for all $x_0 \in \mathcal X^n$ and all $u\in \mathcal U^m$, which in turn implies the following norm-to-integral estimate:
\begin{equation} \label{Iasson__26_}
\int _{0}^{t}\alpha\left(\left\| x_{\tau} (x_0 ,u)\right\| \right)d\tau \le \overline\alpha \left(\left\| x_0\right\| \right)+ \gamma \left(\left\| u \right\| \right) t  ,\quad \forall t\geq 0.
\end{equation}
It can be shown that the integral-to-integral estimate \eqref{Iasson__25_} and the norm-to-integral estimate \eqref{Iasson__26_} are equivalent to ISS under the assumption of RFC. Indeed, combining \cite[Theorem 3.5]{jacob2020noncoercive} and \cite[Lemma 3.3]{Mironchenko:2017wb}, we have the following.

\begin{theorem}[ISS through integral estimates]
Assuming that \eqref{pp:RFDE} is RFC, the following statements are equivalent:
\begin{itemize}
\item \eqref{pp:RFDE} is ISS
\item there exist $\overline\alpha \in \cK_\infty$, $\alpha \in \mathcal K$ and $\gamma \in \mathcal N$ such that \eqref{Iasson__25_} holds for all $x_0 \in \mathcal X^n$ and all $u\in \mathcal U^m$
\item there exist $\overline\alpha \in \cK_\infty$, $\alpha \in \mathcal K $ and $\gamma \in \mathcal N$ such that \eqref{Iasson__26_} holds for all $x_0 \in \mathcal X^n$ and all $u\in \mathcal U^m$.
\end{itemize}
\end{theorem}

The derivation of the norm-to-integral estimate \eqref{Iasson__26_} can be made by means of a non-coercive LKF for \eqref{pp:RFDE} with a history-wise dissipation rate: see \cite[Proposition 4.3]{jacob2020noncoercive}. Here again, it is not known whether a similar result would hold under a point-wise version of these estimates: see Section \ref{sec:open:pointwiseISS} for more discussions about this.

\section{Integral input-to-state stability}\label{sec-iISS}

\subsection{Definition}
The ISS property (Definition \ref{def_Iasson_ISS}) relates the solutions' norm to the magnitude of the applied input. One may rather assess the impact of the input energy on the solutions. This is captured by the following property.

\begin{definition}[iISS]\label{def_antoine_iISS}
The system \eqref{pp:RFDE} is said to be \emph{integral input-to-state stable (iISS)} if there exists $\beta\in\mathcal{KL}$ and $\eta,\mu\in\mathcal N$ such that, for all $x_0\in\mathcal X^n$ and all $u\in\mathcal U^m$, its solution satisfies
\begin{align}\label{eq_antoine_iISS}
|x(t,x_0,u)|\leq\beta(\|x_0\|,t)+\eta\left(\int_0^t \mu(|u(\tau)|)d\tau\right),\quad \forall t\geq 0.
\end{align}
$\mu$ is then referred to as an \emph{iISS gain}.
\end{definition}

By virtue of Lemma \ref{lem-dense}, it is actually sufficient to establish the estimate \eqref{eq_antoine_iISS} for all $x_0\in C^1([-\Delta,0];\mR^n)$ and all $u\in\mathcal U^m$ to conclude iISS, in which case $t\mapsto x_t(x_0,u)$ is locally absolutely continuous (Theorem \ref{pp:lemmahalepeperfde}).

This property is the natural adaptation of the iISS property originally defined for delay-free systems \cite{IISS}. Its first use in a time-delay context goes back to \cite{Pepe:2006ju}. It is worth stressing that $\eta$ and $\mu$ are usually imposed to be of class $\cK_\infty$ (rather than class $\mathcal N$). Like for ISS, there is no qualitative difference between the two formulations, but allowing these functions to be of class $\mathcal N$ turns out to be handy in practice, as it allows to cover the case of robust GAS (when either $\eta$ or $\mu$ is identically zero) and to derive stability results on cascades using small-gain arguments (see Section \ref{sec:interconnection}).

Clearly, iISS ensures that the system is 0-GAS. Unlike ISS, iISS does not guarantee that solutions are bounded in response to bounded inputs. Nevertheless, iISS does ensure a robustness property known as \emph{bounded-energy converging-state (BECS)}. Namely, if solutions of \eqref{pp:RFDE} satisfy the iISS estimate \eqref{eq_antoine_iISS}, then it holds that
\begin{align}\label{eq_antoine_BECS}
\int_0^{+\infty} \mu(|u(\tau)|)d\tau<+\infty\quad \Rightarrow\quad \lim_{t\to+\infty} |x(t,x_0,u)|=0.
\end{align}
In other words, any input with bounded energy (as measured through its iISS gain) generates solutions that converge to the origin.

\subsection{Lyapunov-like conditions}

Here also, depending on the way they dissipate along the solutions, we may consider different types of LKF.

\begin{definition}[iISS LKF]\label{def_antoine_iISS_LKF}
For the time-delay system \eqref{pp:RFDE}, an LKF $V:\mathcal X^n\to \mRp$ is:
\begin{itemize}
\item an \emph{iISS LKF with history-wise dissipation rate} if there exist $\alpha\in\mathcal{P}$ and $\gamma\in\mathcal N$ such that
\begin{align}\label{eq_antoine_103}
D^+V(\phi,f(\phi,v))\leq -\alpha(\|\phi\|)+\gamma(|v|),
\end{align}
\item an \emph{iISS LKF with LKF-wise dissipation rate} if there exist $\alpha\in\mathcal{P}$ and $\gamma\in\mathcal N$ such that
\begin{align}\label{eq_antoine_102}
D^+V(\phi,f(\phi,v))\leq -\alpha(V(\phi))+\gamma(|v|),
\end{align}
\item an \emph{iISS LKF with point-wise dissipation rate} if there exist $\alpha\in\mathcal{P}$ and $\gamma\in\mathcal N$ such that
\begin{align}\label{eq_antoine_101}
D^+V(\phi,f(\phi,v))\leq -\alpha(|\phi(0)|)+\gamma(|v|),
\end{align}
\item an \emph{iISS LKF with $\cKL$ dissipation rate} if there exist $\sigma\in\mathcal{KL}$ and $\gamma\in\mathcal N$ such that
\begin{align}\label{eq_antoine_104}
D^+V(\phi,f(\phi,v))\leq -\sigma(|\phi(0)|,\|\phi\|)+\gamma(|v|),
\end{align}
\end{itemize}
where \eqref{eq_antoine_101}-\eqref{eq_antoine_104} are all meant to hold for all $\phi\in\mathcal X^n$ and all $v\in\mR^m$. $\alpha$ and $\sigma$ are then referred to as a \emph{dissipation rate} whereas $\gamma$ is called a \emph{supply rate}.
\end{definition}

A key difference with the ISS LKFs (Definition \ref{def_iasson_ISS_LKF}) is that, in \eqref{eq_antoine_103}-\eqref{eq_antoine_101}, the dissipation rate is requested to be
a class $\mathcal P$ function (rather than a $\cK_\infty$ one). In particular, $\alpha$ is here allowed to be bounded or even to vanish at infinity. Invoking \cite[Lemma IV.1]{ANG-bigIISS}, which states that, given $\alpha\in\mathcal{P}$, there exist $\mu\in\cK_\infty$ and $\ell\in\mathcal L$ such that $\alpha(s)\geq \mu(s)\ell(s)$ for all $s\geq 0$, it can easily be checked that
\begin{align*}
\eqref{eq_antoine_102}\quad\Rightarrow\quad \eqref{eq_antoine_101}\quad \Rightarrow\quad \eqref{eq_antoine_104}.
\end{align*}
In other words, any iISS LKF with LKF-wise dissipation is also an iISS LKF with point-wise dissipation, which itself constitutes a subclass of iISS LKFs with $\cKL$ dissipation. In the same way, \eqref{eq_antoine_103} readily implies \eqref{eq_antoine_104}, whereas \eqref{eq_antoine_103} and \eqref{eq_antoine_102} turn out to be equivalent if $V$ is coercive. Similar to their ISS counterparts, the distinction between all these notions of iISS Lyapunov functions is irrelevant in a finite-dimensional context. For time-delay systems, it is usually convenient to work with a non-coercive LKF, so the lower bound on $V$ usually involves merely the instantaneous value of the solution's norm $|\phi(0)|$, thus impeding to jungle easily between a point-wise dissipation, an LKF-wise one, and a $\cKL$ one.

The following result, established in \cite{CHGOPE21}, demonstrates the tight links exisiting between these iISS LKF concepts.

\begin{theorem}[iISS LKF characterizations]\label{theo_antoine_iISS_char}
For the system \eqref{pp:RFDE}, the following statements are equivalent:
\begin{itemize}
\item[$i)$] it admits a coercive iISS LKF with history-wise dissipation
\item[$ii)$] it admits an iISS LKF with LKF-wise dissipation
\item[$iii)$] it admits an iISS LKF with history-wise dissipation
\item[$iv)$] it admits an iISS LKF with $\cKL$ dissipation
\item[$v)$] it is iISS.
\end{itemize}
Moreover, if \eqref{pp:RFDE} admits an iISS LKF with point-wise dissipation, then it is iISS.
\end{theorem}

An important aspect of this result is that iISS can be established using an LKF with point-wise dissipation or even with a $\cKL$ dissipation rate (which turns out to be particularly useful when using an LKF of the form $V=\ln(1+W)$, where $W$ denotes an other LKF: a classical trick used in iISS analysis). On the other hand, once iISS is established, the above result ensures the existence of a coercive iISS LKF with history-wise dissipation (or, equivalently, an LKF-wise dissipation), which can prove useful to conduct further robustness analysis.

\begin{example}[Establishing iISS with LKF approach] Consider the scalar time-delay system
\begin{align*}
\dot x(t) = -\frac{x(t)}{1+ x(t-\Delta)^2}+x(t)u(t).
\end{align*}
and let $f:\mathcal X\times \mR\to\mR$ be such that $f(x_t,u(t))$ denotes its right-hand side.
Let $V(\phi)=\ln(1+\phi(0)^2)$, then $V$ is Lipschitz on bounded sets and satisfies, for all $\phi\in\mathcal X$,
\begin{align*}
    \ln(1+|\phi(0)|^2)\leq V(\phi)\leq \ln(1+\|\phi\|^2),
\end{align*}
meaning that $V$ is an LKF. Moreover, for all $\phi\in\mathcal X$ and all $v\in\mR$, it holds that
\begin{align*}
D^+V(\phi, f(\phi, v))& = \frac{2\phi(0)}{1+\phi(0)^2}\left( -\frac{\phi(0)}{1+\phi(-\Delta)^2} +\phi(0)v\right)\\
&\le -\frac{2\phi(0)^2}{(1+\|\phi\|^2)^2} + 2|v|.
\end{align*}
In other words, $V$ is an iISS LKF with $\cKL$ dissipation. By Theorem \ref{theo_antoine_iISS_char}, the system is iISS, regardless of the value of the delay $\Delta\geq 0$. Notice that, for this system, it does not seem easy to find iISS LKFs with point-wise, LKF-wise or history-wise dissipation (although Theorem \ref{theo_antoine_iISS_char} does guarantee their existence). 
\end{example}

We may consider relaxing the requirement of having a negative term in the LKF dissipation inequalities of Definition \ref{def_antoine_iISS_LKF}: this leads to the property known as zero-output dissipativity in the literature of delay-free systems \cite{ANGSGIISS}.

\begin{definition}[Zero-output dissipativity]
The system \eqref{pp:RFDE} is said be \emph{zero-output dissipative} if there exist an LKF $V:\mathcal X^n\to\mRp$ and $\nu\in\mathcal N$ such that, for all $\phi\in\mathcal X^n$ and all $v\in\mR^m$,
\begin{align*}
D^+V(\phi,f(\phi,v))\leq \nu(|v|).
\end{align*}
It is said be \emph{coercively zero-output dissipative} if the above holds with a coercive LKF.
\end{definition}

The following result states that iISS can readily be concluded from 0-GAS provided that the system is zero-output dissipative.

\begin{theorem}[iISS and zero-output dissipativity]\label{theo_antoine_-iISS_ZOdiss}
For the system \eqref{pp:RFDE}, the following properties are equivalent:
\begin{itemize}
\item iISS
\item 0-GAS and coercive zero-output dissipativity
\item 0-GAS and zero-output dissipativity.
\end{itemize}
\end{theorem}

This result is reminiscent of its delay-free counterpart presented in \cite{ANGSGIISS}. Its proof, provided in \cite{CHGOPE21,LIWA18}, relies on the following characterization of 0-GAS, which may be of interest on its own.

\begin{proposition}[0-GAS characterizations]\label{prop_antoine_0GASchar}
The following statements are equivalent:
\begin{itemize}
\item the system \eqref{pp:RFDE} is 0-GAS
\item there exist a coercive LKF $V:\mathcal X^n\to\mRp$, a nondecreasing continuous function $\ell:\mRp\to\mRp$, and $\eta,\gamma\in\cK_\infty$ such that, for all $\phi\in\mathcal X^n$ and all $v\in\mR^m$,
\begin{align}\label{eq_antoine_39bis}
D^+V(\phi,f(\phi,v))\leq -&\eta(\|\phi\|)+\ell(\|\phi\|)\gamma(|v|).
\end{align}
\item there exist a coercive LKF $V:\mathcal X^n\to\mRp$, a continuously differentiable class $\cK$ function $\pi$ satisfying $\pi'(s)>0$ for all $s\geq 0$, $\alpha\in\mathcal{PD}$, and $\gamma\in\cK_\infty$, such that $W:=\pi\circ V$ satisfies, for all $\phi\in\mathcal X^n$ and all $v\in\mR^m$,
\begin{align}
D^+W(\phi,f(\phi,v))\leq -\alpha(\|\phi\|)+\gamma(|v|).\label{eq_antoine_39}
\end{align}
\end{itemize}
\end{proposition}

The finite-dimensional counterpart of this result was originally given in \cite[Lemma IV.10]{ANG-bigIISS}. Due to coercivity, $\|\phi\|$ can equivalently be replaced by $V(\phi)$ in both \eqref{eq_antoine_39bis} and \eqref{eq_antoine_39}. It is worth stressing that the functional $W$ in \eqref{eq_antoine_39} might not be radially unbounded (as $\pi$ might not be a $\cK_\infty$ function), which is why \eqref{eq_antoine_39} does not readily guarantee iISS (otherwise, in view of Theorem \ref{theo_antoine_iISS_char}, any 0-GAS system would be iISS, which is untrue even in finite dimension \cite{ANG-bigIISS}).

For finite-dimensional systems, it is known from \cite{IISS} that internal asymptotic stability implies iISS (and actually the stronger property known as Strong iISS, which combines both iISS and ISS with respect to small inputs \cite{CHANIT13}). A similar result can be derived for time-delay systems.

\begin{proposition}[iISS for bilinear systems]\label{prop-bilinear}
Consider the time-delay system
\begin{align}\label{eq_antoine_111}
    \dot x(t) = Ax(t)+\sum_{i=1}^m\sum_{j=1}^p u_i(t)A_{ij}x(t-\Delta_j)+Bu(t),
\end{align}
where $B\in\mR^{n\times m}$, $\Delta_j\in[0,\Delta]$, and $A,A_{ij}\in\mR^{n\times n}$ for all $i\in\{1,\ldots m\}$ and all $j\in\{1,\ldots,p\}$. Then \eqref{eq_antoine_111} is iISS if and only if $A$ is Hurwitz.
\end{proposition}

This result was established in \cite{Pepe:2006ju}
 for the case of a single delay. Although it is stated here for discrete delays, the above proposition holds also true for distributed delays \cite{CHGOPE21}. More generally, this result remains valid in a wider infinite-dimensional context \cite{Mironchenko:2016kucdc}. In particular, the matrix $A$ can be replaced by a linear bounded operator from $\mathcal X^n$ to $\mR^n$.

%We close this section by providing a Lyapunov-Razumikhin condition %for iISS, see \cite[Theorems 1 \& 2]{TIWAJI12}.
%

As in the ISS case, for systems in the form of \eqref{yw:sys-sg}, considering the delay-free system \eqref{yw:sys-sgw} can be helpful in applying the Razumihkin method for the iISS property as follows.

\begin{theorem}[iISS through Lyapunov-Razumikhin]\label{yw:razum-iISS}
Assume that there exist a positive definite and radially unbounded function $V\in C^1(\R^n; \R_{\ge 0})$, $\alpha\in\mathcal P$, $\gamma\in\mathcal N$, and $ \rho\in\cK_\infty$ such that, for all $\phi\in\mathcal X^n$, $w\in\cX^n$, and all $v\in\mR^m$,
\[
V(\phi(0))\ge \rho\left(\max_{\tau\in[-\Delta, 0]}V(w(\tau))\right)\quad 
\Rightarrow \quad \nabla V(\phi(0))F(x, w,v)\le -\alpha(|\phi(0)|) + \gamma(|v|).
\]
Suppose there exists some $\kappa\in\cK_\infty$ so that $\rho(s) + \kappa\circ\rho(s)\le s$ for all $s\geq 0$. Then the system \eqref{yw:sys-sg} is iISS.
\end{theorem}

The following example illustrates the application of this result.

\begin{example}[Establishing iISS through Razumikhin]
Consider the scalar system
\begin{equation}\label{yw:ex-iISS-Raz}
\dot x(t) = -\arctan\left(\frac{x(t)}{1+ |x(t-\Delta)|} \right) + u(t).
\end{equation}
By treating the term $x(t-\Delta)$ as an uncertainty $w(t)$,
the system is converted to 
\begin{equation}\label{yw:ex-iISSw}
\dot x(t) = -\arctan\left(\frac{x(t)}{1+ |w(t)|}\right) + u(t),
\end{equation}
which is a delay-free system with $(w, u)$ as inputs. In other words, the vector field in Theorem \ref{yw:razum-iISS} can be picked as 
\begin{align*}
    F(x,\omega,v):=-\arctan\left(\frac{x}{1+ |\omega|}\right) + v,\quad \forall x,\omega,v\in\mR.
\end{align*}
The derivative of $V(x) :=\ln(1+x^2)$ along the solutions of this system reads, for all $x,\omega,v\in\mR$,
\begin{align*}
\nabla V(x)F(x, \omega, v) &= \frac{2x}{1+x^2}\left( -\arctan\left(\frac{x}{1+|\omega|}\right) + v\right)\\
&\leq -\frac{2|x|}{1+x^2}\arctan\left(\frac{|x|}{1+|\omega|}\right) + 2|v|.
\end{align*}
In particular, the following implication holds:
\[
|x|\ge \frac{|\omega|}{2} \quad \Rightarrow\quad
\nabla V(x)F(x, \omega, v)\le -\frac{2|x|}{1+x^2}\arctan\left(\frac{|x|}{1+2|x|}\right) + 2|v|.
\]
%It follows that, for the original system \eqref{yw:ex-iISS-Raz}, the following holds for all $\phi\in\mathcal X$ and all $v\in\R$:
%\[
%V(\phi(0))\ge \frac12\,\max_{\tau\in[-\Delta, 0]}V(\phi(\tau))
%\quad \Rightarrow\quad \nabla V(\phi(0))f(\phi, v)\le
%-\frac{2|\phi(0)|}{1+\phi(0)^2}\arctan\left(\frac{|\phi(0)|}{1+2|\phi(0)|}\right) + 2|v|.
%\]
By Theorem \ref{yw:razum-iISS}, with $\rho(s) = s/2$, $\gamma(s)=2s$, $\kappa(s)=s$, and $\alpha(s)=\frac{2s}{1+s^2}\arctan\left(\frac{s}{1+2s}\right)$ for all $s\geq 0$, the system \eqref{yw:ex-iISS-Raz} is iISS, independently of the value of the delay $\Delta$.
\end{example}

%\YW{Will be fine with me to remove the next example}
Observe that Theorem \ref{yw:razum-iISS} implies that if the system \eqref{yw:sys-sgw} is iISS in $u$ and ISS in $w$ with an ISS gain function $\chi$ satisfying the small-gain condition $\chi(s) + \kappa\circ\chi(s) \le s$, then the system \eqref{yw:sys-sg} is iISS.  However,
it should be noted that though treating $x_t$ as a disturbance can be helpful in identifying Lyapunov functions for the Razumikhin approach, it may induce some conservatism, as can be seen through the following example.  

\begin{example}[Limitations of Razumikhin approach]
Consider the scalar system
\begin{align}\label{eq-antoine-112}
\dot x(t) = -x(t) + x(t-\Delta)u(t).
\end{align}
This system being bilinear, iISS follows readily from Proposition \ref{prop-bilinear}. However, the system 
\begin{align}\label{eq-antoine-112n}
\dot x(t) = -x(t) + w(t-\Delta)u(t)
\end{align}
is not ISS in $w$ with a gain function $\chi(s)\le s$ and iISS in $u$, as requested by Theorem \ref{yw:razum-iISS}.  Indeed, if this were the case, then we would have 
\[
|x(t)|\le \max\left\{
\beta(|x(0|)|, t), \ \chi(\|w\|), \ \gamma_1\left(\int_0^t\gamma_2(|u(\tau)|)d\tau\right)\right\},
\]
for some $\beta\in\cKL$ and some $\chi, \gamma_1, \gamma_2\in\cK_\infty$ with $\chi(s)\le s$. Consider the constant inputs $w(t)=\bar w$ and $u(t)=\bar u$ with $\bar w,\bar u>0$. Then it would hold that
\[
|x(t)|\le \max\{\beta(|x(0|)|, t), \ \chi(\overline w), \ \gamma_1(\gamma_2(\overline u)t)\}.
\]
Fix any instant $T>0$.  Then it would hold that, given any $\bar u>0$,
\begin{equation}\label{e-T}
\limsup_{\bar w\to+\infty}\frac{|x(T)|}{\bar w}\le \frac{\chi(\bar w)}{\bar w}\le 1.
\end{equation}
However, the trajectories of \eqref{eq-antoine-112n} are given by $x(t) = x(0)e^{-t} + \bar w \bar u(1-e^{-t})$. In particular, for $T=\ln(2)$ and $\bar u=4/\ln(2)$, we would obtain that
\[
\limsup_{\bar w\to +\infty} \frac{|x(T)|}{\bar w} = \bar u T(1-e^{-T})=2,
\]
which contradicts \eqref{e-T}.
\end{example}

%\begin{example}[Limitations of Razumikhin approach]
%Consider the scalar system
%\begin{align}\label{eq-antoine-112}
%\dot x(t) = -x(t) + x(t-\Delta)u(t).
%\end{align}
%This system being bilinear, iISS follows readily from Proposition \ref{prop-bilinear}. However, it does not seem easy to derive this conclusion from Theorem \ref{yw:razum-iISS}. Indeed, a natural Lyapunov function for this system would be $V(x)=x^2/2$, whose derivative satisfies, for all $\phi\in\mathcal X$ and all $v\in \mR$,
%\begin{align*}
%   \nabla V(\phi(0))f(\phi,v) = \phi(0)\big(\phi(0)+\phi(-\Delta)v\big),
%\end{align*}
%When treating $x(t-\Delta)$ as an input,
%the system is not ISS in $w$ with a gain function $\gamma_w$ satisfying %$\gamma_w(r) < r$. 
%\end{example}

\subsection{Solutions-based conditions}

We have already mentioned that iISS guarantees that solutions converge to the origin in response to any input with bounded energy (as measured through the iISS gain of the system). We may weaken this requirement by merely requesting that, under this bounded input energy assumption, solutions remain uniformly bounded. More precisely, we may consider the following property, originally introduced in  \cite{ANGSGIISS} for delay-free systems.

\begin{definition}[UBEBS]\label{def_antoine_UBEBS}
The system \eqref{pp:RFDE} is said to have the \emph{uniform bounded energy-bounded state (UBEBS) property} if there exist $\alpha,\xi,\zeta\in\mathcal K_\infty$ and $c\geq 0$ such that, for all $x_0\in\mathcal X^n$ and all $u\in\mathcal U^m$, its solution satisfies
\begin{align}\label{eq_antoine_UBEBS}
\alpha(|x(t,x_0,u)|)\leq \xi(\|x_0\|)+\int_0^t\zeta(|u(\tau)|)d\tau+c,\quad \forall t\geq 0.
\end{align}
\end{definition}

In \cite{ANGSGIISS} it was shown that, for finite-dimensional systems, ISS is equivalent to 0-GAS combined with UBEBS. The following result, established in \cite{CHGOPE21}, extends this characterization to time-delay systems.

\begin{theorem}[iISS and UBEBS]\label{theo-iISS2}
For the system \eqref{pp:RFDE}, the following properties are equivalent:
\begin{itemize}
\item iISS
\item 0-GAS and UBEBS property \eqref{eq_antoine_UBEBS} with $c=0$
\item 0-GAS and UBEBS property \eqref{eq_antoine_UBEBS}.
\end{itemize}
\end{theorem}

It is worth noting that the following implication holds:
\begin{align}\label{eq-antoine-UBEBS-Zod}
    \textrm{zero-output dissipativity} \quad \Rightarrow\quad \textrm{UBEBS } \eqref{eq_antoine_UBEBS} \textrm{ with } c=0,
\end{align}
thus providing a convenient LKF-based way to establish UBEBS in practice.

Other solutions-based characterizations of iISS that hold true in finite-dimension have not yet been extended to time-delay systems: see in particular the discussion in Section \ref{sec-open-sol-iISS}.

\section{Input-to-output stability}\label{yw: sec-ios}

We have seen in Section \ref{yw:sec-autonomous} that it is sometimes useful to impose stability properties on only part of the state variables or, more generally, on a specific output of the system. In this section, we provide tools to assess robustness of such output stability properties in the presence of exogenous inputs. To that aim, consider a system as in (\ref{pp:RFDE}) with an output map:
\begin{subequations}\label{yw:syso}
\begin{align}
\dot x(t)&= f(x_t,u(t)),\\ 
y(t) &= h(x_t),
\end{align}
\end{subequations}
where $y(t)\in\R^p$ represents the output of interest. The output map $h: \mathcal X^n\rightarrow\R^p$ is assumed to be Lipschitz on bounded sets with
$h(0)=0$.  %The vector field $f:\mathcal X^n\times\mR^m\to\mR^n$ is also assumed to be Lipschitz on bounded sets. 
%\YW{This should be removed, as in Section 3.3: For the reasons exposed in Section \ref{yw:sec-autonomous}, it is not assumed that $f(0,0)=0$, so Standing Assumption \ref{stand_ass_1} is not requested to hold in this section.} 
For each $x_0\in\mathcal X^n$ and $u\in\mathcal U^m$, we use $y(\cdot, x_0, u)$ to denote the output function given by
$y(t, x_0, u) = h(x_t(x_0, u))$. We assume throughout this section that the system is forward complete, in the sense of Definition \ref{def_FC}, which can be established for instance using Theorem \ref{theo-RFC-1} or \ref{theo-RFC-2}.
The following property will also be needed for several results of this section.

\begin{definition}[UGS]\label{yw:UGS}
We say that the system (\ref{yw:syso}) is {\it uniformly globally stable (UGS)} if there exists $\sigma\in\mathcal K_\infty$ such that, for all $x_0\in \mathcal X^n$ and all $u\in\mathcal U^m$,
\[
|x(t, x_0, u)|\le\max\{\sigma(\|x_0\|), \ \sigma(\|u_{[0,t]}\|)\}, \quad\forall\,t\ge 0.
\]
\end{definition}

This property imposes not only a bounded state in response to any bounded input (the property referred as BIBS in Section \ref{sec:ISS:def}), but also that this bound can be picked uniform on any bounded sets of initials states and inputs. Moreover, it imposes that the solutions norm remains arbitrarily small at all times if the initial state and the input are taken sufficiently small. In particular, the UGS property also implies both RFC and ULS. It is worth mentioning that this property was referred to as \emph{uniform bounded input-bouded state (UBIBS)} in \cite{Kankanamalage:2017ug,kankanamalage2019remarks}.

\subsection{Notions of input-to-output stability}\label{yw:sec-ios-def}

In this section we consider several notions related to input-to-output stability and report the corresponding LKF results to establish them.

\subsubsection{Definitions}\label{sec-IOS-def}

\begin{definition}[IOS, OL-OLIOS, SI-IOS]\label{yw:def-ios}
The forward complete system (\ref{yw:syso}) is said to be
\begin{itemize}
\item
{\it input-to-output stable (IOS)} if there exist $\beta\in \cKL$ and $\mu\in\mathcal N$
 such that, for all $x_0\in\mathcal X^n$ and all $u\in\mathcal U^m$,
\begin{equation}\label{yw:ios}
|y(t, x_0, u)|\le \beta(\|x_0\|, t) + \mu(\|u_{[0, t]}\|), \quad\forall\,t\ge 0
\end{equation}
\item
{\it output-Lagrange input-to-output stable (OL-IOS)} if it is IOS and, additionally, there exists $\sigma\in\mathcal K_\infty$ such that, for all $x_0\in\mathcal X^n$ and all $u\in\mathcal U^m$,
\begin{equation}\label{yw:olios}
|y(t, x_0, u)|\le
\max\left\{\sigma(|h(x_0)|), \ \sigma(\|u_{[0, t]}\|)\right\},\quad \forall t\geq 0
\end{equation}
\item
{\it state-independent input-to-output stable (SI-IOS)} if there exist $\beta\in\cKL$ and $\mu\in\mathcal N$ such that, for all $x_0\in\mathcal X^n$ and all $u\in\mathcal U^m$,
\begin{equation}\label{yw:siios}
|y(t, x_0, u)|\le \beta(|h(x_0)|, \ t) + \mu(\|u_{[0, t]}\|),\quad \forall t\geq 0.
\end{equation}
\end{itemize}
\end{definition}

Here again, by Lemma \ref{lem-dense}, the set of considered initial states can be restricted to $C^1([-\Delta,0],\mR^n)$ with no loss of generality. The difference between these three notions lies in how the initial state may affect the transient overshoots and the decay rate of the output. For the SI-IOS property, both the overshoots and the decay rates are dictated by the magnitude of the initial output value, while for the OL-IOS property the decay rate can be slowed down by large initial state variables that are not present in the output variables. For the IOS property, both the overshoots and the decay rate of the output variable can be affected by all initial state variables, especially when the magnitude of the initial state is much larger than that of the initial output.

As in the delay-free case \cite{SONWANIOS},
the  OL-IOS property is equivalent to the existence of some $\beta\in \cKL$, $\kappa\in \mathcal K$, and $\mu\in\mathcal N$ such that, for all $x_0\in\mathcal X^n$ and all $u\in\mathcal U^m$,
\begin{equation}\label{yw:olgaos2}
|y(t, x_0, u)|\le \beta\left(|h(x_0)|, \ \frac{t}{1+\kappa(\|x_0\|)}\right) + \mu(\|u_{[0, t]}\|), \quad \forall t\geq 0,
\end{equation}
which clearly illustrates that the decay rate may depend on the initial state norm, but the maximal transient overshoot depends only on the initial output. 

Clearly, IOS boils down to ISS when the output map is taken as $h(\phi)=\phi(0)$. It is also straightforward to see that the IOS property implies the following:
\begin{itemize}
\item \emph{uniform global output stability (UGOS)}:
there exists some $\sigma\in\cK_\infty$ such that, for all $x_0\in\cX^n$ and all $u\in\cU^m$,
\begin{align}\label{eq-antoine-BIBO}
|y(t, x_0, u)|\le \max\{\sigma(\|x_0\|), \ \sigma(\|u_{[0,t]}\|)\},\quad\forall\,t\ge0;
\end{align}
\item \emph{output asymptotic gain (OAG)}: %output-AG
there exists some $\mu\in\mathcal N$ such that, for all $x_0\in\cX^n$ and all $u\in\cU^m$,
\begin{align}\label{eq-OAG}
\limsup_{t\rightarrow+\infty}|y(t, x_0, u)|\le \mu(\|u\|).
\end{align}
\end{itemize}
Here again, when the output is the whole state, OAG boils down to the asymptotic gain property introduced in Section \ref{sec:ISS:def}, whereas UGOS boils down to UGS. When $u\equiv0$, the UGOS property becomes the GOS property introduced in Section \ref{yw:sec-autonomous}, whereas OAG corresponds to global output attractivity of the origin (possibly non-uniformly over bounded sets of initial states). From this, it can be seen that
the combination of UGOS and OAG does not imply the IOS property, not even for delay-free systems. 

A related question is whether the combination of the state estimate \eqref{yw:olios} (meaning output-Lagrange stability, OLS)  and OAG implies the OL-IOS property. 
For the finite dimensional systems,
this implication was shown to be true in \cite{SONWANTAC} for the special case of ISS and in \cite[Theorem 1]{Ingalls-ACC01} for the output case. This result also indicates that the OL property can lead to significant improvement on output stability. However, for delay systems, the question whether the combination of OL and OAG imply OL-IOS
still remains open even for the special case when $h(\phi) = \phi(0)$: see Section \ref{sec-conj-sol-IOS}. Nevertheless, the conjunction of UGOS and OAG is related to IOS. For instance, in \cite{POLILU06}, the IOS property was defined this way, and a small-gain theorem was developed for such a property for interconnected systems: see Theorem \ref{yw:sg-thm3} below. This approach was also used in \cite{TIWAJIlargescale12}.

\subsubsection{Lyapunov-like conditions}

As in the case of ISS, the Lyapunov method is indispensable in the study for output stability properties. In this section we adapt the Lyapunov-Krasovskii functionals used in Section \ref{iasson:sec-iss} for the output case.

A frequent application of IOS is the study of the ISS property with respect to a compact set $\Omega\subset \mR^n$ which contains the origin. In that case, the output map can be taken as
\[
h(\phi) := \max_{\tau\in[-\Delta,0]}{\rm dist}(\phi(\tau), \Omega),
\]
where $\textrm{dist}(\cdot\,,\Omega)$ denotes the distance to the set $\Omega$, as defined in Section \ref{sec-Notation}. For such applications, the output map $h$ satisfies the following property:
there exists $r\ge 0$ and $\rho\in\mathcal K_\infty$ such that, for all $\phi\in\mathcal X^n$,
\begin{equation}\label{yw:ots}
\|\phi\|\le \rho(|h(\phi)|) + r.
\end{equation}
Based on this observation, the following LKF characterization of IOS was obtained in \cite[Theorem 3.3]{Karafyllis:2008hc}. 

\begin{theorem}[LKF characterization of IOS under \eqref{yw:ots}]\label{theo-LKF-IOS}
Suppose that the system (\ref{yw:syso}) satisfies condition (\ref{yw:ots}).  Then the following statements are equivalent:
\begin{enumerate}
\item[$i)$] the system \eqref{yw:syso} is IOS
\item[$ii)$]
there exist a functional $V: \mathcal X^n\rightarrow\R_{\ge 0}$ which is Lipschitz on bounded sets,  $\underline\alpha,\overline\alpha\in\mathcal K_\infty$ and $\gamma\in\mathcal N$ such that, for all $\phi\in\mathcal X^n$ and all $v\in\R^m$,
\begin{align}
\underline\alpha(|h(\phi)|)&\le V(\phi)\le \overline\alpha(\|\phi\|), \label{yw:ios-eq2}\\
D^+V(\phi,f(\phi,v)) &\le -V(\phi) + \gamma(|v|) \nonumber %\label{yw:ios-eq3}
\end{align}
\item[$iii)$] there exists a functional $V: \mathcal X^n\rightarrow\R_{\ge 0}$ which is Lipschitz on bounded sets, $\underline\alpha,\overline\alpha\in\mathcal K_\infty$, $\chi\in\mathcal N$, and $\alpha\in \mathcal P$ such that, for all $\phi\in\mathcal X^n$ and all $v\in\R^m$, (\ref{yw:ios-eq2}) holds and
\begin{equation*}%\label{yw:ios-eq4}
V(\phi)\geq\chi(|v|)\quad \Rightarrow\quad D^+V(\phi,f(\phi,v)) \le -\alpha(V(\phi)).
\end{equation*}
\end{enumerate}
%Finally, if $h:\mathcal X^n\rightarrow\R^p$ is equivalent to a finite-dimensional continuous map $h_0:\R^n\rightarrow\R^p$, then inequality (\ref{yw:ios-eq2}) can be replaced by
%\begin{equation}\label{yw:ios-eq5}
%\underline\alpha(|h_0(\phi(0))|)\le V(\phi)\le \overline\alpha(\|\phi\|),\quad \forall \phi \in\mathcal X^n.
%\end{equation}
\end{theorem}

A consequence of this result is that, just like ISS (Proposition \ref{prop:fading}), the IOS property can be alternatively expressed as a state estimate with a fading memory of past inputs.

\begin{proposition}[Fading memory]
Suppose that the system (\ref{yw:syso}) satisfies condition (\ref{yw:ots}). Then it is IOS if and only if there exist $\beta\in\mathcal{KL}$ and $\mu\in\mathcal N$ such that, for all $x_0\in\mathcal X^n$,
\begin{equation*}%\label{yw-ios-eq1}
|y(t, x_0, u)|\le\max\left\{\beta(\|x_0\|, t), \ \sup_{\tau\in[0,t]}
\beta\bigl(\mu(|u(\tau)|), \ t-\tau\bigr)\right\},\quad \forall t\geq 0.
\end{equation*}
\end{proposition}

In \cite{Kankanamalage:2017ug}, the work on LKF characterization of IOS was developed along a different line: the UGS condition  (Definition \ref{yw:UGS}) was imposed instead of requiring condition (\ref{yw:ots}). In this statement, we rely on the following notation: for a continuous function $V:\cX^n\rightarrow\R$, $\textrm{Ker}(V):=\{\phi\in\cX^n\,:\, V(\phi) = 0\}$.

\begin{theorem}[LKF characterization of IOS under UGS]\label{yw:yw-ios-thm}
Assume that (\ref{yw:syso}) satisfies the UGS property. Then the following statements hold:
\begin{enumerate}
\item[$i)$] the system
(\ref{yw:syso}) is IOS if and only if there exist a functional $V:\mathcal X^n\rightarrow\R_{\ge 0}$ which is locally Lipschitz on $\mathcal X^n\setminus\textrm{Ker}(V)$, $\underline\alpha,\overline\alpha
    \in \cK_\infty$, $\chi\in\mathcal N$ 
     and $\sigma\in\cKL$ such that, for all $\phi\in\mathcal X^n$ all $v\in\R^m$,
\begin{align}
\underline\alpha(|h(\phi)|)&\le \,V(\phi)\le \overline\alpha(\|\phi\|), \label{yw:yw-lyap1}\\
V(\phi)\ge\chi(|v|)\quad \Rightarrow\quad &D^+V(\phi,f(\phi,v))\le -\sigma(V(\phi), \|\phi\|) \label{yw:yw-lyap2}
\end{align}
\item[$ii)$]
the system is OL-IOS if and only if
there exist a functional $V:\mathcal X^n\rightarrow\R_{\ge 0}$ that is locally Lipschitz on  $\mathcal X^n\setminus \textrm{Ker}(V)$, $\underline\alpha, \overline\alpha\in\mathcal K_\infty$, $\chi\in\mathcal N$, and $\sigma\in\cKL$ such that, for all $\phi\in\mathcal X^n$ and all $v\in\mR^m$, property (\ref{yw:yw-lyap2}) holds and
\begin{equation}\label{yw:yw-lyap-ol-1}
\underline\alpha(|h(\phi)|)\le V(\phi)\le \overline\alpha(|h(\phi)|).
\end{equation}
\end{enumerate}
\end{theorem}

The above characterization is stated for locally Lipschitz functional (everywhere they do not vanish). It is likely that a similar characterization holds for functionals that are Lipschitz on bounded sets, although this has not yet been formally established yet.

Note the difference between the bounds (\ref{yw:yw-lyap1}) and (\ref{yw:yw-lyap-ol-1}): while $V(\phi)$ is bounded above by a function of the state norm $\|\phi\|$ for the IOS case, it is upper-bounded by a function of the output norm $|h(\phi)|$ for the OL-IOS case.  This difference illustrates how a functional for OL-IOS guarantees that the output overshoot is dominated by magnitudes of the initial output variable and the input signals.
In \cite{SONWANIOS-LYA}, it was shown that for any $\sigma\in\mathcal{KL}$, there exist $\kappa_1, \kappa_2\in\mathcal K$ such that
\[
\sigma(r, s)\ge \frac{\kappa_1(r)}{1+\kappa_2(s)},\quad \forall r,s\geq 0.
\]
Hence, the implication \eqref{yw:yw-lyap2} can be restated as, for some $\chi\in\mathcal N$ and some $\kappa_1, \kappa_2\in\mathcal{K}$:
\[
V(\phi)\ge\chi(|v|)\quad \Rightarrow\quad D^+V(\phi,f(\phi,v))\le -\frac{\kappa_1(V(\phi))}{1+\kappa_2(\|\phi\|)}.
\]

The following Razumikhin-type result is a specialization of \cite[Proposition 4.1]{Karafyllis:2008hc}.

\begin{theorem}[Razumikhin conditions for IOS]
Consider system (\ref{yw:syso}) and suppose that $h:\mathcal X^n \to \mathbb R^{p} $ is equivalent to a finite-dimensional mapping $h_{0} :\mathbb R^{n} \to \mathbb R^{p} $. Assume that there exist a continuously differentiable function $V:\mR^n\to\mRp$, $\underline\alpha, \overline\alpha,\rho\in\cK_\infty$, $\alpha\in\mathcal P$, and $\gamma\in\mathcal N$ such that, for all $\phi \in \mathcal X^n$ and all $v \in \mathbb R^{m}$,
\begin{align*} %\label{Raz-bounds}
\underline\alpha (\left|h_{0} (\phi(0))\right|)\le V(\phi(0))&\le \overline\alpha (\left|\phi(0)\right|)\\
V(\phi (0))\ge \max\left\{\rho\left(\max_{\tau\in [-\Delta ,0]} V(\phi (\tau))\right)\,,\,\gamma \left(\left|v\right|\right)\right\}\quad &\Rightarrow\quad \nabla V(\phi (0))f(\phi ,v)\le -\alpha(V(\phi (0))). %\label{Raz_diff_inequality}
\end{align*}
Assume further that either \eqref{pp:RFDE} is RFC or there exist $\eta\in \mathcal K_{\infty }$ and $R\ge 0$ such that $\eta\left(\left|\phi(0)\right|\right)\le V(\phi(0))+R$ for all $\phi \in \mathcal X^n$. Then, provided that $\rho(s)<s$ for all $s>0$, the system \eqref{pp:RFDE} is IOS and RFC.
\end{theorem}

\subsection{Notions of integral input-to-output stability}\label{yw:section-iIOS}

In this section we consider the integral version of the properties associated with IOS in which, similarly to iISS, the impact of the input is measured through its energy rather than its magnitude. Some of the results may apply to the case when $h$ is a map from $\mathcal X^n$ to $\R^p$, but we will focus on the special case when $h$ is a map from $\R^n$ to $\R^p$ so that, along trajectories of the system, $ y(t) = h(x(t))$. We use $H:\cX^n\rightarrow\cX^p$ to denote the map defined by \begin{align*}
   H(\phi)(\tau)=(h\circ\phi)(\tau), \quad \forall \tau\in[-\Delta, 0].
\end{align*}

\subsubsection{Definitions}

The following notions were introduced in \cite{Nawarathna21}.

\begin{definition}[iIOS, OL-iIOS, SI-iIOS]\label{yw:def-iios}
The forward complete system (\ref{yw:syso}) is:
\begin{itemize}
\item {\it integral input-to-output stable} (iIOS) if there exist $\beta\in\mathcal{KL}$, and $\eta, \mu\in\mathcal N$ such that for all $x_0\in\mathcal X^n$, and all $u\in\mathcal U^m$, 
\begin{equation}\label{yw:iios}
|y(t, x_0, u)|\le\beta(\|x_0\|, t) + \eta\left(\int_0^t\mu(|u(\tau)|)\,d\tau\right),\quad\forall\,t\ge 0
\end{equation}
\item
{\it output-Lagrange iIOS}
(OL-iIOS) if it is iIOS and additionally, there exist $\sigma\in\mathcal K_\infty$ and $\eta,\mu\in\mathcal N$ such that, for all $x_0\in\mathcal X^n$ and all $u\in\mathcal U^m$, 
\begin{equation}\label{yw:e-oliios-1}
|y(t,x_0, u)|) \le
\sigma(\|H(x_0)\|) +
\eta\left(\int_0^t\sigma(|u(\tau)|)\,d\tau\right), \quad\forall\,t\ge 0
\end{equation}
\item {\it state-independent iIOS}\ (SI-iIOS)
if there exist $\beta\in\mathcal{KL}$ and $\eta, \mu\in\mathcal N$ such that, for all $x_0\in\mathcal X^n$, and all $u\in\mathcal U^m$, 
\begin{equation}\label{yw:e-siiios}
|y(t, x_0, u)|\le\beta(\|H(x_0)\|, t) + \eta\left(\int_0^t\mu(|u(\tau)|)\,d\tau\right),\quad\forall\,t\ge 0.
\end{equation}
\end{itemize}
\end{definition}

Similar to the IOS case, the OL-iIOS property (\ref{yw:e-oliios-1}) can be restated as: for some $\beta\in\mathcal{KL}$, $\alpha\in\mathcal K_\infty$, $\kappa\in\cK$, and $\mu\in\mathcal N$,
\begin{equation*}%\label{yw:e-olios3}
\alpha(|y(t,x_0,u)|)\le\beta\left(
\|H(x_0)\|, \ \frac{t}{1+\kappa(\|x_0\|)}\right) + \int_0^t\gamma(|u(\tau)|)\,d\tau,
\end{equation*}
which highlights the fact that, for the OL-iIOS property, the output decay rate may be slower for larger initial states.

In the special case when $h(z)=z$ for all $z\in\R^p$, that is, when $h(x(t)) = x(t)$, all the notions
of iIOS, OL-iIOS, and SI-iIOS become the iISS property. For the general case, the three notions obey the implications:
\begin{align*}
    \textrm{SI-iIOS} \quad \Rightarrow \quad \textrm{OL-iIOS} \quad \Rightarrow\quad \textrm{iIOS},
\end{align*}
but none of the implications can be reversed. The output estimate \eqref{yw:e-oliios-1} is more than just a superficial property for mathematical aesthetics. As in the IOS case, the integral-OL property as in \eqref{yw:e-oliios-1} can be related to further properties such as point-wise dissipation rate of Lyapunov functions, and possibly properties related to 
uniform convergence, as demonstrated for finite-dimensional systems in \cite{Ingalls-ACC01}. 

\subsubsection{Lyapunov-like conditions}

The following LKF characterization of iIOS properties is based on the UBEBS assumption, introduced in Definition \ref{def_antoine_UBEBS}, which can be established for instance using the implication \eqref{eq-antoine-UBEBS-Zod}. This result was stated in \cite{Nawarathna20} in the case when the UBEBS property \eqref{eq_antoine_UBEBS} holds with $c=0$, but the result can be proved for the general case when $c>0$ with minor modifications.

\begin{theorem}[Lyapunov sufficiency conditions for iIOS, OL-iIOS, SI-iIOS]
\label{yw:thm-iios}
Assume that (\ref{yw:syso}) owns the UBEBS property. Then the following statements hold:
\begin{enumerate}
\item[$i)$] (\ref{yw:syso}) is iIOS if there exist a functional $V:\mathcal X^n\rightarrow\R_{\ge 0}$, Lipschitz on bounded sets, $\underline\alpha,\overline\alpha\in\cK_\infty$, $\gamma\in\mathcal N$ and $\sigma\in\cKL$ such that, for all $\phi\in\mathcal X^n$ and all $v\in\mR^m$,
\begin{align}
\underline\alpha(|h(\phi(0)|)&\le V(\phi)\le \overline\alpha(\|\phi\|), \label{yw:iios-lyap1}\\
D^+V(\phi,f(\phi,v))&\le-\sigma(V(\phi), \|\phi\|) + \gamma(|v|) \label{yw:iios-lyap2}
\end{align}
\item[$ii)$] (\ref{yw:syso}) is OL-iIOS if there exist a functional $V:\mathcal X^n\rightarrow\R_{\ge 0}$, Lipschitz on bounded sets, $\underline\alpha,\overline\alpha\in\cK_\infty$, $\gamma\in\mathcal N$ and $\sigma\in\cKL$ such that, for all $\phi\in\mathcal X^n$ and all $v\in\mR^m$, (\ref{yw:iios-lyap2}) holds and
\begin{equation}\label{yw:iios-lyap1n}
\underline\alpha(|h(\phi(0)|)\le V(\phi)\le \overline\alpha\left(\|H(\phi)\|\right)
\end{equation}
\item[$iii)$] (\ref{yw:syso}) is SI-iIOS if there is a functional $V:\mathcal X^n\rightarrow\R_{\ge 0}$, Lipschitz on bounded sets, $\underline\alpha,\overline\alpha\in\cK_\infty$, $\gamma\in\mathcal N$ and $\alpha\in\mathcal P$ such that (\ref{yw:iios-lyap1n}) holds and
\begin{equation}\label{yw:iios-lyap2n}
D^+V(\phi,f(\phi,v))\le-\alpha(V(\phi)) + \gamma(|v|).
\end{equation}
\end{enumerate}
\end{theorem}

A special feature about the considered functional for OL-iIOS is that, in (\ref{yw:iios-lyap1n}), $V$ is upper-bounded by a function of supremum norm of the output history. For this property, a point-wise (and even a $\cKL$) dissipation is actually sufficient, as stated by the following result taken from \cite{Nawarathna21}.

\begin{proposition}[OL-iIOS under $\cKL$ dissipation]\label{yw:prop-lyap}
Assume that (\ref{yw:syso}) satisfies the UBEBS property and that there exist a functional $V:\mathcal X^n\to\mRp$ which is  Lipschitz on bounded sets,  $\underline\alpha,\overline\alpha\in\cK_\infty$, $\sigma\in\mathcal{KL}$, and $\gamma\in\mathcal N$ such that, for all $\phi\in\mathcal X^n$ and all $v\in\mR^n$, \eqref{yw:iios-lyap1n} holds, and
\begin{align}
%\underline\alpha(|h(\phi(0)|)&\le V(\phi)\le \overline\alpha(\|H(\phi)\|)\nonumber\\
D^+V(\phi,f(\phi,v))&\le -\sigma\left(|h(\phi(0))|, \|\phi\|\right) + \gamma(|v|).\label{yw:iios15-cdc21}
\end{align}
Then the system is OL-iIOS. 
\end{proposition}

Note that, by \cite[Lemma IV.1]{ANG-bigIISS}, \eqref{yw:iios15-cdc21} is equivalent to the following:
\begin{equation}\label{yw:iios-lyap2nn}
D^+V(\phi,f(\phi,v))\le -\frac{\alpha(|h(\phi(0))|)}{1+\kappa(\|\phi\|)} + \gamma(|v|),
\end{equation}
where $\alpha\in\cK$ and $\kappa, \gamma\in\cN$. The proof of Proposition \ref{yw:prop-lyap} depends critically on the assumption that $V$ is upper-bounded in terms of the output history. This is the main road block to extend Proposition \ref{yw:prop-lyap} to the iIOS case. Nevertheless, in case the considered LKF can be sandwiched between functions involving a particular output, a point-wise (and even a $\cKL$ dissipation) is sufficient to establish iIOS, as stated next.

\begin{corollary}[iIOS under $\cKL$ dissipation]\label{yw:cor-iios}
Assume that (\ref{yw:syso}) owns the UBEBS property and that there exists a locally Lipschitz 
map $h_1: \R^n\to\R_{\ge 0}$ with $h_1(0)=0$ and $\pi\in\mathcal K$ such that
\begin{equation}\label{yw:e-h1h0}
|h(\phi)|\leq \pi(h_1(|\phi(0)|)), \quad \forall \phi\in\mathcal X^n.
%|h(z)|\le \pi(h_1(z)),\quad\forall\,z\in\R^n.
\end{equation}
Suppose that there exist a functional $V:\mathcal X^n\rightarrow\R_{\ge 0}$ which is Lipschitz on bounded sets,  $\underline\alpha,\overline\alpha\in\cK_\infty$, $\sigma\in\mathcal{\cKL}$, and $\gamma\in\mathcal N$ such that, for all $\phi\in\mathcal X^n$ and all $v\in\R^m$,
\begin{align*}
\underline\alpha(h_1(\phi(0)))&\le V(\phi)\le \overline\alpha\left(\|H_1(\phi)\|\right),\\
D^+V(\phi,f(\phi,v))&\le -\sigma(h_1(\phi(0)), \|\phi\|) + \gamma(|v|),
\end{align*}
where $H_1:\cX^n\rightarrow\cX$ is defined by
$H_1(\phi)(\tau)= (h_1\circ\phi)(\tau)$ for all $\tau\in[-\Delta,0]$.
Then the system (\ref{yw:syso}) is iIOS. 
\end{corollary}

It can be seen that under the conditions of Corollary \ref{yw:cor-iios}, the system is OL-iIOS with $h_1$ as output. By property \eqref{yw:e-h1h0}, it follows that the system is iIOS with $h$ as output. %\AC{I'm not sure what you mean here: OL-iIOS always implies IOS (even without \eqref{yw:e-h1h0}). I don't see the added value of Corollary \ref{yw:cor-iios} with respect to Proposition \ref{yw:prop-lyap}.}
%\YW{Correct, OL-iIOS always implies iIOS. However, OL-iIOS with $h_1$ as output doesn't imply iIOS with $h$ as output. The point is that
%the pointwise dissipation rate does not apply to iIOS.  So, the point was that if $V$ is bounded above and below by the output map $h_1$ and $H_1$ (not merely bounded above by $\|\phi\|$) that dominates the original output map $h$, then the pointwise dissipation rate in terms of the new output $h_1$ will lead to iIOS. The motivation was that it can be hard to find a Lyapunov-dissipation rate or a dissipation rate in terms of $\|phi\|$, but a pointwise dissipation in terms of a "bigger" output map will do for iIOS.
%}

\section{Systems interconnection}\label{sec:interconnection}

\subsection{Feedback interconnection}\label{sec:ISS:smallgain}

A possible strategy to analyze ISS of a complex system is to decompose it into ISS subsystems and make sure that their interconnection does not compromise ISS of the overall system. This can be achieved by invoking small gain arguments. The first result in that direction was proposed for finite-dimensional systems in \cite{JIATEEPRA}. It has then been extended to allow for Lyapunov-based conditions \cite{JIMAWA96}, possibly non-ISS subsystems \cite{ItoTAC06,ITJI09}, and interconnections involving more than two subsystems \cite{DARUWI09,ITJIDARU13}. Such a small-gain approach can also be adopted for time-delay systems. To that aim, consider the feedback interconnection of two subsystems with output maps:
\begin{subequations}\label{yw:eq-sg-2}
\begin{align}
\dot x_1(t) =\,&f_1(x_{1t}, w_{1t}, u(t)),\quad
y_1(t) = h_1(x_{1t})\label{yw:eq-sg-21}\\
\dot x_2(t) =\,&f_2(x_{2t}, w_{2t}, u(t)),\quad
y_2(t) = h_2(x_{2t}),\label{yw:eq-sg-22}
\end{align}
\end{subequations}
where, for each $i\in\{1, 2\}$, $x_i(t)\in\R^{n_i}$ and $x_{it}\in \cX^{n_i}$ respectively represent the current value of the solution and the state histories of the $x_i$-subsystem, $y_i(t)\in\R^{p_i}$ represents its output, $w_i\in\cU^{n_{3-i}}$ represents its feedback input (it will be taken as $w_i=y_{3-i}$), and $u\in\cU^m$ represents the exogenous input to the overall system. The maps
$f_i:\cX^{n_i}\times\cX^{n_{3-i}}\times\R^m\rightarrow\cX^{n_i}$ and $h_i:\cX^{n_i}\rightarrow\R^{p_i}$ are assumed to be Lipschitz on bounded sets with $f_i(0, 0, 0)=0$ and $h_i(0) = 0$. The following result is based on a small-gain theorem developed in \cite{Karafyllis:2007kz} for a much wider class of systems, including abstract systems satisfying a semigroup property.

\begin{theorem}[Solutions-based small gain for IOS]\label{yw: sg-thm1}
Assume that, for each $i\in\{1, 2\}$, the $x_i$-subsystem is RFC (with both $w_i$ and $u$ seen as inputs) and there exist $\beta_i\in\cKL$ and $\gamma_i, \mu_i\in\mathcal N$ such that, for all $x_{0i}\in\mathcal X^{n_i}$, all $w_i\in\mathcal U^{n_{3-i}}$ and all $u\in\mathcal U^m$, the following IOS estimate holds:
\[
|y_i(t, x_{0i}, (w_i, u))|\le \beta(\|x_{0i}\|, \ t) + \gamma_i(\|w_{i[0,t]}\|) + \mu_i(\|u_{[0,t]}\|), \quad \forall\, t\ge 0.
\]
Assume that there exists $\rho\in\cK_\infty$ such that, with $\tilde \gamma_i(s) := \gamma_i(s) + \rho(\gamma_i(s))$, it holds that
\begin{align}\label{eq_antoine_112}
\tilde \gamma_1\circ \tilde \gamma_2(s)\le s, \quad\forall\,s\ge 0.
\end{align}
Then the system \eqref{yw:eq-sg-2} under the interconnection feedback
\[
w_1(t) = y_1(t), \quad w_2(t)=y_2(t)
\]
is RFC and IOS with $(x_1, x_2)$ as the state, $(y_1, y_2)$ as the output, and $u$ as the input.
\end{theorem}

The above result therefore provides a small-gain condition \eqref{eq_antoine_112} under which the feedback interconnection of two IOS subsystems is itself IOS. Notice that this small-gain conditions involves only $\gamma_i$, $i\in\{1,2\}$, meaning the internal IOS gains.

\begin{remark}[Solutions-based small gain for ISS]
In the special case when $h_i(\phi_i)=\phi_i(0)$ for each $i\in\{1, 2\}$, the RFC assumption of the individual subsystems becomes redundant and the above small-gain theorem concludes ISS for the overall feedback system.
\end{remark}

In applications, the requirement on ISS (or IOS) can be too stringent.  The next result from \cite{Ito:2010hm} provides a relaxation on the ISS-gain functions by mixing ISS and iISS gains.  The result is here reported in a simplified version for ease of exposition, but can be found in full generality in \cite[Theorem 8]{Ito:2010hm}.

\begin{theorem}[LKF-based small-gain for iISS and ISS]\label{ppyw:sg-thm2}
Suppose that, for each $i\in\{1, 2\}$, there exist a functional $V_i: \cX^{n_i}\rightarrow\R_{\ge 0}$ which is Lipschitz on bounded sets,  $\underline{\alpha}_i, \overline{\alpha}_i\in\cK_\infty$, $\alpha_i\in\cK$ and $\gamma_i, \mu_i\in\cN$ such that, for all $\phi_i\in\mathcal X^{n_i}$, all $\omega_i\in\mathcal X^{n_{3-i}}$ and all $v\in\mR^m$
\begin{align}
    \underline{\alpha}_i(\vert \phi_i(0)\vert)&\le V_i(\phi_i)\le  \overline{\alpha}_i(\|\phi_i\|) \label{yw:eq-sg-3o}\\
   D^+V_i(\phi_i, f_i(\phi_i, \omega_i, v))&\le -\alpha_i(V_i(\phi_i)) + \gamma_{i}(V_{3-i}(\omega_i))
+ \mu_i(|v|). \nonumber 
\end{align}
Assume further that the following three conditions hold:
\begin{itemize}
\item[$i)$] $\alpha_1\in\cK_\infty$ or $\lim_{s\rightarrow+\infty}\gamma_2(s)<+\infty$
\item[$ii)$] $\alpha_2\in\cK_\infty$ or $\lim_{s\rightarrow+\infty}\gamma_1(s)<+\infty$
%$\displaystyle{\lim_{s\rightarrow\infty} \alpha_2(s)=\infty \ \vee\ 
%\lim_{s\rightarrow\infty}\gamma_1(s)<\infty}$, and
\item [$iii)$] $\alpha_2\in\cK_\infty$ or $\displaystyle{\lim_{s\rightarrow+\infty}\alpha_2(s)>\lim_{s\to+\infty} \gamma_2(s)}$.
%$\displaystyle{\lim_{s\rightarrow\infty}\alpha_2(s) = \infty \ \vee\
%\lim_{s\rightarrow\infty}\alpha_2(s) > \lim_{s\rightarrow\infty} \gamma_2(s)}$;
\end{itemize}
Then, under the small-gain condition that there exist $c_1, c_2>1$ such that
\[
(c_1\gamma_1)\circ\alpha_2^{-1}
\circ
(c_2\gamma_2)(s)\le\alpha_1(s), \quad\forall\,s\ge 0.
\]
the system \eqref{yw:eq-sg-2} under the feedback interconnection
\[
w_1(t)=x_2(t), \quad w_2(t)=x_1(t)
\]
is iISS with $x=(x_1, x_2)$ as output and $u$ as input.  Furthermore, if $\alpha_1, \alpha_2\in\cK_\infty$, then the system is ISS. 
\end{theorem}

In the process of developing Theorem \ref{ppyw:sg-thm2}, an iISS LKF with LKF-wise dissipation was actually constructed for the overall system, based on the knowledge of $V_1$ and $V_2$.

In \cite{POLILU06}, an IOS property for a system with output was defined in terms of the UGOS and the OAG properties (see Section \ref{yw:sec-ios-def}), and the following small-gain theorem was developed for the interconnected system \eqref{yw:eq-sg-2} subject to constraints on inputs $(w_1, w_2)$ with $w_1(t)=w_2(t)=0$ for $t<0$ and
\begin{align}\label{yw:eq-sg-101}
|w_1(t)|\le \chi_1(|y_1(t)|), \quad
|w_2(t)|\le \chi_2(|y_2(t)|),
\end{align}
for almost all $t\ge 0$.

\begin{theorem}[Solutions-based small-gain for OAG]\label{yw:sg-thm3}
Consider a forward complete system as in \eqref{yw:eq-sg-2}. Suppose that, for each $i\in\{1, 2\}$, there exist $\sigma_i, \rho_i, \gamma\in\cK$ such that, for all $x_{0i}\in\mathcal X^{n_i}$, all $w_i\in\mathcal U^{n_{3-i}}$, and all $u\in\mathcal U^m$,
\begin{itemize}
\item (UGOS): $|y_i(t,x_{0i},(w_i,u))| \le \max\{\sigma_i(\|x_{0i}\|), \ \rho_i(\|w_{i[0,t]}\|), \ \gamma(\|u_{[0,t]}\|)\}$
\item (OAG):
$\displaystyle{\limsup_{t\to+\infty}|y_i(t,x_{0i},(w_i,u))|\le \max\left\{
\rho_i\left(\limsup_{t\to+\infty}|w_i(t)|\right), \
\gamma\left(\limsup_{t\to+\infty}|u(t)|\right)\right\}.
}
$
\end{itemize}
Then under the small-gain condition $\chi_1\circ\rho_1\circ\chi_2\circ\rho_2(s) < s$ for all $s>0$, the interconnected system \eqref{yw:eq-sg-2} subject to the constraint \eqref{yw:eq-sg-101} satisfies the UGOS and the OAG properties as defined in Section~\ref{yw:sec-ios-def} with $y=(y_1, y_2)$ as output and $u$ as input.
\end{theorem}

These results constitute just a small sample of the extensive work that has been conducted in the past two decades on small-gain theorems for time-delay systems. Without claiming to be exhaustive, we now list some of them. In \cite{Enciso2006}, a small-gain theorem was developed for monotone time-delay systems. In \cite{ItoPepeJiang:Automatica2013}, a Lyapunov-based small-gain theorem was presented for iISS. In \cite{Karafyllis:2011fk}, an ISS small-gain theorem was developed based on vector Lyapunov functions for large-scale systems, meaning involving possibly more than two subsystems.  In \cite{Dashkovskiy-MTNS09},
an ISS small-gain theorem was obtained for large-scale systems based on both the Lyapunov-Krasovskii and the Lyapunov-Razumikhin approaches and was extented to reset systems in \cite{dashkovskiy2012stability}. In \cite{TIWAJIlargescale12} both Lyapunov-based and solutions-based small-gain results were developed for IOS. Notably, extensive work on small-gain theorems has been developed for abstract systems and
infinite-dimensional systems that can be applied to time-delay systems, see for instance \cite{SontagIngalls02,LiuFren-ECC13,Mironchenko:2021,Mironchenko:SIAM2021,MirKawGlu-2021}. 

\subsection{Cascade interconnection}

A cascade pertains to the unidirectional interconnection of two subsystems, namely:
\begin{subequations}\label{eq_antoine_10}
\begin{align}
		\dot x(t)=\,&f_1(x_{t},z_t,u(t))
		\label{eq_antoine_10a}\\
		\dot z(t)=\,&f_2(z_t,u(t)).
		\label{eq_antoine_10b}
\end{align}
\end{subequations}
$x(t)\in\mathbb{R}^{n_1}$ and $z(t)\in\mathbb{R}^{n_2}$ denote respectively the current values of the states of the driven and driving systems, $x_{t}\in\mathcal X^{n_1}$ and $z_{t}\in\mathcal X^{n_2}$ are the corresponding state histories, and $u\in\mathcal U^m$ denotes the exogenous input. The functions $f_1:\mathcal{X}^{n_1}\times\mathcal{X}^{n_2}\times \mR^m\to\mathbb{R}^{n_1}$ and $f_2:\mathcal{X}^{n_2}\times \mathbb R^m\to\mathbb{R}^{n_2}$ are assumed to be Lipschitz on bounded sets and to satisfy $f_1(0,0,0)=0$ and $f_2(0,0)=0$.

In this setting, the driving subsystem \eqref{eq_antoine_10b} influences the dynamics of the driven one \eqref{eq_antoine_10a}, but the latter has no influence on the former. Note that a feed through term is allowed in \eqref{eq_antoine_10}, meaning that the exogenous input may affect both the driving and the driven subsystems. 

In finite dimension, a vast literature has been devoted to the preservation of stability under cascade interconnections: see for instance \cite{SEISUA,CASCSCL,CASCAUT,MURAT-CASC,CHANscl05,CHANIT14}. In particular, for delay-free systems, the ISS property is known to be preserved under cascade interconnection. This has been established in \cite{SONTEETAC95} using a change of supply rate technique. An alternative proof strategy consists in seeing  a cascade interconnection as a particular type of feedback interconnection and to invoke small-gain theorems in which an internal input gain is set to zero. This was the path followed in \cite{Karafyllis:2007kz} and later on in \cite[Remark 13]{Ito:2010hm}, which yields the following corollary of Theorem \ref{yw: sg-thm1}.

\begin{corollary}[ISS preservation under cascade]
\label{cor-cascade}
Assume that \eqref{eq_antoine_10a} is ISS with respect to its inputs $z$ and $u$ and that \eqref{eq_antoine_10b} is ISS. Then the cascade \eqref{eq_antoine_10} is ISS.
\end{corollary}

For cascades without exogenous inputs (meaning $u\equiv 0$ in \eqref{eq_antoine_10}), an immediate consequence of this result is that, if the driven subsystem is ISS with respect to its input $z$ and the driven subsystem is GAS, then the whole cascade is GAS. We illustrate the application of this result through the following example taken from \cite{Adimy-2006,Adimy-2008}.

\begin{example}[Blood cells dynamics]\label{yw:ex-cascade}
A simplified model describing the production process of blood cells is
\begin{align}
\dot x_i(t) =&\, -\left(\delta_i + q_i(x_i(t))\right)x_i(t)
+ 2(1-k_i)\int_{-\Delta_i}^0\alpha_i(\tau)q_i(x_i(t+\tau))x_i(t+\tau)d\tau\nonumber\\
&\, + 2k_{i-1}\int^0_{-\Delta_{i-1}}\alpha_{i-1}(\tau)q_{i-1}(x_{i-1}(t+\tau))x_{i-1}(t+\tau)d\tau, \label{yw:eq-ex1}
\end{align}
where, for each $i\in\{1,\ldots,n\}$, $x_i$ denotes the density of hematopoietic stem cells in the $i^{\textrm{th}}$ generation. The first term in the right-hand side of \eqref{yw:eq-ex1} describes the loss of resting cells of the $i^{\textrm{th}}$ generation either by death or introduction in the proliferating phase, the second term accounts for  self-renewing cells re-entering the generation, and the last term models differentiation cells of the previous generation entering the $i^{\textrm{th}}$ generation. The following assumptions can be made for each $i\in\{1,\ldots,n\}$:
\begin{itemize}
    \item $\alpha_i:[-\Delta,0]\to\mRp$ is integrable over $[-\Delta_i, 0]$
    \item $q_i:\R_{\ge 0}\rightarrow\R_{\ge 0}$ is locally Lipschitz, decreasing, and $q_i(s)\rightarrow 0$ as $s\rightarrow +\infty$ 
    \item the constants $\delta_i, k_i$ are positive, with $k_i\le 1$ and $k_0=0$.
\end{itemize}
For convenience, let $q_i(s) = q_i(-s)$ for $s<0$ so that $q(s) = q(|s|)$.  It can be seen that the trajectories of the $x_i$-subsystem satisfy the following estimate:
\begin{align}
|x_i(t)|\le |x_i(0)|e^{-\delta_it} & + \frac{2(1-k_i)\int^0_{-\Delta_i}\alpha_i(\tau)d\tau}{\delta_i} q_i(0)\max_{\tau\in[t-\Delta_i, t]} |x_i(\tau)|
\label{yw:eq-ex0}\\
& + \frac{2(1-k_{i-1})\int^0_{-\Delta_{i-1}} \alpha_{i-1}(\tau)d\tau}{\delta_i} q_{i-1}(0)\max_{\tau\in [t-\Delta_{i-1}, t]}|x_{i-1}(\tau)|\nonumber
\end{align}
on the maximum interval of existence of the solutions.
By the Razumikhin approach for ISS (Theorem \ref{thm-iss-raz}), the $x_i$-subsystem is ISS with $x_{i-1}$ as input if
\begin{equation}\label{yw:eq-ex2}
2(1-k_i) q_i(0)\int^0_{-\Delta_i}\alpha_i(\tau)d\tau
< \delta_i.
\end{equation}
By Corollary \ref{cor-cascade}, we conclude that the overall cascade is GAS if \eqref{yw:eq-ex2} holds for all $i\in\{1,\ldots,n\}$. 
\end{example}

On some occasions, ISS may happen to be a too demanding requirement for cascade analysis. In finite dimension, the weaker property of Strong iISS (which amounts to the combination of ISS with respect to small inputs and iISS \cite{CHANIT13}) was also shown to be preserved under cascade interconnection \cite{CHANIT14}. If both individual subsystems are assumed to be merely iISS, then the cascade is not guaranteed to be iISS (or even 0-GAS) as shown though an example in \cite{MURAT-CASC}. Extra conditions are needed. Such conditions have been proposed in \cite{CHANscl05} for finite-dimensional cascades and have recently been extended to time-delay systems in \cite{CHGO22}, as stated next.

\begin{theorem}[Cascade of iISS systems]\label{thm_antoine_Cascade2}
Let $p\in\mN$ and $\mathcal I\subset\{1,\ldots, ,n_2\}$. Assume there exist two LKF $V_1:\mathcal X^{n_1}\to\mathbb R_{\geq 0}$ and $V_2:\mathcal X^{n_2}\to\mathbb R_{\geq 0}$, $\sigma\in\cKL$, $\gamma\in\cK_\infty$ and, for each $i\in\mathcal I$ and each $j\in\{1,\ldots,p\}$, $\sigma_i\in\mathcal{KL}$, $\gamma_{ij}\in\cK_\infty$ and $\Delta_{ij}\in[0,\Delta]$ such that, for all $\phi\in\mathcal X^{n_1}$, all $\varphi=(\varphi_1,\ldots,\varphi_{n_2})^\top\in \mathcal X^{n_2}$, and all $v\in\mR^m$,
\begin{align}
		D^+V_1(\phi,f_1(\phi,\varphi,v))&\le -\sigma\big(|\phi(0)|,V_1(\phi)\big)
				+\sum_{i\in\mathcal I}\sum_{j=1}^p \gamma_{ij}\big(|\varphi_i(-\Delta_{ij})|\big)+\gamma(|v|)\label{eq_antoine_48b}\\
D^+V_2(\varphi,f_2(\varphi,v))&\le
				-\sum_{i=1}^{n_2} \sigma_i\big(|\varphi_i(0)|,V_2(\varphi)\big)+\gamma(|v|).
		\label{eq_antoine_48d-bis}
\end{align}
Then the cascade (\ref{eq_antoine_10}) is iISS provided that the following growth condition holds:
\begin{align}\label{eq_antoine_growth-bis}
\limsup_{s\rightarrow 0^+} \frac{\gamma_{ij}(s)}{\sigma_i(s,0)}<+\infty,\quad \forall j\in\{1,\ldots,p\},\,i\in\mathcal I.
%\gamma_{ij}(s)=\mathcal{O}_{s\rightarrow 0^+}\big(\sigma_i(s,0)\big),\quad \forall j\in\{1,\ldots,p\},\,i\in\mathcal I.
\end{align}
\end{theorem}

In view of Theorem \ref{theo_antoine_iISS_char}, it can be easily checked that condition \eqref{eq_antoine_48d-bis} is equivalent to requiring that the driving subsystem is iISS (this is formally established in \cite[Proposition 2]{CHGO22}).  Similarly, condition \eqref{eq_antoine_48b} is equivalent to imposing that the driven subsystem \eqref{eq_antoine_10a} is iISS. The extra condition to ensure that the cascade is iISS lies in the growth constraint \eqref{eq_antoine_growth-bis}, which imposes that the $\cKL$ dissipation rate of $V_2$ dominates around zero the input rates $\gamma_{ij}$ through which the state variables $\varphi_{i}$ enter the dissipation inequality \eqref{eq_antoine_48b}.

Note that not all state variables $\varphi_i$ of the driving system \eqref{eq_antoine_10b} may appear in \eqref{eq_antoine_48b}, as allowed by the fact that $\mathcal I$ is just a subset of $\{1,\ldots,n_2\}$. Accordingly, no growth condition is imposed in \eqref{eq_antoine_growth-bis} for the variables $\varphi_i$ that are not involved in the dissipation inequality \eqref{eq_antoine_48b}.

\section{Open questions}\label{sec:open}

The results listed in this survey show that the ISS framework for time-delay systems is now a rather mature subject. Nonetheless, some important questions remain unsolved. We list some of them below.

\subsection{Is GAS necessarily uniform?}\label{sec:open:UGAS}
In this paper, we have decided to define GAS through a $\cKL$ estimate of the solution's norm (see Definition \ref{pp:defstabetc}). This terminology hides a crucial property of GAS: both the transient overshoot and the decay rate of solutions are requested to be uniform over bounded sets of initial conditions. This explains why this property is sometimes referred to as UGAS in the literature (``U'' standing for ``uniform'').

For finite-dimensional systems, this uniformity comes for free. In other words, the $\cKL$ estimate boils down to the combination of Lyapunov stability and convergence of solutions to the origin. Theorem \ref{pp:equivGASsolita} states that the same holds for time-delay systems, provided that the RFC property is satisfied. On the other hand, Theorem \ref{pp:equivASsolita} shows that, as far as local asymptotic stability is concerned, the RFC assumption is not needed. It is not clear yet whether this extends to global properties, thus leading to the following conjecture.

\begin{conj}[GAS $\Leftrightarrow$ UGAS]\label{conj:GAS}
The time-delay system \eqref{pp:RFDE_autonomous} is GAS if and only if the following two conditions hold:
\begin{itemize}
    \item (stability) for every $\varepsilon>0$ there exists $\delta>0$ such that, for any $x_0\in \mathcal X^n$ with $\Vert x_0\Vert\le \delta$, the corresponding solution of \eqref{pp:RFDE_autonomous} satisfies $\vert x(t,x_0)\vert\le \varepsilon$ for all $t\ge 0$
    \item (global attractivity) for any $x_0\in {\cal X}^n$, the
corresponding solution of \eqref{pp:RFDE_autonomous} satisfies
$\lim_{t\to +\infty}x(t,x_0)=0$.
\end{itemize}
\end{conj}

Since global convergence of solutions to the origin trivially ensures forward completeness, a positive answer would be given to the above conjecture if one was able to establish the following.

\begin{conj}[FC $\Leftrightarrow$ RFC]\label{conj_antoine_RFC}
If the time-delay system \eqref{pp:RFDE_autonomous} is forward complete (in the sense that its solutions exist at all positive times: see Definition \ref{def_FC}), then it is RFC in the sense of Definition \ref{def_RFC}.
\end{conj}

 Note that the equivalence between FC and RFC does not hold for general infinite-dimensional systems (see the example in \cite[Section VI]{mironchenko2017characterizations}), but the peculiarities of time-delay systems leave some hope for Conjecture \ref{conj_antoine_RFC} to hold true. The interested reader is referred to \cite{karafyllis2022global} for more in-depth discussion on this matter.

%\AC{Do we want to leave the text below?} A possible avenue to prove Conjecture \ref{conj_antoine_RFC} would be to consider the following property:
% \begin{align}\label{eq_antoine_RFC2}
%     \sup\{\|x_t(x_0,u)\|\,:\, t\in[0,T],\|u\|\leq R\}<+\infty, \quad \forall x_0\in\mathcal X^n, \quad \forall \,T,R\geq 0.
% \end{align}
%This property means that the reachable set from each given initial state with bounded input in finite time is bounded. Unlike in RFC, no uniformity requirement is imposed with respect to $x_0$. Clearly, the following relationships hold:
%\begin{align*}
%    \textrm{RFC}\quad \Rightarrow\quad \eqref{eq_antoine_RFC2}\quad \Rightarrow\quad \textrm{FC}.
%\end{align*}
%It seems reasonable to show that RFC is actually equivalent to \eqref{eq_antoine_RFC2}. If yes, the remaining question would be to show that FC $\Rightarrow$ \eqref{eq_antoine_RFC2}. This constituted the critical step in \cite{LINSONWAN} that was solved by considering backward solutions (which does not appear as a relevant strategy in a time-delay context).

\subsection{GES under point-wise dissipation}\label{sec:open:GES}

We have seen with Theorem \ref{pp:teoremafondamentalekargas} that GAS is equivalent to the existence of an LKF that dissipates in a point-wise manner along the system's solutions. Theorem \ref{pp:gesteoremapointwise} shows that, under additional quadratic constraints on the LKF, a point-wise dissipation inequality is also sufficient to establish GES. However, this result holds so far only under a growth rate assumption on the vector field, thus significantly reducing its range of applicability. It is not known whether this growth rate condition is indeed needed to conclude GES, thus leading to the following conjecture.

\begin{conj}[GES under point-wise dissipation]\label{conj_antoine_GES}
Let the map $f_0$ in \eqref{pp:RFDE_autonomous} be Lipschitz on bounded sets with $f_0(0)=0$ and assume that there exist an LKF $V:\mathcal {X}^n\to \mathbb R_{\ge 0}$ and $\underline a,\overline a, a>0$ such that, for all $\phi \in \mathcal {X}^n$,
\begin{align*}
\underline a\vert \phi(0)\vert^2\le V(\phi)\le \overline a\Vert \phi\Vert^2 \\
D^+V(\phi,f_0(\phi))\le -a\vert \phi(0) \vert^2.
\end{align*}
Then the system \eqref{pp:RFDE_autonomous} is GES.
\end{conj}

Note that the converse holds from Theorem \ref{pp:gesteoremapointwise}. The validity of Conjecture \ref{conj_antoine_GES} would thus provide a complete characterization of GES through a point-wise dissipation.

\subsection{ISS under point-wise dissipation}\label{sec:open:pointwiseISS}

We have seen with Theorem \ref{theo_antoine_iISS_char} that iISS is guaranteed under the existence of an iISS LKF with point-wise dissipation. As discussed in \cite{CHPEMACH17}, whether or not this extends to the ISS property remains an open question. Theorem \ref{theo_antoine_ISS_pointwise2} provides a preliminary answer to this question, but imposes growth restrictions on the dissipation rate.

We believe that solving the next conjecture would thus constitute a significant contribution to the ISS theory for time-delay systems, not only to simplify computations in practical applications (a point-wise dissipation being usually easier to obtain than an LKF-wise one) but also for the sake of homogeneity with the input-free stability analysis.

\begin{conj}[ISS under point-wise dissipation]\label{conj:ISSpointwise}
Assume that there exist an LKF $V:\mathcal X^n\to\mRp$, $\alpha\in\cK_\infty$ and $\gamma\in\mathcal N$ such that, for all $\phi\in\mathcal X^n$ and all $v\in\mR^m$,
\begin{align*}
    D^+V(\phi,f(\phi,v))\leq -\alpha(|\phi(0)|)+\gamma(|v|).
\end{align*}
Then the system \eqref{pp:RFDE} is ISS.
\end{conj}

It is worth noting that the converse of this result straightforwardly holds from Theorem \ref{theo_antoine_ISS_char}. If valid, this conjecture would thus constitute an alternative characterization of ISS.

It is possible to show that, under the assumptions of Conjecture \ref{conj:ISSpointwise}, the following property holds: there exists $\mu\in\mathcal N$ such that, for every $\varepsilon, r>0$ there exists $T >0$ such that, for all $x_0 \in \mathcal X^n$ and all $u\in \mathcal U^m$ with $\left\| \phi \right\| \le r$,
\begin{align}\label{eq_antoine_20}
\min_{t\in[0,T]} \left| x(t,x_0 ,u)\right|\le \varepsilon +\mu \left(\|u\|\right).
\end{align}
This property shows that any solution \emph{point-wisely} visits a $\mu(\|u\|)$-neighborhood of the origin and that the time needed for such a visit is uniform on bounded sets of initial states and inputs. However, unlike the ULIM property employed in \eqref{eq_antoine_19}, it does not guarantee that the whole history segment visits such a neighborhood. Since, in view of Theorem \ref{pp:teoremafondamentalekargas}, the assumptions of Conjecture \ref{conj:ISSpointwise} clearly ensure 0-GAS, the latter would be solved if one could establish the following.

\begin{conj}[0-GAS + point-wise ULIM $\Rightarrow$ ISS]\label{conj-ULIM}
Assume that the system \eqref{pp:RFDE} is 0-GAS and satisfies the above point-wise version of the ULIM property. Then it is ISS.
\end{conj}

An alternative strategy to establish Conjecture \ref{conj:ISSpointwise} is the following. Integrating both sides of \eqref{eq_antoine_20} along the system's solutions, one easily get that, for some $\overline\alpha\in\cK_\infty$,
\begin{equation} \label{eq_antoine_21}
\int_{0}^{t}\alpha\left(\left| x(\tau,x_0 ,u)\right| \right)d\tau \le \overline\alpha \left(\left\| x_0 \right\| \right)+\int _{0}^{t}\gamma \left(\left|u(\tau)\right|\right)d\tau   ,\quad \forall t\geq 0.
\end{equation}
This constitutes a point-wise version of the integral-to-integral estimate \eqref{Iasson__25_}. Thus, Conjecture \ref{conj:ISSpointwise} would be proved if one was able to establish the following.

\begin{conj}[Point-wise integral-to-integral estimate $\Rightarrow$ ISS]\label{conj-int-pointwise}
Assume there exist $\alpha,\overline\alpha\in\cK_\infty$ and $\gamma\in\mathcal N$ such that \eqref{eq_antoine_21} holds for all $x_0\in\mathcal X^n$ and all $u\in\mathcal U^m$. Then the system \eqref{pp:RFDE} is ISS.
\end{conj}

\subsection{Solutions-based characterization of ISS}

A vast catalog of properties equivalent to ISS have been proposed in \cite{SONWANTAC} for finite-dimensional systems. Some of them were already extended to time-delay systems, as summarized by Theorem \ref{theo-sol-char-ISS}. Yet, important equivalences valid in finite dimensions are still missing for time-delay systems.

Some of these characterizations have been proposed in \cite[Theorem 3.4]{Mironchenko:2017wb}, but are only valid under the RFC assumption: a positive answer to Conjecture \ref{conj_antoine_RFC} would thus have important consequences on solutions-based characterizations of ISS. Yet, the results in \cite{Mironchenko:2017wb} are specializations of other ISS characterizations developed in a wider infinite-dimensional context \cite{mironchenko2017characterizations}. In particular, most notions employ the sup norm of the state history, thus paying little attention to point-wise properties, such as the ones used in Conjectures \ref{conj-ULIM} and \ref{conj-int-pointwise}. 

We believe that further work is needed in those directions to get a full view of the solutions-based characterizations of the ISS property.

\subsection{Converse LKF with point-wise dissipation for iISS}

While Theorem \ref{theo_antoine_iISS_char} shows that iISS holds under the existence of an iISS LKF with point-wise dissipation, it is not yet known whether the converse also holds true.

\begin{conj}[iISS $\Rightarrow$ iISS LKF with point-wise dissipation]
If the system \eqref{pp:RFDE} is iISS then it admits an iISS LKF with point-wise dissipation.
\end{conj}

While Theorem \ref{theo_antoine_iISS_char} ensures  the (seemingly more demanding) existence of a coercive iISS LKF with history-wise dissipation rate, it is not clear whether this is enough to ensure a point-wise dissipation due to the $\mathcal P$ nature of the considered dissipation rate.

\subsection{Solutions-based characterizations of iISS}\label{sec-open-sol-iISS}

We have seen in Theorem \ref{theo-iISS2} that iISS is equivalent to the combination of 0-GAS and the UBEBS property. The latter imposes a uniform bound on the solutions' norm in terms of the norm of the initial state and the energy brought by the input.

For delay-free systems, an even weaker requirement was shown to ensure iISS when combined to 0-GAS \cite{ANGING04}. This requirement is known as the \emph{bounded input -  frequently bounded state} property (BEFBS) and imposes, as its name suggests, that under a bounded energy assumption on the input, the $\liminf$ of the solutions' norm is finite. This result probably constitutes the weakest solutions-based requirement to establish iISS for finite-dimensional systems, but its validity in a time-delay context is not known yet.

\begin{conj}[0-GAS + BEFBS $\Rightarrow$ iISS]\label{conj:BEFBS}
Assume that \eqref{pp:RFDE} is 0-GAS and satisfies the following \emph{bounded input -  frequently bounded state} property\footnote{This implication implicitly imposes that the system is forward complete}:
\begin{align*}
    \int_0^{+\infty} \zeta(|u(\tau)|)d\tau<+\infty\quad\Rightarrow\quad \liminf_{t\to+\infty} |x(t,x_0,u)|<+\infty
\end{align*}
for some $\zeta\in\cK_\infty$ and for any $x_0\in\mathcal X^n$ and any $u\in\mathcal U^m$. Then it is iISS.
\end{conj}

The proof of this result for delay-free systems, proposed in \cite{ANGING04}, relies on the notion of \emph{input/output-to-state stability} (IOSS), which essentially quantifies zero-state detectability in the presence of exogenous inputs \cite{KRISONWAN} and whose theory is still insufficiently developed for time-delay systems (see Section \ref{sec-conj-IOSS}).

Note that the converse of Conjecture \ref{conj:BEFBS} trivially holds from the definition of iISS. 

\subsection{Strong iISS}

For finite-dimensional systems, a property was introduced half-way between ISS and iISS. This property, called \emph{Strong iISS}, imposes that the system is both iISS and ISS with respect to inputs with small amplitude \cite{CHANIT13}. For bilinear systems, Strong iISS happens to be equivalent to 0-GAS and some Lyapunov-based conditions have been proposed to guarantee it in practice. Similarly to ISS, and contrarily to iISS, Strong iISS happens to be naturally preserved under cascade interconnection \cite{CHANIT14}. This property is also at the core of several small-gain results for non-ISS systems \cite{ItoTAC06,ITJI09,ITJIDARU13}, it has been studied for Lur'e systems using the circle criterion \cite{guiver2020circle} and has been shown to be achievable through saturated feedback for some classes of systems \cite{AZCHCHGR15}.

To date, the extension of this notion to time-delay systems remains rather limited, which explains why it is not covered by this survey. We believe that developing theory for the Strong iISS property in the context of time-delay systems constitutes a relevant direction of research.

\subsection{IOS without LKF-wise dissipation}

We have seen in Section \ref{yw:sec-autonomous} that a functional that dissipates only in terms of the output norm is not enough to conclude GAOS, as such a dissipation does not ensure a uniform convergence rate over bounded sets of initial states. In Theorem \ref{yw:thm2-gaos}, we have also seen that an extra condition on the dissipation rate can be imposed (namely, that it does not increase along the system's solutions) to conclude GAOS. So far, this approach is confined to input-free systems.

An interesting line of research would be to investigate what additional condition can be imposed to conclude IOS from a functional whose derivative does not dissipate in terms itself (as in Theorem \ref{theo-LKF-IOS}), but rather in terms of an output of the system, without requesting that the functional is sandwiched between $\cK_\infty$ functions of the output norm (as imposed by Corollary \ref{yw:cor-iios}).

\subsection{Solutions-based characterizations of IOS}\label{sec-conj-sol-IOS}

As we have seen in Section \ref{yw:sec-ios-def}, the example provided in \cite{ORCHSI-lcss20} shows that, even for finite-dimensional systems,
\begin{align*}
    UGOS \ \& \ OAG \quad \nRightarrow \quad IOS.
\end{align*}
Here again, this is due to the fact that IOS requires some uniformity of the output convergence to zero, which is not captured the OAG property (see \eqref{eq-OAG}). This leads naturally to the question whether, when combined with OAG, the output-Lagrange stability introduced in \eqref{yw:olios}, namely:
\begin{equation}\label{yw:oliosbis}
|y(t, x_0, u)|\le
\max\left\{\sigma(|h(x_0)|), \ \sigma(\|u_{[0, t]}\|)\right\},\quad \forall t\geq 0
\end{equation}
for some $\sigma\in\cK_\infty$, is enough to ensure the OL-IOS property (see Definition \ref{yw:def-ios}) . More precisely, we conjecture the following.

\begin{conj}[Output-Lagrange stability {\rm \&} OAG $\Rightarrow$ OL-IOS]
Let (\ref{yw:syso}) be forward complete and assume that there exist $\sigma\in\cK_\infty$ and $\mu\in\mathcal N$ such that, for all $x_0\in\mathcal X^n$ and all $u\in\mathcal U^m$, its solution satisfies \eqref{yw:oliosbis} and
\begin{align*}
    \lim_{t\to +\infty} |y(t,x_0,u)|\leq \mu(\|u\|).
\end{align*}
Then the system is OL-IOS.
\end{conj}

Such a result would be the natural counterpart of the sufficient condition for OL-iIOS provided in Proposition \ref{yw:prop-lyap}. The key point in the above statement is that \eqref{yw:oliosbis} requests that the output is uniformly bounded in terms of the initial output value (but not the other state variables) and the applied input. This feature did not hold in the counter-example of \cite{ORCHSI-lcss20}.

\subsection{Input/output-to-state stability}\label{sec-conj-IOSS}

For finite dimensional systems, the notion of \emph{output-to-state stability (OSS)} has been introduced in \cite{SONSCLOSS}. This notion considers systems with an output that represents measurements that can be made on the system and thus constitutes the only information available on the state variables (whereas, in IOS, the output should be thought of as the state variables one wishes to stabilize). The OSS property then characterizes ``zero state detectability'', meaning the fact that, whenever the output is identically zero, the state necessarily vanishes. In a time-delay context, this property would consist in the existence of $\beta\in\cKL$ and $\gamma\in\mathcal N$ such that, for all $x_0\in\mathcal X^n$,
\begin{align*}
    |x(t,x_0)|\leq \beta(\|x_0\|,t)+\gamma(\|y_{[0,t]}\|)
\end{align*}
over the maximal interval of existence of the solution of \eqref{yw:sys-auto}. In the presence of an input, we may additionally impose that this property remains valid up to an error which is small whenever the input is sufficiently small. This leads to the notion of \emph{input/output-to-state stability (IOSS)}, which was also introduced for finite-dimensional systems in \cite{SONSCLOSS}. For time-delay systems, the IOSS property would impose that there exist $\beta\in\cKL$ and $\gamma,\mu\in\mathcal N$ such that, for all $x_0\in\mathcal X^n$ and all $u\in\mathcal U^m$, 
\begin{align*}
    |x(t,x_0,u)|\leq \beta(\|x_0\|,t)+\gamma(\|y_{[0,t]}\|)+\mu(\|u_{[0,t]}\|)
\end{align*}
over the maximal interval of existence of the solution of \eqref{yw:syso}. Lyapunov characterizations of both OSS and IOSS have been proposed in \cite{SONSCLOSS,KRISONWAN}, whereas as solutions-based characterizations can be found in \cite{ANGING04}. These properties were also at the basis of the design of state norm estimators \cite{KRISONWAN} and were used to derive the ISS property based on the knowledge of a Lyapunov function that dissipates only in terms of the output norm \cite{ANGPD}. The IOSS property was also fundamental to show that iISS is equivalent to the combination of 0-GAS and BEFBS, as discussed in Section \ref{sec-open-sol-iISS}.

To date, we are not aware of any attempt to extend the OSS and IOSS theory to time-delay systems. 

\section*{Acknowledgement}

The authors would like to warmly thank G\"okhan G\"oksu and Epiphane Loko for their careful reading of this survey.

%\section{Conclusion}

%In this survey, we have tried to summarize in a condensed and homogeneous way the main results existing so far in the ISS framework for time-delay systems. We have also underlined remaining open questions in this field, which would deserve further investigations.

%Although this framework was at the basis of several control designs for time-delay systems, we have decided not to present them here. In particular, for control methodologies that account for delays in the input/output channel, the reader is invited to consult the monographs \cite{BEKIARISKRSTIC2013, KARKRS2017, KRSTICBOOK2009, RICOCAMACHO2007, ZHANGXIE2007, ZHONG2006, ZHUQIMACHEN18}, as well as the papers \cite{mazenc2006backstepping,mazenc2008further} \PP {there are many papers for control with input/output delays in both the linear and the nonlinear case.... we should pay some attention here if we cite just a few ones and not others.} \AC{You are right, it is just that I feel unfair not to cite Mazenc's work\ldots} \PP{I Would add also recent book \cite{LIUFRIXIA2021} and the tutorial on stabilization of delay systems \cite{KARMALMAZPEP2016}}

\bibliographystyle{apalike}
\bibliography{refs}

\end{document}